\documentclass[12pt]{article}
\usepackage[utf8]{inputenc}
\usepackage[T1]{fontenc}
\usepackage[english]{babel}
\usepackage{amssymb,amsmath,amsthm,amsfonts,bm,bbm}
\usepackage{graphicx,graphics, enumerate}
\usepackage{booktabs, caption, multirow, float, url}
\usepackage{dsfont}
\usepackage{comment}
\usepackage[titletoc,title]{appendix}
\usepackage{adjustbox}
\usepackage{color}
\usepackage{indentfirst}
\usepackage[ruled]{algorithm2e}
\usepackage{hyperref}
\usepackage[a4paper, margin=1.25in, nohead]{geometry}

\usepackage{tikz}
\usetikzlibrary{trees,positioning,arrows,arrows.meta}

\usepackage{setspace} 
\newtheorem{hyp}{Assumption}
\newtheorem{prop}{\bf Proposition}[section]

\newtheorem{rem}[prop]{\bf Remark}
\newtheorem{thm}[prop]{\bf Theorem}
\newtheorem{lemma}[prop]{\bf Lemma}

\newtheorem{example}[prop]{\bf Example}

\def\simiid{\buildrel {\text{i.i.d.}} \over \sim}
\newcommand{\convD}{\underset{n\to +\infty}{\overset{d}{\longrightarrow}}}
\newcommand{\equalD}{\stackrel{d}{=}}
\newcommand{\convP}{\underset{n\to +\infty}{\overset{P}{\longrightarrow}}}
\newcommand{\convAS}{\underset{n\to +\infty}{\overset{a.s.}{\longrightarrow}}}

\newcommand{\eps}{\varepsilon}
\newcommand{\epstilde}{\tilde{\varepsilon}}

\newcommand{\meanU}{\overline{U}_n}
\newcommand{\meanX}{\overline{X}_n}
\newcommand{\meanY}{\overline{Y}_n}

\newcommand{\Prob}{\mathbbm{P}}
\newcommand{\Expec}{\mathbbm{E}}
\newcommand{\Var}{\mathbbm{V}}
\newcommand{\expecX}{\Expec[X]}
\newcommand{\expecY}{\Expec[Y]}

\newcommand{\expecXn}{\Expec[X_{1,n}]}
\newcommand{\expecYn}{\Expec[Y_{1,n}]}

\newcommand{\expecXnsq}{\Expec[X_{1,n}^2]}
\newcommand{\expecYnsq}{\Expec[Y_{1,n}^2]}

\newcommand{\VarX}{\Var[X]}
\newcommand{\VarY}{\Var[Y]}
\newcommand{\VarYn}{\Var[Y_{1,n}]}

\newcommand{\Cov}{\mathbbm{C}ov}
\newcommand{\Corr}{\mathbbm{C}orr}

\newcommand{\upperYn}{u_{Y,n}}
\newcommand{\upperXn}{u_{X,n}}
\newcommand{\lowerYn}{l_{Y,n}}
\newcommand{\aXn}{a_{X,n}}
\newcommand{\bXn}{b_{X,n}}
\newcommand{\aYn}{a_{Y,n}}
\newcommand{\bYn}{b_{Y,n}}

\newcommand{\Nstar}{\mathbbm{N}^{*}}

\newcommand{\Indicator}{\mathds{1}}

\newcommand{\distributionXY}{P_{X,Y}}
\newcommand{\distributionXYn}{P_{X,Y,n}}
\newcommand{\distributionX}{P_X}
\newcommand{\distributionY}{P_Y}
\newcommand{\Nc}{\mathcal{N}}
\newcommand{\Uc}{\mathcal{U}}
\newcommand{\Rb}{\mathbb{R}}

\newcommand{\Rstarplus}{\mathbb{R}_+^*}
\newcommand{\Rstar}{\mathbb{R}^*}
\newcommand{\INT}{\mathbb{N}}
\newcommand{\INTST}{\mathbb{N}^*}
\newcommand{\lsin}{\left\{}
\newcommand{\rsin}{\right\}}

\newcommand{\Bz}{\mathcal{B}\left(\mathcal{Z}\right)}
\newcommand{\Zc}{\mathcal{Z}}
\newcommand{\Pfam}{\mathcal{P}}

\newcommand{\Csetstheta}{\mathcal{F}_{\Theta}}
\newcommand{\BCsetstheta}{\mathcal{B}\left(\mathcal{F}_{\Theta}\right)}
\newcommand{\eventAXeps}{A_{\eps}^{X_{1,n}}}
\newcommand{\eventAYepstilde}{\tilde{A}_{\epstilde}^{Y_{1,n}}}
\newcommand{\SNRn}{\text{SNR}_n}
\newcommand{\SNRntilde}{\widetilde{\text{SNR}}_n}
\newcommand{\thetahatn}{\widehat{\theta}_n}
\newcommand{\Cnalpha}{C_{n,\alpha}}

\usepackage[colorinlistoftodos]{todonotes}

\title{On the construction of confidence intervals for ratios of expectations\thanks{
We would like to thank Laurent Davezies, Xavier D'Haultf\oe{}uille,  and the participants of the CREST internal seminar (Nov. 2018) for their valuable comments. This research has been supported by the Labex Ecodec.}}

\author{Alexis Derumigny\thanks{CREST,
5, avenue Henry Le Chatelier,
91764 Palaiseau cedex, France.
\newline
E-mail adresses: firstname.lastname@ensae.fr for the three authors.},
Lucas Girard\footnotemark[2], Yannick Guyonvarch\footnotemark[2]
}

\date{\today}

\begin{document}

\begin{titlepage}
\clearpage
\maketitle
\thispagestyle{empty}



\begin{abstract}
    \noindent
    In econometrics, many parameters of interest can be written as ratios of
    expectations.
    The main approach to construct confidence intervals for such parameters is the delta method.
    However, this asymptotic procedure yields intervals that may not be relevant for small sample sizes or, more generally, in a sequence-of-model framework that allows the expectation in the denominator to decrease to~$0$ with the sample size.
    In this setting, we prove a generalization of the delta method for ratios of expectations and the consistency of the nonparametric percentile bootstrap.
    We also investigate finite-sample inference and show a partial impossibility result: nonasymptotic uniform confidence intervals can be built for ratios of expectations but not at every level.
    Based on this, we propose an easy-to-compute index to appraise the reliability of the intervals based on the delta method.
    Simulations and an application illustrate our results and the practical usefulness of our rule of thumb.


    \hspace{0.7em}
\end{abstract}
\end{titlepage}

\newpage

\section{Introduction} 
\label{sec:introduction}

In applied econometrics, the prevalent method for constructing confidence intervals (CIs) is asymptotic: the theoretical guarantees for most CIs used in practice hold only when the number of observations tends to infinity.
For a large class of parameters, the construction of asymptotic CIs also relies on the delta method.
In this paper, we focus on parameters that can be expressed as ratios of expectations for which the delta method is a standard procedure to conduct inference.
The objective is twofold: study the behavior of the delta method and other confidence intervals in some difficult settings and provide tools to detect cases in which the delta method may behave poorly.

\medskip

Many popular parameters in economics take the form of ratios of expectations.
Typical examples are conditional expectations since any conditional expectation with a discrete conditioning variable, or a conditioning event, can be written as a ratio of unconditional expectations.
For instance, assume that we observe an
independent and identically distributed
(i.i.d.) sample of individuals indexed by $i \in \{1,\ldots,n\}$ with $W_i$ the wage of an individual and $D_i$ an indicator equal to~$1$ whenever individual~$i$ belongs to some treatment group, say a training program; $0$ otherwise.
Suppose you are interested in the average wage of participants in the program.
We have
$\Expec\left[W \mid D=1\right]
= \Expec\left[W D\right] /\,\Expec\left[D\right]$ as $D$ is binary.

\medskip

Most confidence intervals used in practice are based on asymptotic justifications, hence possible concerns as regards their finite-sample reliability.
For ratios of expectations, we document this issue on simulations (see Section~\ref{sub:problem_slow_convergence}).
One of our findings is that the coverage of the CIs based on the delta method happens to be far below their nominal level, even for large sample sizes, when the expectation in the denominator is close to~$0$.\footnote{The definitions of coverage and other fundamental properties of confidence intervals are recalled in Appendix~\ref{appendix:sec:definitions} with the conventions that we use.}
For some scenarios, these asymptotic CIs require above 100,000 observations to get reasonably close to their nominal level.
Yet, denominators close to~$0$ are not unusual in practice.
Coming back to the treatment/wage example,
a small denominator would correspond to a binary treatment with a low participation rate.

\medskip

In order to deal with that issue, we consider sequences of models, namely we authorize the distribution of the observations to change with the sample size.
This framework enables to formalize in an asymptotic way the idea of a denominator close to~$0$.
Indeed, in a standard asymptotic viewpoint, with the expectation in the denominator different from~$0$, all parameters are fixed and well-defined.
Hence, $n$~always grows large enough so that empirical means are close to their expectations and the CIs based on the delta method are valid.
In other words, the signal that we want to estimate is constant while the noise goes to~$0$, and therefore the problem vanishes in this asymptotic perspective.
We would like to model more difficult cases, in which the signal can go to~$0$ as well.
This is precisely what the sequence-of-model set-up allows.\footnote{This can also rationalize the practice of applied social researchers (see Example~\ref{example:practice_socialresearchers}). The heuristic idea is that researchers can consider narrower effects as the data gets richer.}
This is similar to some frameworks that have been developed for weak instrumental variables (IV), see notably~\cite{staiger1997,stock2005testing,andrews2018_weak}.

\medskip

In this literature, another approach does not consider sequences of models but designs ``robust'' procedures that allow to be exactly in the problematic case, namely a null covariance between the instrument and the endogenous regressor (see~\cite{anderson1949estimation}).
In this case, the parameter of interest is unidentified.
In contrast with the weak IV framework, it is worth noting that for ratios in general the parameter of interest is not even defined when the denominator is exactly equal to~$0$. As a consequence, such an approach seems difficult to extend to our problem.

\medskip

In our setting, it is unclear, even asymptotically, what the properties of the CIs based on the delta method are when the expectation in the denominator tends to~$0$.
We show that usual CIs can fail and the limiting law of ${\thetahatn - \theta_n}$ may not be Gaussian anymore, denoting by $\theta_n$ the ratio of expectations and $\thetahatn$ its empirical counterpart.
In some cases, the difference ${\thetahatn - \theta_n}$ may actually have a Cauchy limit, as can be found in the weak IV literature.

\medskip

We show in this sequence-of-model framework that confidence intervals provided by the nonparametric percentile bootstrap have the same asymptotic properties as the ones obtained with the delta method.
Simulations support that claim and even suggest the former have better coverage than the latter in finite samples.

\medskip

Even in standard settings with a fixed but small denominator, simulations document that asymptotic-based CIs may require very large sample sizes to attain their nominal level.
This suggests to study more in details nonasymptotic inference.
More precisely, we construct finite-sample CIs, extending old-established concentration inequalities for means to ratios of means.
Concentration inequalities for the mean refer to upper bounds on the probability that an empirical mean departs from its expectation more than a given threshold.
Such inequalities permit to construct confidence intervals valid for any sample size and for large classes of probability distributions (see in particular \cite{boucheron2013concentration}).
To our knowledge, there is no such result for ratios.
We consider distributions within a class characterized by a lower bound on the first moment for the denominator variable, and an upper bound on the second moment for both the numerator and denominator variables.\footnote{We refer to this setting as the ``Bienaymé-Chebyshev'' (BC) case. In Appendix~\ref{appendix:sec:results_hoef_framework}, we present similar results for distributions whose supports are bounded (``Hoeffding'' case).}

\medskip

One additional result highlights there exists a critical confidence level, above which it is not possible to construct nonasymptotic CIs, uniformly valid on such classes, and that are almost surely bounded under every distribution of those classes.
More precisely, we exhibit explicit upper and lower bounds on this critical confidence level: the former is a threshold above which we show it is impossible to construct such CIs; the latter is a threshold below which we show how to construct them.

\medskip

These ideas closely relate to some impossibility results as regards the construction of confidence intervals.
A large share of the research effort has concentrated on the problem of constructing confidence intervals for expectations.
In an early contribution, \cite{bahadur1956} show that, when $\Pfam$ is the set of all distributions on the real line with finite expectation, the parameter of interest~$\theta(P)$ is the expectation with respect to a distribution~${P \in \Pfam}$ and $\Theta=\Rb$, a confidence interval built from an i.i.d. sample of ${n \in \Nstar}$ observations that has uniform coverage ${1-\alpha}$ over~$\Pfam$ must contain any real number with probability at least~${1-\alpha}$.
Broadly speaking, any confidence interval must have infinite length with positive probability for every $P\in\Pfam$ to ensure a coverage of~${1-\alpha}$.

\medskip

Stronger results can be derived when one further restricts $\Pfam$ or $\Theta$.
When $\Pfam$ is taken to be the set of all distributions on the real line with variance uniformly bounded by a finite constant,
it is possible to show (using the Bienaym\'e-Chebyshev inequality) that for every $n \in \Nstar$ and every $\alpha\in(0,1)$, there exists a confidence interval that is almost surely bounded under every $P\in\Pfam$ and has coverage $1-\alpha$.
In this case, the obtained CIs have the advantage that their length shrinks to~$0$ at the optimal rate~$1/\sqrt{n}$.
But on the downside, they are not of size~$1-\alpha$, even asymptotically, except for some extreme distributions.
This means that they tend to be conservative in practice.

\medskip

A strand of the literature has also investigated more complex problems in which $\theta(P)$ is not restricted to being an expectation.
For general parameters, \cite{dufour1997} derives a generalization of \cite{bahadur1956}.
An implication of the results in \cite{dufour1997} is the existence of an impossibility theorem for ratios of expectations.
Let $P$ be a distribution on $\Rb^2$ with marginals $\distributionX$ and $\distributionY$.
If $\theta(P)=\Expec_{\distributionX}\left[X\right] / \,
\Expec_{\distributionY}\left[Y\right]$, then for every $\alpha\in(0,1)$, it is impossible to build nontrivial CIs of coverage $1-\alpha$ when $\Pfam$ is the set of all distributions on $\Rb^2$ with finite second moments and
$\Theta = \left\{\theta=\Expec_{P_X}\left[X\right] / \,
\Expec_{P_Y}\left[Y\right]:\left(\Expec_{P_X}\left[X\right],\Expec_{P_Y}\left[Y\right]\right)\in\Rb\times\Rstar\right\}$.
As will be explained below, this impossibility result disappears as soon as $\Pfam$ is chosen such that $\left|\Expec_{P_Y}\left[Y\right]\right|$ is bounded away from~$0$ uniformly over~$\Pfam$.
Interestingly, the impossibility breaks down only partly in the sense that there remains an upper bound on confidence levels (that depends on $n$) above which it is impossible to build nontrivial CIs.

\medskip

Other interesting results can be found in~\cite{romano2000} and~\cite{pinelis2016}.
\cite{romano2000} construct nonasymptotic valid confidence intervals that happen to be also asymptotically optimal. However, they only consider expectations.
\cite{pinelis2016} study smooth functions of a vector of means and give bounds on the distance between the distribution of the normalized and centered estimator and its Gaussian limiting distribution.
Nonetheless, the authors do not link their results to the construction of confidence intervals.

\medskip

In the light of that existing literature, our nonasymptotic findings can be interpreted as a partial impossibility result.
Indeed, even if we assume a known positive lower bound on the expectation in the denominator, the limitation on the attainable coverage of our nonasymptotic CIs remains.
That point complements \cite{dufour1997}: for a given sample size $n$, interesting CIs can be built but not at every confidence level.
By contrast, provided the expectation in the denominator is not null, the delta method gives CIs at every confidence level, but their coverage is only asymptotic.

\medskip

To bridge this gap, we suggest a rule of thumb to assess the reliability of the delta method for ratios of expectations in finite samples.
The heuristic idea is simply, for a given sample, to compute an estimator of the lower bound on the above-mentioned critical confidence level.
This lower bound can be seen as a conservative value for the unknown critical level, which is a necessary criterion to conduct valid inference in finite samples uniformly over a given class of distributions.
Hence, for any desired level higher than this bound, the CIs based on the delta method cannot reach this desired uniform level in finite samples.
We illustrate the empirical usefulness of that rule of thumb on simulations and with an application to gender wage disparities in France for the years 2010-2017.

\medskip

The rest of the paper is organized as follows.
Section~\ref{sec:framework} details our framework and assumptions.
In Section~\ref{sec:about_delta_method}, we illustrate the weaknesses of the CIs based on the delta method with a denominator ``close to 0'' on simulations and detail the asymptotic behavior of the delta method and of the nonparametric percentile bootstrap in our sequence-of-model setting.
Section~\ref{sec:construct_cis} is devoted to the construction of nonasymptotic confidence intervals and presents a lower bound on the aforementioned critical confidence level.
In Section~\ref{sec:about_alpha}, we derive an upper bound on the critical confidence level as well as a lower bound on the length of nonasymptotic CIs.
This section also includes the description of a practical index to gauge the soundness of the CIs based on the delta method in finite samples.
Section~\ref{sec:numerical_app} present simulations and an application to a real dataset to illustrate our methods.
Section~\ref{sec:conclusion} concludes.
General definitions about confidence intervals are recalled in Appendix~\ref{appendix:sec:definitions}.
The proofs of all results are postponed to Appendix~\ref{appendix:sec:proofs_BC_case}.
Additional results under an alternative set of assumptions (``Hoeffding'' case) are detailed in Appendix~\ref{appendix:sec:results_hoef_framework}.
Appendix~\ref{appendix:sec:additional_simulations} presents supplementary simulations.


\section{Our framework} 
\label{sec:framework}

Throughout the paper, for any random variable $U$ and $n$ i.i.d. replications $(U_{1,n},\ldots,U_{n,n})$, we denote by $\overline{U}_n$ the empirical mean of $U$, that is ${n^{-1} \sum_{i=1}^n U_{i,n}}$.
Assumption~\ref{hyp:basic_dgp_XY} defines our sequence-of-model framework and provides the basic requirements to state our asymptotic results.

\begin{hyp}
\label{hyp:basic_dgp_XY}
    For every ${n \in \Nstar}$, we observe a sample
    $(X_{i,n}, Y_{i,n})_{i=1,\ldots,n} \simiid \distributionXYn$, where $\distributionXYn$ is a given distribution on $\Rb^2$ that satisfies $\expecYn > 0$, $\expecXnsq < + \infty$, and $\expecYnsq < + \infty$.
\end{hyp}

Remark that~$n$ indexes both the distribution~$\distributionXYn$ of the observations in this model and the number of observations~$n$.
This encompasses the standard i.i.d. set-up if the distribution does not change with~$n$: for every ${n \in \Nstar}$, ${\distributionXYn = \distributionXY}$ for some given distribution~$\distributionXY$.
As we assume the existence of a finite expectation, we can consider ${\expecYn \geq 0}$ without loss of generality.\footnote{Otherwise, we simply replace $Y_{i,n}$ by its opposite $-Y_{i,n}$.}
In order to have properly defined ratios of interest, we need to assume away a null denominator, namely suppose that for every ${n \in \Nstar}$, ${\expecYn > 0}$.

\medskip

\begin{example}[Sequences of models and the practice of applied researchers]
    \label{example:practice_socialresearchers}
    ~ \newline
    Researcher may look at the average value of a variable~$A_{i,n}$ of interest in a subgroup of the data. Subgroups could be defined as the intersections of, say, time, geographical area, gender, age, income brackets and so on. As the number of observations~$n$ grows, it is possible to consider subgroups~$g_n$ that become thinner and thinner (intersection of more and more variables for instance). This practice could be modelled as estimating $\theta_n :=
    \Expec\left[A_{i,n} \mid G_{i,n} = 1\right]
    = \Expec\left[A_{i,n} G_{i,n}\right] / \, \Prob\left(G_{i,n}=1\right)$ where~$G_{i,n}$ is a binary variable that is equal to 1 if an individual $i$ belongs to the subgroup $g_n$.
    This corresponds to our framework denoting $X_{i,n} := A_{i,n} \times G_{i,n}$ and $Y_{i,n} := G_{i,n}$.
\end{example}

\medskip

To derive our nonasymptotic results,  Assumption~\ref{hyp:basic_dgp_XY} has to be strengthened.

\begin{hyp}
    \label{hyp:bt_dgp_XY}
    For every~${n \in \Nstar}$, there exist positive finite constants $\lowerYn$, $\upperXn$, and $\upperYn$ such that (i) ${\expecYn \geq \lowerYn > 0}$, (ii) $\expecXnsq \leq \upperXn$ and $\,\expecYnsq \leq \upperYn$.
\end{hyp}

Note that in practice, the value of the constants~$\lowerYn$, $\upperXn$, and $\upperYn$ may not be available for practitioners.
This is the reason why, in Section~\ref{subsec:plug_in_estimators}, we propose heuristic methods that palliate the lack of knowledge of those constants.

\medskip

The first part of the assumption bounds the expectation of~$Y_{1,n}$ away from~$0$ while the second states that the second moments of $X_{1,n}$ and $Y_{1,n}$ are bounded.
These are necessary to derive nonasymptotic CIs with maintained coverage uniformly over a class of distributions and that are not trivial.
Otherwise, if ${\lowerYn=0}$ or in the absence of the upper bounds $\upperXn$ and $\upperYn$, the impossibility theorem of \cite{dufour1997} applies and prevents from constructing nontrivial CIs for any confidence level.
In a way, given this result, Assumption~\ref{hyp:bt_dgp_XY} can be seen as close to the minimal hypothesis that allows for the possibility of nontrivial confidence intervals with finite-sample guarantees for ratios of expectations.
Furthermore, the sequence-of-model framework allows $\lowerYn$ to decrease to~$0$, which enables us to study limiting cases close to but different from the problematic case~${\lowerYn=0}$.

\medskip

This set-up, where Assumptions~\ref{hyp:basic_dgp_XY} and~\ref{hyp:bt_dgp_XY} hold, is named the \textit{BC case} since it is possible under these assumptions to construct nonasymptotic CIs using the Bienaym\'e-Chebyshev inequality.
In Appendix~\ref{appendix:sec:results_hoef_framework}, we present an adapted version of our results under the assumption that $X_{1,n}$ and $Y_{1,n}$ have a bounded support instead of bounded second moments; a setting we call the \textit{Hoeffding case}.

\medskip

To sum up, Assumptions~\ref{hyp:basic_dgp_XY} and~\ref{hyp:bt_dgp_XY} define a set~$\Pfam$ of distributions for some constants $\lowerYn$, $\upperXn$ and~$\upperYn$.
For a distribution $\distributionXYn$ in $\Pfam$, the parameter of interest $\theta(\distributionXYn)$ is denoted ${\theta_n := \expecXn /\,\expecYn}$ with values in~$\Rb$.
To estimate this parameter, we consider its empirical counterpart~${\thetahatn := \meanX /\,\meanY}$.
We seek to construct confidence intervals~$\Cnalpha$ for~${\theta_n}$ with nominal level~${1-\alpha}$ based on this estimator.

\medskip

In practice, it is possible that ${\meanY = 0}$ and it may even happen with a strictly positive probability for non-continuous distributions of~$Y$.
The estimator $\thetahatn$ does not exist for such samples.
In such a case, it is difficult to construct meaningful confidence intervals.
Different conventions are possible:
\begin{itemize}
    \item We could choose to define ${\Cnalpha = \Rb}$. This entails that $\theta_n$ belongs to~$\Cnalpha$ by construction.
    We believe that such a choice would artificially improve the coverage of $\Cnalpha$ as it induces that the higher $\Prob(\meanY = 0)$, the better the interval in terms of coverage.
    \item We could choose ${\Cnalpha = \emptyset}$.
    The hypothesis ${\theta_n = \theta_0}$ would then be rejected for every ${\theta_0 \in \Rb}$ using the duality between tests and confidence intervals.
    We would also like to avoid this situation because it may not be reasonable to always reject for the mere reason that $\theta_n$ cannot be estimated in the sample.
    \item Other choices are possible, for example ${\Cnalpha = \{0\}}$, but they do not seem sensible either since there is no reason to select only~$0$ in our confidence interval, especially if~${\meanX \neq 0}$.
\end{itemize}

For these considerations, we choose to let $\Cnalpha$ \textit{undefined} whenever~${\meanY = 0}$, following the convention that ratios $x/0$ are undefined for any real~$x$.\footnote{When facing ${\meanY = 0}$, applied researchers may use other estimators. For instance, one could consider sub-samples (possibly several and combine them in some way) of the data for which the empirical mean in the denominator differs from~$0$.
Nevertheless, the construction of satisfactory estimators in this case lies beyond the scope of this paper.}
In practice, when given a realization $\omega \in \Omega$ and a real $a \in \Rb$, we either know that $a$ belongs to $C_{n, \alpha}(\omega)$, or we know that $a$ does not belong to $C_{n, \alpha}(\omega)$, or $C_{n, \alpha}(\omega)$ is undefined.
As a consequence, we have the decomposition $\Omega = \{\omega: a \in C_{n, \alpha}(\omega) \} \sqcup
\{\omega: a \notin C_{n, \alpha}(\omega) \} \sqcup
\{\omega: C_{n, \alpha}(\omega) \text{ undefined}\},$
where $\sqcup$ denotes the disjoint union of sets.
This means that $\Prob \{ a \in C_{n, \alpha} \} + \Prob \{ a \notin C_{n, \alpha} \} + \Prob \{ C_{n, \alpha} \text{ undefined}\} = 1$.



\section{Limitations of the delta method: when are asymptotic confidence intervals valid?} 
\label{sec:about_delta_method}

In practice, for a sample of size~$n$, the coverage of asymptotic CIs may be well below their nominal level~${1-\alpha}$.
Intuitively, this phenomenon should be driven by ``problematic'' distributions in $\Pfam$ in the following sense:
when the true distribution $P$ is close to the boundary of the class $\Pfam$, the probability $c(n,P) := \Prob_{P^{\otimes n}}\left(\Cnalpha \ni \theta(P)\right)$ may be much smaller than~${1-\alpha}$.\footnote{Recall that in the nonasymptotic approach, the \textit{coverage} of any given confidence interval $\Cnalpha$ is defined as the infimum of $c(n,P)$ for $P$ ranging over the studied class~$\Pfam$ of distributions.}

\medskip

In Section~\ref{sub:problem_slow_convergence}, with $\Cnalpha$ the confidence interval based on the delta method, we illustrate on simulations that $c(n,P)$ can fail to match~${1-\alpha}$
when the expectation in the denominator is fixed close to~$0$.
In other words, it may require a very large number of observations to make reasonable the asymptotic approximation.
In Section~\ref{sub:problem_sequence_models}, we investigate a more serious issue: in the sequence-of-model framework, we let the expectation in the denominator not only be small but converge to~$0$ as $n$~increases.
We show on simulations that depending on the speed at which the denominator goes to~$0$, $c(n,P)$ can either converge to the nominal level (more or less quickly) or even not converge at all to this target.
This sheds light on a partial failure of the delta method when the denominator goes to~$0$ that we derive formally in Section~\ref{sub:extension_delta_method_sequence_models}.
Finally, in Section~\ref{sub:bootstrap_CI}, we show the asymptotic consistency of the nonparametric percentile bootstrap (also known as Efron's percentile bootstrap) in this sequence-of-model framework.

\subsection{Asymptotic approximation takes time to hold} 
\label{sub:problem_slow_convergence}

In this subsection, we consider the i.i.d. case.\footnote{For every $n \in \Nstar$, $\distributionXYn$ is identical, hence denoted  $\distributionXY$. To simplify notations, we also denote by $(X,Y)$ a random vector following~$\distributionXY$.}
Under Assumption~\ref{hyp:basic_dgp_XY}, asymptotic confidence intervals are easily obtained combining the multivariate central limit theorem (CLT) and the delta method:
\begin{equation}
\label{eq:delta_method_simple_case_n_for_sample_size}
    \sqrt{n}\left(\frac{\overline{X}_n}{\overline{Y}_n}-\frac{\expecX}{\expecY}\right) \convD \mathcal{N}\left(0, \Sigma\right),
\end{equation}
where ${\Sigma = \Var[X]/ \,\expecY^2
+ \expecX^2\VarY / \,\expecY^4
- 2 \, \Cov\left[X,Y\right] \expecX / \,\expecY^3}$ and in practice is replaced by a consistent estimate (Slutsky's lemma).

\medskip

\begin{figure}[ht]
\centering
\includegraphics[width=0.8\textwidth]{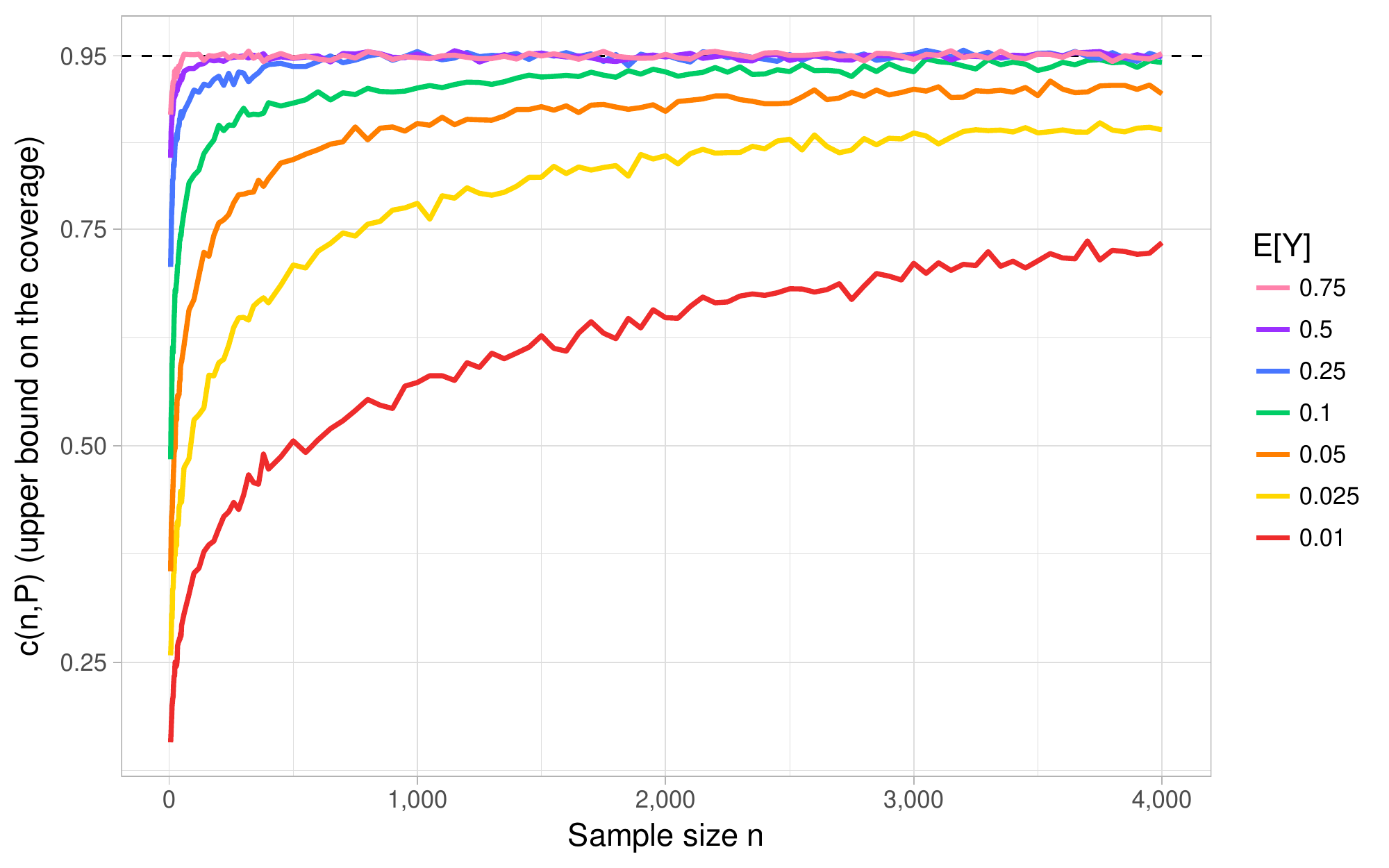}
\caption{\small{$c(n,P)$ of the asymptotic CIs based on the delta method as a function of the sample size~$n$. \newline
Specification: $\forall n \in \Nstar, {\distributionXYn = \mathcal{N}(1,1) \otimes \mathcal{N}(\expecY,1)}$. The nominal pointwise asymptotic level~is set to $0.95$. For each pair $(\expecY, n)$, the coverage is obtained as the mean over 5,000 repetitions.}}
\label{fig:slow_convergence_normal_case_EX1_VX1_VY1_CorrXY0}
\end{figure}

To assess the quality of the CI based on~\eqref{eq:delta_method_simple_case_n_for_sample_size}, we compute its $c(n,P)$ using simulations for different sample sizes~$n$ and distributions~$P$ and compare it to the nominal level.
By definition, the pointwise coverage $c(n,P)$ forms an upper bound on the uniform coverage.
In our simulations, we choose the level ${1-\alpha = 95\%}$.
For different sample sizes~$n$ and values of~$\expecY$, we draw $M=$ 5,000 i.i.d. samples of size $n$ following $\mathcal{N}(1,1) \otimes \mathcal{N}(\expecY,1)$.
We compute $c(n,P)$ for the interval based on the delta method for every pair $(n, \, \expecY)$ using the 5,000 replications.
The expectation $\expecY$ ranges from $0.01$ (the denominator is close to~$0$) to $0.75$ (the denominator is far from~$0$).
Figure~\ref{fig:slow_convergence_normal_case_EX1_VX1_VY1_CorrXY0} sums up the results.
For every~$n$, it turns out that the closer $\expecY$ to~$0$, the smaller the $c(n,P)$ of the delta method.
When $\expecY=0.01$, we observe that $c(n,P)$ gets close to the nominal level only for $n$ above~300,000.
Additional simulations indicate that the phenomenon is robust across different choices of the distribution~$\distributionXY$ (see Section~\ref{appendix:sec:additional_simulations}).


\subsection{Asymptotic results may not hold in the sequence-of-model framework} 
\label{sub:problem_sequence_models}

\begin{figure}[htb]
    \centering
    \includegraphics[width=0.8\textwidth]{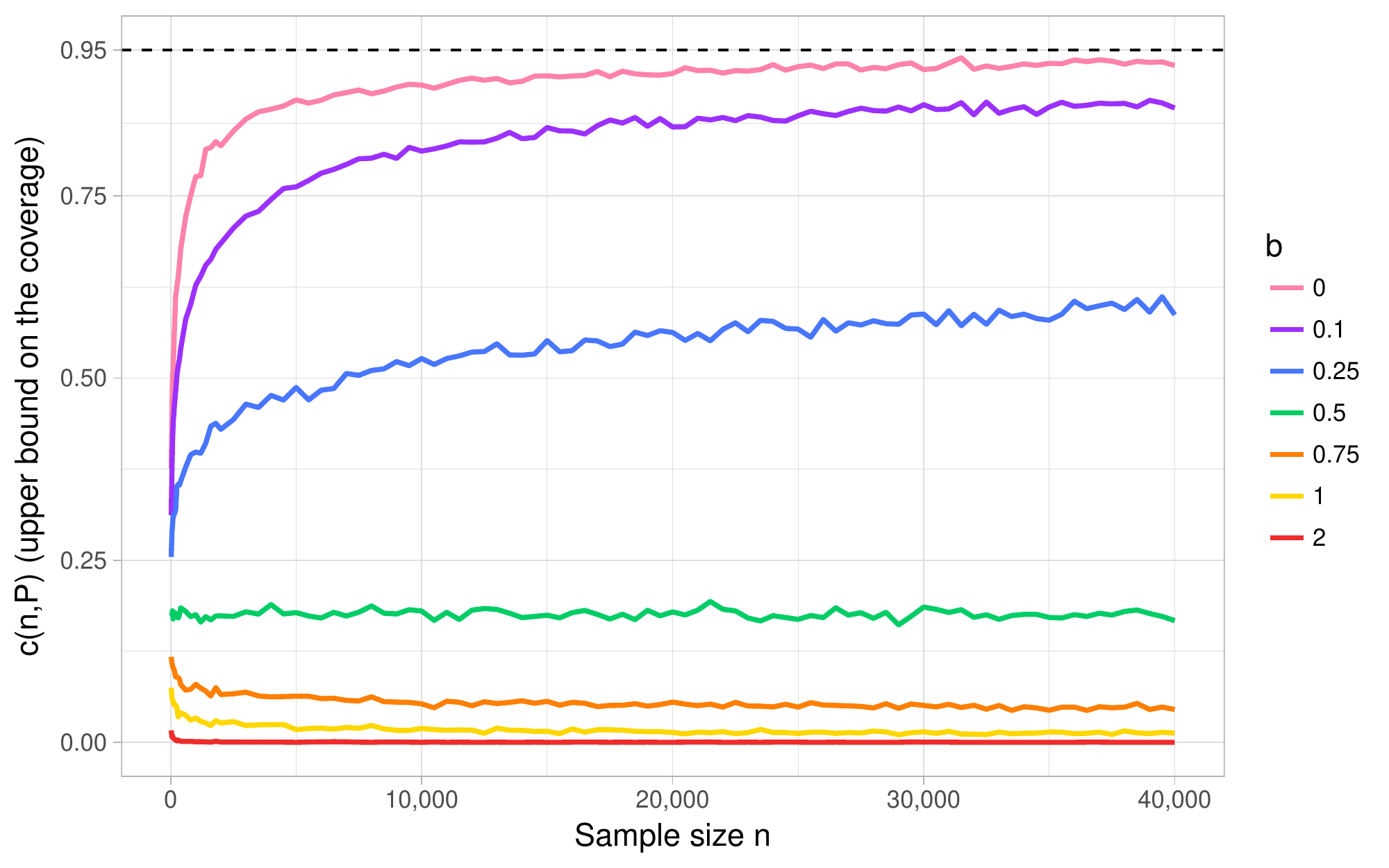}
    \caption{\small{$c(n,P)$ of the asymptotic CIs based on the delta method as a function of the sample size~$n$. \newline
    Specification: $\forall n \in \Nstar, {\distributionXYn = \mathcal{N}(1,1) \otimes \mathcal{N}(Cn^{-b},1)}$, with $C=0.025$. The nominal \mbox{pointwise} asymptotic level is set to $0.95$. For each pair ${(b, n)}$, the coverage is obtained as the mean over 5,000 repetitions.}}
    \label{fig:sequence_modele_normal_case_EX1_VX1_VY1_CorrXY0_C0025}
\end{figure}

Unlike the result displayed in~\eqref{eq:delta_method_simple_case_n_for_sample_size}, it is unclear how
$\sqrt{n}\left(\meanX/\,\meanY-\expecX/\,\expecY\right)$ behaves asymptotically when we consider sequences of models such that the expectation in the denominator tends to~$0$ as $n$~increases.
For a given specification, Figure~\ref{fig:sequence_modele_normal_case_EX1_VX1_VY1_CorrXY0_C0025} shows the $c(n,P)$ of the CIs based on the delta method when ${\expecYn = Cn^{-b}}$ where $C$ is set to~$0.025$ and $b$~varies.
For a speed ${b \geq 1/2}$ (i.e. faster than the usual rate of the CLT), the pointwise coverage $c(n,P)$ of the asymptotic CIs obtained by~\eqref{eq:delta_method_simple_case_n_for_sample_size} is not good in the sense that it is far lower than the nominal level $1-\alpha$ and it does not converge to the latter.
Our simulations even suggest that the coverage tends to $0$ for $b>1/2$.
For $b<1/2$, the upper bound $c(n,P)$ on the coverage of the delta method seems to tend to $1-\alpha$.
Yet, in line with Figure~\ref{fig:slow_convergence_normal_case_EX1_VX1_VY1_CorrXY0}, the validity of the asymptotic approximation requires very large sample sizes.

At this stage, Figure~\ref{fig:sequence_modele_normal_case_EX1_VX1_VY1_CorrXY0_C0025} presents some evidence that the CIs based on the delta method need to be adapted for sequences of models and that the rate of decrease toward~$0$ of the expectation $\expecYn$ matters.
The next subsection details formal results in this set-up.



\subsection{Extension of the delta method for ratios of expectations in the sequence-of-model framework}
\label{sub:extension_delta_method_sequence_models}

We are interested
in the asymptotic distribution, as $n$ tends to infinity, of the real random variable
$S_n := \sqrt{n} \left({\meanX}/\,{\meanY} - {\expecXn}/\,{\expecYn} \right)$.
The following theorem states the asymptotic behavior of $S_n$ according to the comparison of $\VarYn \, /\sqrt{n}$ and $\expecYn$ under a multivariate Lyapunov condition.
It is proved in Section~\ref{proof:thm:delta_method_sequence_models}.

\medskip

We show that in some cases $|S_n| \convAS + \infty$. It is then impossible to state the limiting distribution $S_n$ in the traditional sense.
Despite that, we can still get a more precise result looking at the subsequent terms in the asymptotic expansion of~$S_n$.
Such an asymptotic expansion is complicated to state, especially in our sequence-of-model framework, since the distributions $\distributionXYn$ change with $n$ without any link from one to the next.
To overcome this problem, we consider equivalents in distribution of $S_n$ in the following sense.
%
%
%
%
We say that two sequences of random variables $S_n$ and $T_n$ are \textit{equivalent in distribution} if there exist a probability space $\tilde \Omega$ and two sequences of random variables $\tilde S_n, \tilde T_n$ such that $\forall n \in \Nstar$,
$S_n \equalD \tilde S_n$ and $T_n \equalD \tilde T_n$, and $\tilde S_n$ is equivalent to $\tilde T_n$ almost surely as~$n \to \infty$.
This means that for almost every $\tilde \omega \in \tilde \Omega$, $\tilde S_n(\tilde \omega)$ is equivalent to $\tilde T_n(\tilde \omega)$ (considered as deterministic sequences of real numbers).
This notion enables to formalize the link between $S_n$ and a simpler expression~$T_n$.

\begin{thm}
\label{thm:delta_method_sequence_models}


    Let Assumption~\ref{hyp:basic_dgp_XY} hold and
    (i) $\Var[(\gamma_{X,n} X_{1,n} \, , \, \gamma_{Y,n} Y_{1,n})] \to V$ as ${n \to \infty}$ for some positive sequences $\{\gamma_{X,n}\}_{n \in \Nstar}$ and $\{\gamma_{Y,n}\}_{n \in \Nstar}$ where $V$ is a definite positive $2 \times 2$ matrix,
    (ii)~$\sup_{n\in\Nstar}\Expec\left[|X_{1,n}|^3  \gamma_{X,n}^3 + |Y_{1,n}|^3 \gamma_{Y,n}^3\right] < +\infty$, and
    (iii) $\Prob(\meanY = 0) \to 0$ as $n \to \infty$.

    \noindent Denote the signal-to-noise-ratio by
    ${\SNRn := \expecYn / (V_{2,2}^{1/2} n^{-1/2} \gamma_{Y,n}^{-1})}$.

    Then, the sequence of random variables
    $S_n := \sqrt{n} \left(\meanX/\,\meanY
    - \expecXn/\,\expecYn \right)$ satisfies as~${n \to \infty}$:
    \begin{enumerate}
        \item If ${\SNRn \to + \infty}$,
        then $S_n$ is equivalent in distribution to:
        \begin{equation*}
            \frac{\sqrt{n} \gamma_{X,n} (\meanX - \expecXn)}{\expecYn \gamma_{X,n}}
            - \frac{\sqrt{n} \gamma_{Y,n} (\meanY - \expecYn) \expecXn}{\expecYn^2 \gamma_{Y,n}}.
        \end{equation*}
    \item If there exists a finite constant ${C \neq 0}$ such that ${\SNRn \to C}$,
    then $S_n$ is equivalent in distribution to:
    \begin{align*}
        n \gamma_{Y,n} \expecXn
        &\Bigg( \frac{1}{C + \sqrt{n} \gamma_{Y,n} (\meanY - \expecYn)} - \frac{1}{C} \Bigg) \quad \\
        & \qquad \qquad + \frac{n \gamma_{X,n} (\meanX - \expecXn) \times \gamma_{Y,n}}
        {\big(C + \sqrt{n} \gamma_{Y,n} (\meanY - \expecYn) \big) \times \gamma_{X,n}}.
    \end{align*}
    \item If ${\SNRn \to 0}$,
    then $S_n$ is equivalent in distribution to:
    \begin{equation*}
        \sqrt{n} \left(
        \frac{\sqrt{n} \gamma_{X,n} (\meanX - \expecXn)} {\sqrt{n} \gamma_{Y,n} (\meanY - \expecYn)}
        \times \frac{\gamma_{Y,n}}{\gamma_{X,n}}
        - \frac{\expecXn}{\expecYn}
        \right).
    \end{equation*}
    \end{enumerate}
\end{thm}

\begin{figure}[p]
    {\centering
    \includegraphics[width=0.75\textwidth]{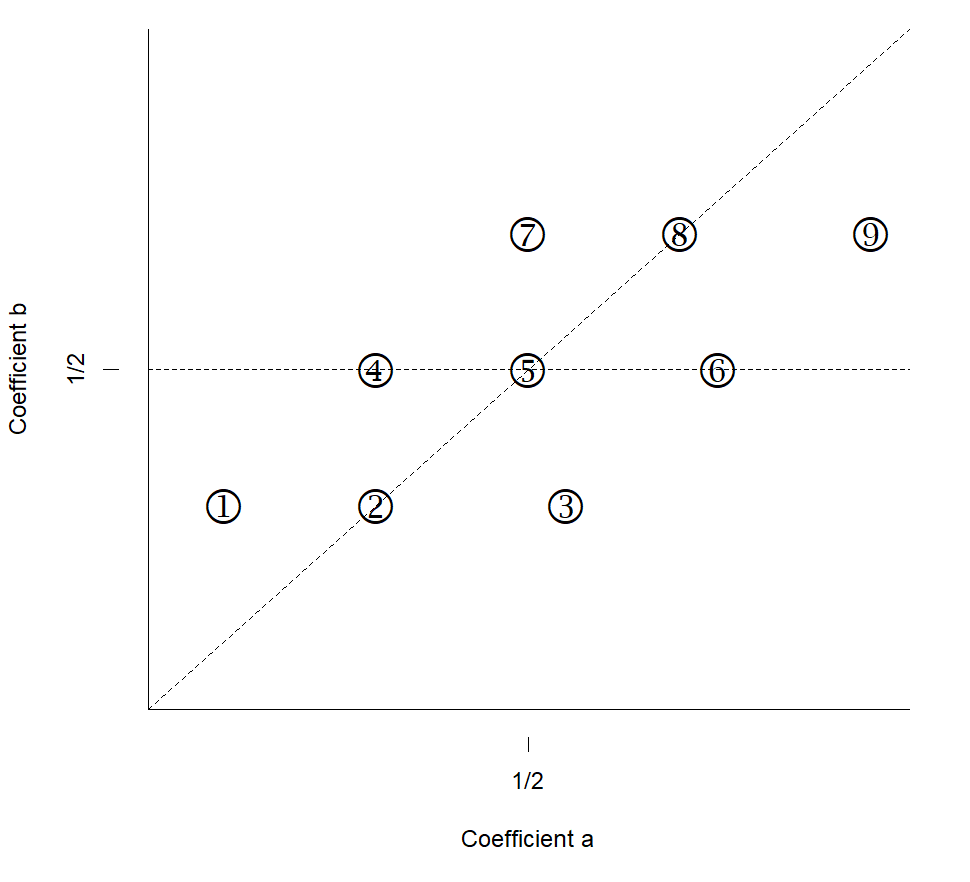}
    \caption{Separation between the different asymptotic regimes as a function of $(a,b)$ for fixed $(a',b')=(0,0)$, in the case where $\expecXn = C_1 / n^a$, $\VarX = 1 / n^{a'}$, $\expecYn = C_2 / n^b$, and $\VarY = 1 / n^{b'}$, ${(a,a',b,b') \in \Rb_+^4}$.}
    \label{fig:graph_transition_phase}

    \vspace{1cm}

    \renewcommand{\arraystretch}{2}

    \resizebox{\textwidth}{!}{%
    \begin{tabular}{c|ccc}
        & $a + b' < b + a'$ & $a + b' = b + a'$ & $a + b' > b + a'$ \\
        \hline
        $b > 1/2 + b'$
        & $n^{1/2 + b' - a'} W_1 / W_2$
        & $n^{1/2 + b' - a'} \big(W_1 / W_2 - C_1 / C_2 \big)$
        & $- n^{1/2 + b - a} C_1 / C_2$ \\
        $b = 1/2 + b'$
        & $n^{1-a+b'} \Big( C_1 / (C_2 + W_2) - C_1/C_2 \Big)$
        & $n^{1/2+b' - a'} \big(C_1 / (C_2 + W_2) $
        & $n^{1/2+b' - a'} \big(W_1 / (C_2 + W_2 n^{a'}) \big)$ \\
        && $- C_1/C_2
        + W_1 / (C_2 + W_2 n^{a'}) \big)$ & \\
        $b < 1/2 + b'$ & $n^{2b-a-b'} C_1 W_2 / C_2^2$
        & $n^{b - a'} (W_1 / C_1 - C_1 W_2 / C_2^2)$
        & $n^{b - a'} W_1 / C_1$
    \end{tabular}
    \label{tab:graph_transition_phase}
    }
    \captionof{table}{Limiting law of $S_n := \sqrt{n} \left(\meanX/\,\meanY - \expecXn/\,\expecYn \right)$ in the nine different regimes. The couple of variables $(W_1, W_2)$ follow the distribution $\Nc(0, V)$, where
    $V = \lim_{n\to+\infty} \Var\left[(n^{a'} X_{1,n} , n^{b'} Y_{1,n})\right]$.}
    }

\end{figure}

Theorem~\ref{thm:delta_method_sequence_models} can thus be interpreted as a generalization of the result given by the CLT and the delta method for ratios of expectations.
The sequence-of-model framework allows both the expectation and the variance in the denominator to tend to $0$.
In particular, this happens whenever $Y_{i,n}$ follows a Bernoulli distribution with a parameter $p_n$ tending to $0$, as detailed in Example~\ref{example:sequence_Bernoulli}.
For instance, when we estimate a conditional expectation with a discrete conditioning variable or a conditioning event, the denominator is an average of indicator variables that follow a Bernoulli distribution.
Figure~\ref{fig:graph_transition_phase} and its companion table highlight the different asymptotic regimes depending on the behaviors of $\{\expecXn\}_{n \in \Nstar}$, $\{\expecYn\}_{n \in \Nstar}$, $\{\gamma_{X,n}\}_{n \in \Nstar}$ and $\{\gamma_{Y,n}\}_{n \in \Nstar}$.

\medskip

The main takeaway of Theorem~\ref{thm:delta_method_sequence_models} is that when $\expecXn=C_1/n^a$, $\expecYn=C_2/n^b$ and $\VarY = C_3 / n^{b'}$ for some constants $C_1, C_2, C_3 \neq 0$, and $b<1/2+b'$,
$S_n$ properly renormalized by $n$ to some power still converges in distribution to a Normal random variable.
This can be explained using the signal-to-noise ratio (SNR) defined in Theorem~\ref{thm:delta_method_sequence_models}. Indeed, in this first case, the $\SNRn$ tends to $+\infty$:
the signal in the denominator (that is the expectation of $Y_{1,n}$) is asymptotically bigger than the noise (which is $1/(\gamma_{Y,n} n^{1/2})$ up to a constant factor).
Asymptotic inference based on the Normal approximation remains valid, even if the length of such confidence intervals may not decrease with the sample size~$n$.

\medskip

In all other cases, when the noise dominates in the denominator, $S_n$ converges weakly to a non-Gaussian distribution, in some cases to a generalized Cauchy distribution with parameters that depend on the data generating process (up to a normalization of some power of $n$).
By construction, when the noise dominates, we do not have much information and thus may not be able to conduct inference in these settings.
This echoes the impossibility results presented in Section~\ref{sec:about_alpha}, notably Remark~\ref{rem:snrtilde_limit}.
In the next section, we provide another method for constructing confidence intervals using the nonparametric percentile bootstrap.

\begin{example}
    When $Y_{1,n}$ follows a Bernoulli distribution with parameter $p_n$ in $(0,1)$, we are always in the first case of
    Theorem~\ref{thm:delta_method_sequence_models}, meaning that its expectation $p_n$ is always larger than the noise $\sqrt{p_n (1-p_n)/n}$.
    This latter formula is obtained by remarking that the standard deviation of $Y_{i,n}$ is $\sqrt{p_n (1-p_n)}$ so that
    $\gamma_{Y,n} = 1/\sqrt{p_n (1-p_n)}$.
    However, in order to satisfy the constraint $\Prob(\meanY = 0) \to 0$, we have to impose that $n p_n \to + \infty$.
    Therefore, when $p_n = n^{-b}$, confidence intervals based on the delta method will be pointwise consistent if $b < 1$.
    \label{example:sequence_Bernoulli}
\end{example}


\subsection{Validity of the nonparametric bootstrap for sequences of models}
\label{sub:bootstrap_CI}

In this part, we construct confidence intervals for ratios of expectations using Efron's percentile bootstrap.
This technique relies on the nonparametric bootstrap resampling scheme that we now recall.
We fix a number $B > 0$ of bootstrap replications.
For a given initial sample $(X_{i,n}, Y_{i,n}), i=1, \dots, n$, and a given integer $b$ smaller than $B$, we define the bootstrapped sample $(X_{i,n}^{(b)}, Y_{i,n}^{(b)}), i=1, \dots, n$, which is obtained by $n$ i.i.d. resampling from the initial sample, i.e. with replacement.
Let $\overline{X}_n^{(b)} := n^{-1} \sum_{i=1}^n X_{i,n}^{(b)}$ be the empirical mean of the numerator in the $b$-th bootstrapped sample (resp. $\overline{Y}_n^{(b)}$ for the denominator).

\medskip

Then, Efron's percentile bootstrap, also known as the nonparametric percentile bootstrap, consists in using the quantiles of the bootstrapped distribution conditional on the data to conduct inference.
More precisely, for every $\tau \in (0,1)$, let $q_\tau^{boot}$ denote the
quantile at level~$\tau$ of $\overline{X}_n^{(1)} / \, \overline{Y}_n^{(1)}$, which is estimated in practice by the empirical quantile at level $\tau$ of the bootstrapped statistics
$\big( \overline{X}_n^{(b)} / \, \overline{Y}_n^{(b)}\big)_{b = 1,\ldots,B}$.
For a given nominal level ${1-\alpha \in (0,1)}$, the confidence interval we consider is defined as $C_{n,\alpha}^{boot} := \big[ q_{\alpha/2}^{boot} \, , \, q_{1-\alpha/2}^{boot} \big]$.
The following theorem states the consistency of this interval.
It is proved in Section~\ref{proof:thm:validity_bootstrap_sequenceModels}.

\begin{thm}
\label{thm:validity_bootstrap_sequenceModels}



    Let Assumption~\ref{hyp:basic_dgp_XY} hold and
    (i) $\Var[(\gamma_{X,n} X_{1,n} \, , \, \gamma_{Y,n} Y_{1,n})] \to V$ as ${n \to \infty}$ for some positive sequences $\{\gamma_{X,n}\}_{n \in \Nstar}$ and $\{\gamma_{Y,n}\}_{n \in \Nstar}$ where $V$ is a definite positive $2 \times 2$ matrix,
    (ii) ${\sup_{n \in \Nstar} \Expec \Big[ (\gamma_{X,n} X_{1,n})^{4+\delta}
    + (\gamma_{Y,n} Y_{1,n})^{4+\delta} \Big] < + \infty}$ for some $\delta > 0$,
    (iii) $\Prob(\meanY = 0) \to 0$ as $n \to \infty$, and
    (iv) $\Prob(\meanY^{(1)} = 0) \to 0$ as $n \to \infty$.

    \noindent Denote the signal-to-noise-ratio by
    ${\SNRn := \expecYn / (V_{2,2}^{1/2} n^{-1/2} \gamma_{Y,n}^{-1})}$.

    If $\SNRn \to +\infty$, then for every~${\alpha \in (0,1)}$, the confidence interval $C_{n,\alpha}^{boot}$ is
    pointwise consistent
    at level~${1-\alpha}$, viz.
    $\Prob \big( C_{n,\alpha}^{boot} \ni \expecXn/\,\expecYn \big) \to 1 - \alpha \text{ as } n \to \infty.$
\end{thm}

The assumption $\Prob(\meanY^{(1)} = 0) \to 0$ is satisfied for a large set of cases, for instance when the variables $Y_{i,n}$ are continuous or when they follow a Bernoulli distribution with a parameter decreasing to~$0$ not too fast (see Example~\ref{example:sequence_Bernoulli_bootstrap} below).

\medskip

\begin{figure}[htb]
    \begin{minipage}[b]{0.5\linewidth}
    \includegraphics[width=1\linewidth]{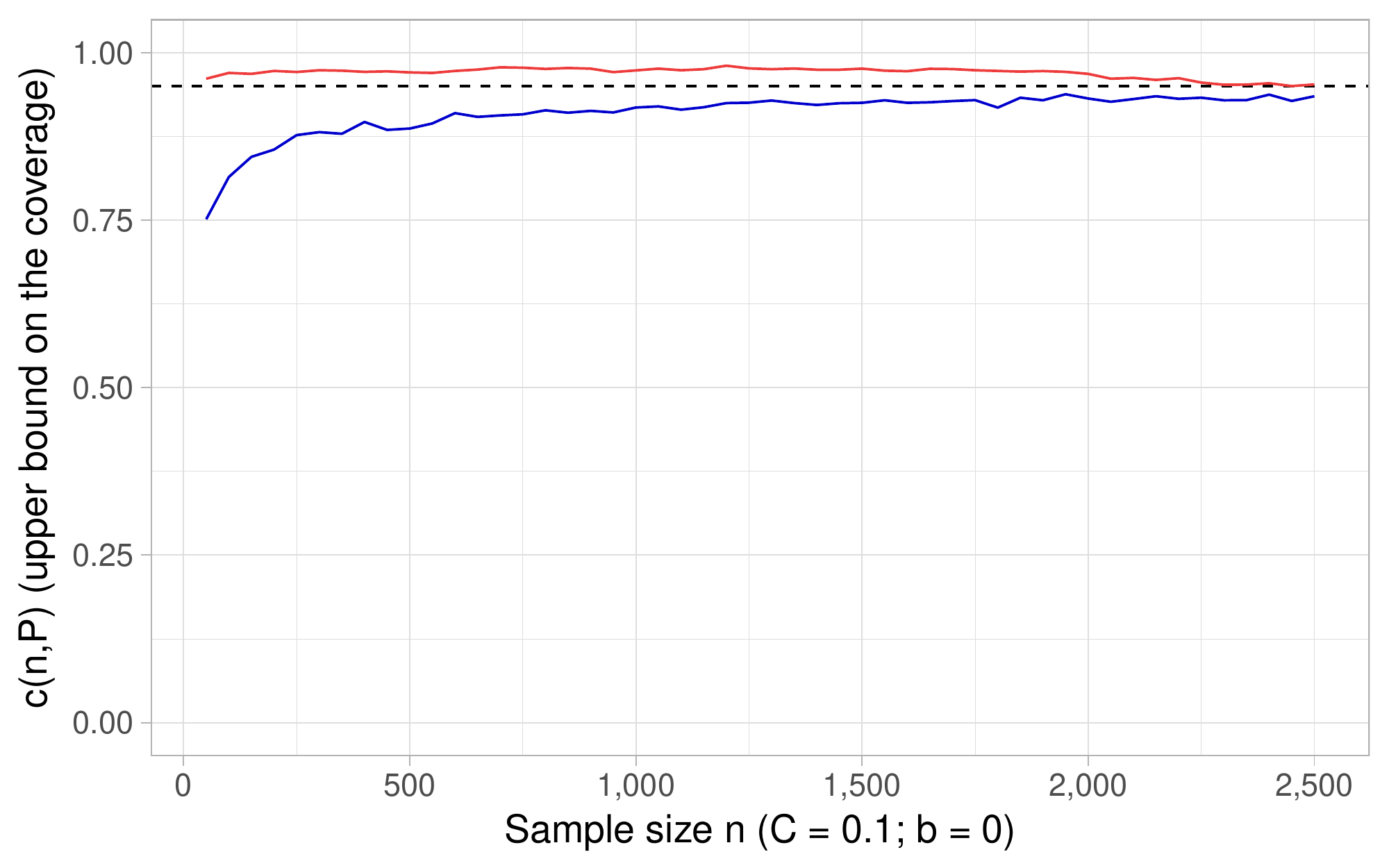}
    \end{minipage}
    \begin{minipage}[b]{0.5\linewidth}
    \includegraphics[width=1\linewidth]{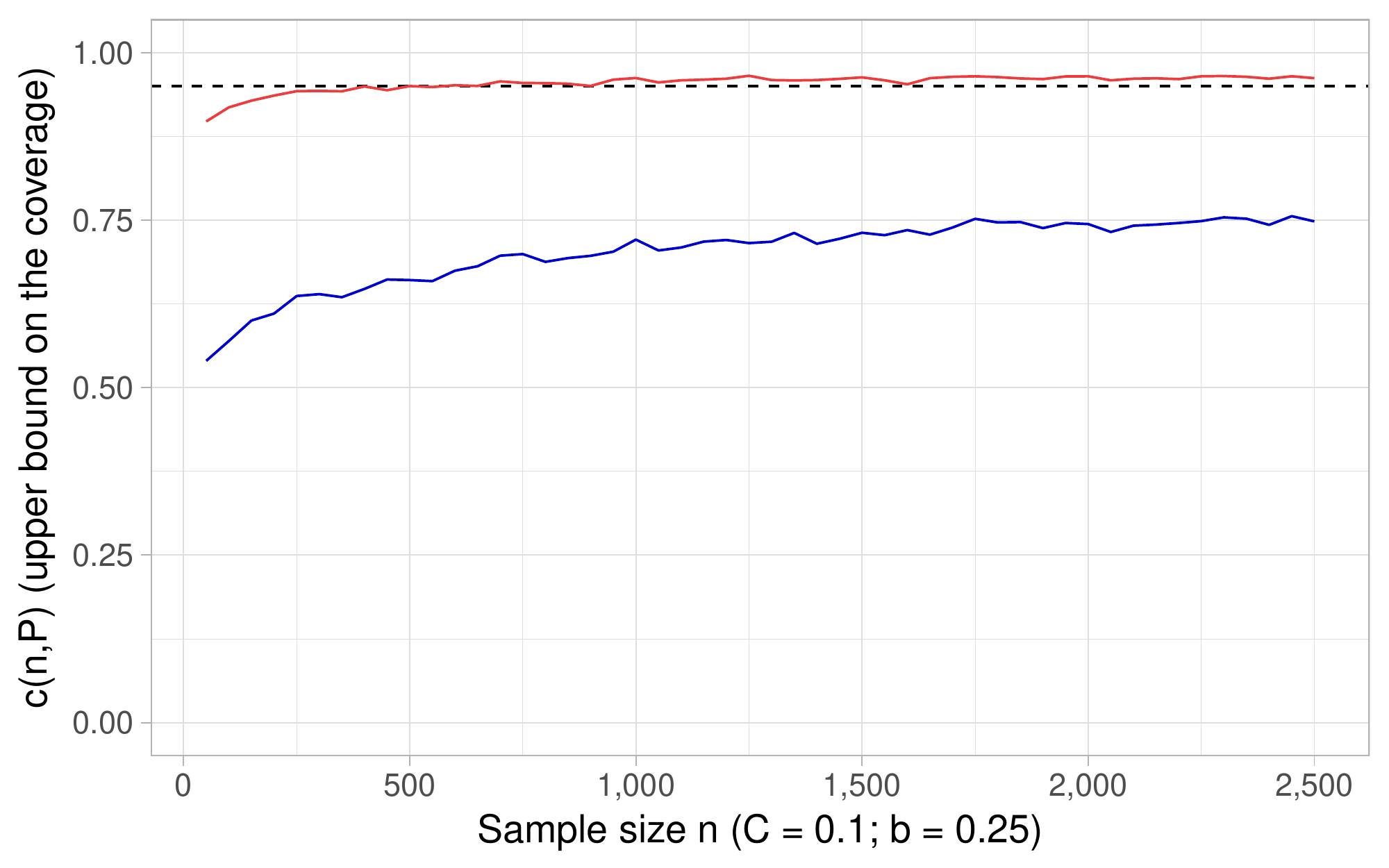}
    \end{minipage}
    \begin{minipage}[b]{0.5\linewidth}
    \includegraphics[width=1\linewidth]{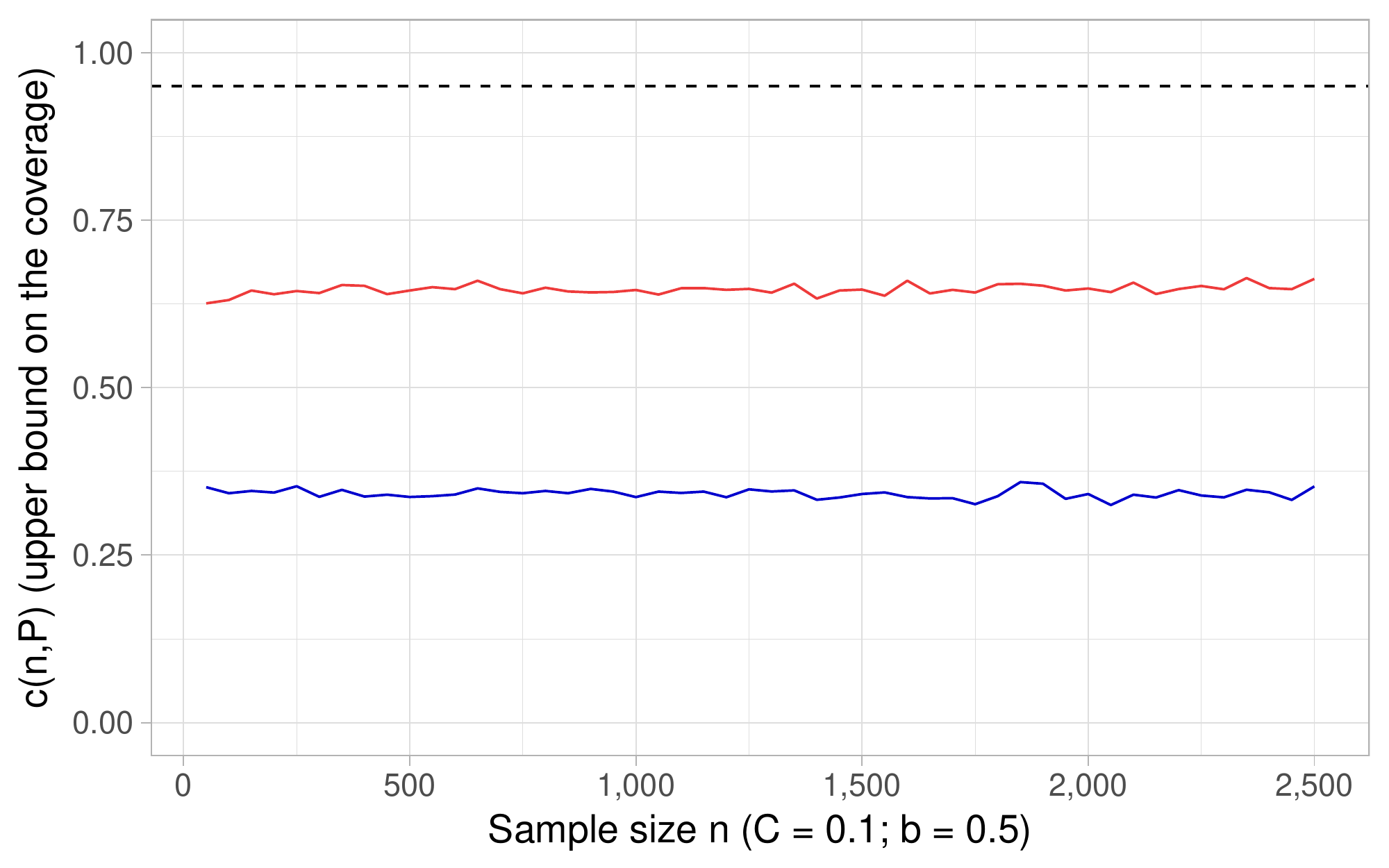}
    \end{minipage}
    \hfill
    \begin{minipage}[b]{0.5\linewidth}
    \includegraphics[width=1\linewidth]{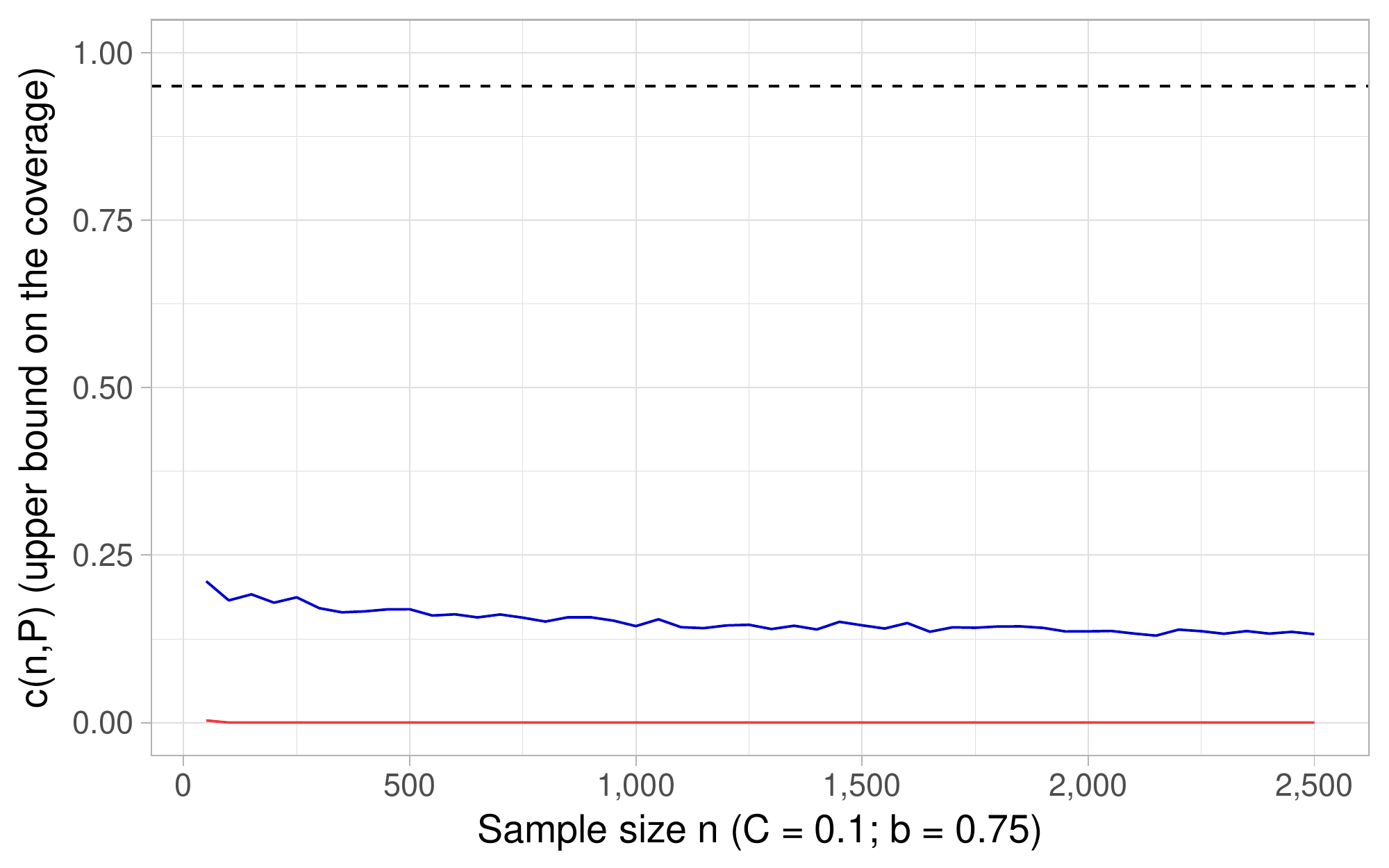}
    \end{minipage}
    \caption{\small{$c(n,P)$ of the asymptotic CIs based on the delta method (blue) and of the CIs constructed with Efron's percentile bootstrap using 2,000 bootstrap replications (red).\newline
    Specification: $\forall n \in \Nstar, {\distributionXYn = \mathcal{N}(1,1) \otimes \mathcal{N}(Cn^{-b},1)}$, with $C=0.1$ and $b \in \{0, 0.25, 0.5, 0.75\}$. The nominal pointwise asymptotic level is set to $0.95$. For each pair~${(b, n)}$, the coverage is obtained as the mean over 5,000 repetitions.}}
    \label{fig:bootstrap_DM_sequences_of_models_C01}
\end{figure}

Note that the moment condition of order $4+\delta$ is nearly sharp.
Indeed, the proofs require the strong law of large numbers for ${n^{-1}\sum_{i=1}^n X_{1,n}^2}$ and ${n^{-1}\sum_{i=1}^n Y_{1,n}^2}$.
As we are dealing with a triangular array of random variables,
Theorem 3.1 of~\cite{gut1992complete} shows that moments of order at least~$4$ are necessary, even in the simpler case where the distribution~$\distributionXYn$ does not depend on~$n$.

\begin{example}[Example~\ref{example:sequence_Bernoulli} continued]
    \label{example:sequence_Bernoulli_bootstrap}
    When $Y_{1,n}$ follows a Bernoulli distribution with parameter
    $p_n = 1/n^b$ for a given $b > 0$, the condition $\Prob(\meanY^{(1)} = 0) \to 0$ is satisfied when $b < 1$. We refer the reader to Section~\ref{proof:example:sequence_Bernoulli_bootstrap} for a proof of this claim.
\end{example}

In practice, even if the theoretical results of the delta method and of the bootstrap are valid under nearly the same set of assumptions, we observe in the simulations in Figure~\ref{fig:bootstrap_DM_sequences_of_models_C01} a gap between their pointwise coverage.\footnote{Additional simulations comparing the two types of asymptotic confidence intervals are presented in Appendix~\ref{appendix:sub:additional_simulations_DM_bootstrap}.}
This fact appears even when $\distributionXYn$ does not depend on $n$ (i.e. $b=0)$.
Nonetheless, the coverage gap between these two methods shrinks as $n$ increases provided~${b<0.5}$.
In the sequence of models where the denominator decreases slowly (i.e. $b=0.25$) in Figure~\ref{fig:bootstrap_DM_sequences_of_models_C01}, the bootstrap's coverage is much higher than the one of the delta method.
Therefore, the CI provided by the nonparametric percentile bootstrap may be an interesting alternative compared to the delta method when conducting inference with a given sample.
This is all the more so as the mean in the denominator is close to~$0$ (in Figure~\ref{fig:bootstrap_DM_sequences_of_models_C01}, of the size of $n^{-0.25}/10$ for a variance normalized to $1$) and the number of observations is moderately large (a few thousands here).
\section{Construction of nonasymptotic confidence intervals for ratios of expectations}
\label{sec:construct_cis}

To construct nonasymptotic confidence intervals,
we rely on the possibility to ensure that with large probability (i) $\meanX$ is close to $\expecXn$, and (ii) $\meanY$ is both close to $\expecYn$ and bounded away from 0. Under Assumptions~\ref{hyp:basic_dgp_XY} and~\ref{hyp:bt_dgp_XY}, the Bienaym\'e-Chebyshev inequality can be applied to obtain (i) and (ii).
On the other hand, without further restrictions, we are only able to build nonasymptotic CIs at nominal levels that are not too close to~$1$ (see Section~\ref{sub:bt_robust_ci_general_case}).

\medskip

This limitation does not arise with nonasymptotic confidence intervals for expectations.
In that sense, we can say that building  nonasymptotic CIs for ratios of expectations is more demanding.
Intuitively, the extra difficulty of the latter task comes from the need to ensure~(ii). To stress that point, we show in the next subsection that when $\meanY$ is bounded away from~$0$ and positive almost surely, we can build nonasymptotic CIs at every nominal level.

\subsection{An easy case: the support of the denominator is well-separated from~\texorpdfstring{$0$}{0}}
\label{subsec:bt_easy}


We present a simple framework in which it is possible to build nonasymptotic CIs, valid for every $n\in\Nstar$, and with coverage $1-\alpha$ for every $\alpha\in(0,1)$. To do so, we restrict further the set $\Pfam$ of admissible distributions with the following assumption.

\begin{hyp}
    \label{hyp:bt_easy_dgp_well_separated}
    For every~${n \in \Nstar}$, there exists a positive finite constant $\aYn$ such that $Y_{1,n}\geq\aYn$ almost surely.
\end{hyp}

Under Assumption~\ref{hyp:bt_easy_dgp_well_separated}, for every $n \in \Nstar$, $\meanY \geq \aYn > 0$ almost surely under every distribution in $\Pfam$ and $\meanY^{-1}$ is bounded from above.
This assumption obviously rules out binary $\{0,1\}$ random variables in the denominator of the ratio, which can be quite restrictive in practice.
Under this assumption, the following theorem gives a concentration inequality for our ratio of expectations. It is proved in Section~\ref{proof:thm:bt_finite_ci_easy}.






\begin{thm}\label{thm:bt_finite_ci_easy}
    Let Assumptions~\ref{hyp:basic_dgp_XY}, \ref{hyp:bt_dgp_XY} and~\ref{hyp:bt_easy_dgp_well_separated} hold.
    For every $n\in\Nstar$, $\varepsilon > 0$,
    we have
    \begin{align*}
        \sup_{P \in \Pfam} \Prob_{P^{\otimes n}} \Bigg(
        \bigg| \frac{\meanX}{\meanY} - \frac{\expecXn}{\expecYn} \bigg|
        & > \frac{\big(\varepsilon+\sqrt{\upperXn}\big)\varepsilon}{\aYn\lowerYn}
        + \frac{\varepsilon}{\lowerYn} \Bigg)
        \leq \frac{\upperXn}{n\varepsilon^2}+\frac{\upperYn-\lowerYn^2}{n\varepsilon^2}.
    \end{align*}
    As a consequence, $\inf_{P \in \Pfam} \Prob_{P^{\otimes n}} \Big( \expecXn / \, \expecYn \in
    \left[\, \meanX /\,\meanY \pm t \right] \Big) \geq 1 - \alpha$,
    with the choice $$t:= \frac{1}{\lowerYn}\sqrt{\frac{\upperXn+\upperYn-\lowerYn^2}{n\alpha}}\left(1+\frac{1}{\aYn}\left\{\sqrt{\frac{\upperXn+\upperYn-\lowerYn^2}{n\alpha}}+\sqrt{\upperXn}\right\}\right),$$ for every $\alpha \in (0,1)$.
\end{thm}

The theorem shows that it is possible to construct nonasymptotic CIs for ratios of expectations, with guaranteed coverage at every confidence level, that are almost surely bounded under every distribution in $\Pfam$ characterized by Assumptions~\ref{hyp:basic_dgp_XY}, \ref{hyp:bt_dgp_XY} and~\ref{hyp:bt_easy_dgp_well_separated}.
In Section~\ref{sub:bt_robust_ci_general_case}, we give an analogous result that only requires Assumptions~\ref{hyp:basic_dgp_XY} and~\ref{hyp:bt_dgp_XY} to hold, so that it encompasses the case of $\{0,1\}$-valued denominators.
However, the cost to pay will be an upper bound on the achievable coverage of the confidence intervals.

\subsection{General case: no assumption on the support of the denominator}
\label{sub:bt_robust_ci_general_case}

We seek to build nontrivial nonasymptotic CIs under Assumptions~\ref{hyp:basic_dgp_XY} and \ref{hyp:bt_dgp_XY} only.
Under Assumption \ref{hyp:basic_dgp_XY}, ${\expecYn\neq 0}$, so that there is no issue in considering the fraction $\expecXn/\,\expecYn$.
However, without Assumption~\ref{hyp:bt_easy_dgp_well_separated}, $\left\{\meanY=0\right\}$ has positive probability in general so that $\meanX/\,\meanY$ is well-defined with probability less than one.
Note that when $P_{Y,n}$ is continuous with respect to Lebesgue's measure, there is no issue in defining $\meanX/\,\meanY$ anymore since the event $\left\{\meanY = 0\right\}$ has probability zero.
This is not an easier case from a theoretical point of view though since, without more restrictions, $\meanY$ can still be arbitrarily close to~$0$ with positive probability.


\begin{thm}
    \label{thm:concentration_inequality_BC_general_case}
    Let Assumptions~\ref{hyp:basic_dgp_XY} and~\ref{hyp:bt_dgp_XY} hold.
    For every $n\in\Nstar$, $\varepsilon > 0, \tilde \varepsilon \in (0,1)$,
    we have
    \begin{align*}
        \sup_{P \in \Pfam} \Prob_{P^{\otimes n}}  \Bigg( \bigg| \frac{\meanX}{\meanY} - \frac{\expecXn}{\expecYn} \bigg|
        > \bigg( \frac{\big(\sqrt{\upperXn}+\varepsilon\big) \tilde \varepsilon}{(1 - \tilde \varepsilon)^2} + \varepsilon \bigg) \frac{1}{\lowerYn}  \Bigg)
        \leq \frac{\upperXn}{n\varepsilon^2}+\frac{\upperYn-\lowerYn^2}{n\tilde\varepsilon^2\lowerYn^2}.
    \end{align*}
    As a consequence, $\inf_{P \in \Pfam} \Prob_{P^{\otimes n}} \Big( \expecXn / \, \expecYn \in
    \left[\, \meanX /\,\meanY \pm t \right] \Big) \geq 1-\alpha$,
    with the choice
    \begin{align*}
        &t= \frac{1}{\lowerYn}\left(\frac{\left(\sqrt{\upperXn}+\sqrt{2\upperXn/(n\alpha)}\right)\sqrt{2(\upperYn-\lowerYn^2)/(n\alpha\lowerYn^2)}}{\left(1-\sqrt{2(\upperYn-\lowerYn^2)/(n\alpha \lowerYn^2})\right)^2}+\sqrt{\frac{2\upperXn}{n\alpha}}\right),
    \end{align*}
    for every $\alpha > \overline{\alpha}_n := \frac{2(\upperYn-\lowerYn^2)}{n\lowerYn^2}$.\footnote{Equivalently, it means that for a given $\alpha$, the above choice of $t$ is valid for every integer~$n > \overline{n}_\alpha :=2(\upperYn-\lowerYn^2)/(\alpha\lowerYn^2)$.}
    \label{thm:bt_finite_ci}
\end{thm}

This theorem is proved in Section~\ref{proof:thm:concentration_inequality_BC_general_case}.
It states that when $\lowerYn > 0$, it is possible to build valid nonasymptotic CIs with finite length up to the confidence level $1-\overline{\alpha}_n$. This is a more positive result than \cite{dufour1997} which states that it is not possible to build nontrivial nonasymptotic CIs when $\lowerYn$ is taken equal to 0, no matter the confidence level. 
Note that Theorem~\ref{thm:bt_finite_ci} is not an impossibility theorem since it only claims that considering confidence levels smaller than ${1-\overline{\alpha}_n}$ is \textit{sufficient} to build nontrivial CIs under Assumptions~\ref{hyp:basic_dgp_XY} and \ref{hyp:bt_dgp_XY}. The remaining question is to find out whether it is \textit{necessary} to focus on confidence levels that do not exceed a certain threshold under Assumptions~\ref{hyp:basic_dgp_XY} and~\ref{hyp:bt_dgp_XY}.
We answer this in Section~\ref{subsec:bt_upper_bound_level}.

\medskip

Theorem~\ref{thm:bt_finite_ci} has two other interesting consequences: for every confidence level up to $1-\overline{\alpha}_n$, a nonasymptotic interval of the form $\left[\meanX/\,\meanY\pm \tilde{t}\right]$ with $\tilde{t}>t$ has coverage ${1-\alpha}$ but is unnecessarily conservative.
Moreover, if the data generating process does not depend on $n$ (i.e. in the standard i.i.d. set-up), the length of the confidence interval shrinks at the optimal rate $1/\sqrt{n}$ for every fixed $\alpha$.
Note that the coefficient~$2$ in the definition of $\overline{\alpha}_n$ defined above can be reduced to any number ${w > 1}$, at the expense of increasing the length of the confidence interval (this length actually tends to infinity when $w$ tends to~$1$).

\section{Nonasymptotic CIs: impossibility results and practical guidelines}
\label{sec:about_alpha}

In this section, we prove two impossibility results: a maximum confidence level above which it is impossible to build nontrivial nonasymptotic CIs and a necessary lower bound on the length of nonasymptotic CIs.

\subsection{An upper bound on testable confidence levels}
\label{subsec:bt_upper_bound_level}


\begin{prop}
    \label{thm:bt_necessary_cond_alpha_t}
    Let~$\Pfam$ be the class of all distributions satisfying Assumptions~\ref{hyp:basic_dgp_XY} and~\ref{hyp:bt_dgp_XY} and $\underline{\alpha}_n :=  \big(1-l_{Y,n}^2/u_{Y,n}\big)^n$.
	For every $n\in\Nstar$ and every $\alpha \in
	\left(0, \underline{\alpha}_n \right)$,
	if~$\lowerYn^2/\upperYn < 1$, there is no finite $t>0$ such that $\left[\meanX/\,\meanY\pm t\right]$ has coverage $1-\alpha$ over $\Pfam$.
\end{prop}

This theorem asserts that confidence intervals of the form $\left[\meanX/\,\meanY\pm t\right]$ with coverage higher than ${1-\underline{\alpha}_n}$ under Assumptions~\ref{hyp:basic_dgp_XY} and~\ref{hyp:bt_dgp_XY} are not defined (or are of infinite length) with positive probability for at least one distribution in~$\Pfam$.
This is due to the fact that $\underline{\alpha}_n$ is a lower bound on $\Prob(\meanY = 0)$ over all distributions in~$\Pfam$.

\medskip

Remark that when $\upperYn/\lowerYn^2=1$, there is no impossibility result anymore: assume that $\upperYn/\lowerYn^2=1$ and let $Q$ be a distribution on $\Rb^2$ that satisfies Assumptions~\ref{hyp:basic_dgp_XY} and \ref{hyp:bt_dgp_XY}.
Let $(X_{i,n},Y_{i,n})_{i=1}^n\simiid Q$.
We have that $\Var[Y_{1,n}]=0$, which implies that $Y_{1,n}=\expecYn$ almost surely.
Assumption~\ref{hyp:basic_dgp_XY} further ensures that $Y_{1,n}\neq 0$ almost surely. Consequently, the results of Section~\ref{subsec:bt_easy} apply and allow us to conclude that under Assumptions~\ref{hyp:basic_dgp_XY}, \ref{hyp:bt_dgp_XY} and $\upperYn/\lowerYn^2=1$, it is possible to build nontrivial nonasymptotic CIs at every confidence level.
Indeed, in that case, we are in fact only estimating a simple mean, and therefore there is no constraint on~$\alpha$.

\medskip

Proposition~\ref{thm:bt_necessary_cond_alpha_t} is actually a corollary of the more general Theorem~\ref{thm:bt_necessary_cond_alpha}.
It states it is impossible to construct confidence intervals that contain $\meanX/\,\meanY$ almost surely and are almost surely bounded over $\Pfam$ with coverage greater than ${1 - \underline{\alpha}_n}$.
It is proved in Section~\ref{proof:thm:bt_necessary_cond_alpha}.

\begin{thm}
    \label{thm:bt_necessary_cond_alpha}
    Let~$\Pfam$ be the class of all distributions satisfying Assumptions~\ref{hyp:basic_dgp_XY} and~\ref{hyp:bt_dgp_XY}.
    Let $n\in\INTST$, and a random set $I_n$ that contains $\meanX/\,\meanY$ almost surely whenever it is defined and is undefined if $\meanY=0$.
    Then $\sup_{P \in \Pfam}
    \Prob_{P^{\otimes n}} \big( I_n \, undefined \big)
    \geq \underline{\alpha}_n.$
\end{thm}

Combining Theorems~\ref{thm:bt_finite_ci} and~\ref{thm:bt_necessary_cond_alpha}, we conclude that there exists some critical level $1-\alpha_n^c$ belonging to the interval $[1-\overline{\alpha}_n,1-\underline{\alpha}_n]$ such that it is impossible to build nontrivial nonasymptotic confidence intervals if and only if their nominal level is above $1-\alpha_n^c$. Finally, it is worth remarking that with a sample of size $n$, the CIs based on the delta method with a nominal level $1-\alpha>1-\alpha_n^c$ cannot have coverage~${1-\alpha}$ uniformly over~$\Pfam$ as such CIs verify the condition of Theorem~\ref{thm:bt_necessary_cond_alpha}.

\medskip

Figure~\ref{fig:schema_comparison_levels} below shows the critical level and its bounds obtained in our nonasymptotic results.

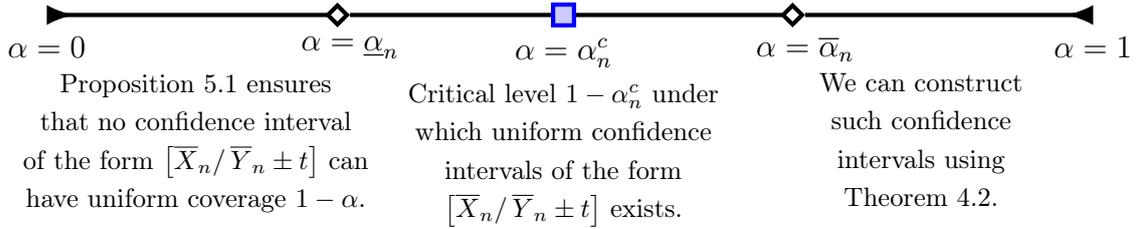
\begin{figure}[ht]
    \vspace{1.5em}
    \centering
    \begin{tikzpicture}
        \node[] (alpha_0) {} ;
        \node[right = of alpha_0, xshift = 2.7cm] (alpha_nec) {} ;
        \node[right = of alpha_nec, xshift = 4.7cm] (alpha_sup) {} ;
        \node[right = of alpha_sup, xshift = 2.5cm] (alpha_1) {} ;

        \begin{scope}[>=Latex, ultra thick]
        \draw[>-{Turned Square[open]}]
        (node cs:name=alpha_0, anchor=center)
        to (node cs:name=alpha_nec, anchor=center) ;
        \draw[-{Turned Square[open]}]
        (node cs:name=alpha_nec, anchor=center)
        to (node cs:name=alpha_sup, anchor=center) ;
        \draw[-<]
        (node cs:name=alpha_sup, anchor=center)
        to (node cs:name=alpha_1, anchor=center) ;
        \end{scope}

        \node[below = of alpha_0, yshift = 1cm] {$\alpha = 0$} ;
        \node[below = of alpha_nec, yshift = 1cm] {$\alpha = \underline{\alpha}_n$} ;
        \node[below = of alpha_sup, yshift = 1cm] {$\alpha = \overline{\alpha}_n$} ;
        \node[below = of alpha_1, yshift = 1cm] {$\alpha = 1$} ;

        \node[right = of alpha_0, yshift = -1.7cm, xshift = -1.6cm,
        align = center]
        {\footnotesize Proposition~\ref{thm:bt_necessary_cond_alpha_t} ensures \\
        \footnotesize that no confidence interval \\
        \footnotesize of the form $\left[\meanX/\,\meanY\pm t\right]$ can \\
        \footnotesize have uniform coverage $1-\alpha$.} ;

        \node[right = of alpha_nec, xshift = 1.5cm, fill = blue!20, draw = blue, ultra thick] (alpha_c) {} ;

        \node[below = of alpha_c, yshift = 1cm] {$\alpha = {\alpha}_n^c$} ;

        \node[below = of alpha_c, yshift = 0.4cm, xshift = 0cm,
        align = center]
        {\footnotesize Critical level $1-{\alpha}_n^c$ under \\
        \footnotesize which uniform confidence \\
        \footnotesize intervals of the form \\
        \footnotesize  $\left[\meanX/\,\meanY\pm t\right]$ exists.} ;

        \node[right = of alpha_sup, yshift = -1.65cm, xshift = -1.1cm,
        align = center]
        {\footnotesize We can construct \\
        \footnotesize such confidence \\
        \footnotesize intervals using \\
        \footnotesize Theorem~\ref{thm:concentration_inequality_BC_general_case}.} ;
    \end{tikzpicture}
    \caption{The critical level and its bounds.}
    \label{fig:schema_comparison_levels}
\end{figure}

\begin{rem}
\label{rem:snrtilde_limit}
    In the same spirit as in Theorem~\ref{thm:delta_method_sequence_models}, we consider a modified version of the signal-to-noise ratio defined by $\SNRntilde := l_{Y,n}/(u_{Y,n}^{1/2} n^{-1/2})$.
    When $\SNRntilde \to + \infty$ $($resp. $0)$ as $n \to \infty$, $\underline{\alpha}_n$ and $\overline{\alpha}_n$ tend to~$0$ $($resp. $+\infty)$.
    When we have enough information ${(\SNRntilde \to + \infty)}$, the critical level~${1-\alpha_n^c}$ tends to~$1$.
    Therefore, for every~${\alpha \in (0,1)}$, nonasymptotic confidence intervals can be constructed at every level for $n$ large enough.
    On the contrary, when ${\SNRntilde \to 0}$, the critical level~${1-\alpha_n^c}$ tends to~$0$, which means that it is impossible to construct uniformly valid CIs for $n$ large enough.
    Finally, when ${\SNRntilde \to C}$ for a positive constant~$C$, a critical level remains as in the nonasymptotic case since $\underline{\alpha}_n \to \exp(-C)$.
\end{rem}

\subsection{A lower bound on the length of nonasymptotic confidence intervals}
\label{subsec:bt_lower_bound_length}

The following theorem is an extension of \cite{catoni2012challenging}[Proposition 6.2] to ratios.
It is proved in Section~\ref{proof:thm:bt_l_bound_ci_length}.

\begin{thm}\label{thm:bt_l_bound_ci_length}
    For every integer $n\geq 7$,
    $\alpha\in \big( 0, 1\,\wedge\, n/\big(\lowerYn+\sqrt{\upperYn-\lowerYn^2}\big)^2\big)$,
    and $\xi < 1$ there exists a distribution $Q$ on $\Rb^2$ that satisfies Assumptions~\ref{hyp:basic_dgp_XY} and \ref{hyp:bt_dgp_XY} such that for $\left(X_{i,n},Y_{i,n}\right)_{i=1}^n\overset{i.i.d}\sim~Q$, we have
    \begin{align*}
        &\Prob_{Q^{\otimes n}} \left(\left|
        \frac{\meanX}{\meanY} - \frac{\expecXn}{\expecYn} \right|
        > \xi \sqrt{\frac{v_n}{3n\alpha}} \right) >\alpha,
    \end{align*}
    where $v_n:= \upperXn / \big(\lowerYn+\sqrt{\upperYn-\lowerYn^2}\big)^2$.
\end{thm}

With this theorem, we can claim that
CIs of the form $\left[\meanX/\,\meanY\pm t\right]$ cannot have uniform coverage $1-\alpha$,
for every $\alpha\in\big(0,1\wedge n/\big(\lowerYn+\sqrt{\upperYn-\lowerYn^2}\big)^2\big)$,
under Assumptions~\ref{hyp:basic_dgp_XY} and~\ref{hyp:bt_dgp_XY} if they are shorter than $\sqrt{v_n / (3n\alpha)}$.
%
%
By a careful inspection of the proof (see Lemma~\ref{lem:tech_lemma_lower_bounds}), we can in fact replace the value $3$ in the theorem by any number strictly larger than $e=\exp(1)$, at the price of assuming $n \geq n_0$ for $n_0$ large enough. It is interesting to note that the distributions $Q$ that are built in the proof of the theorem are on the boundary of~$\Pfam$ in the sense that they satisfy $\expecXnsq=\upperXn$, $\expecYn=\lowerYn$ and $\expecYnsq=\upperYn$.

\subsection{Practical methods and plug-in estimators}
\label{subsec:plug_in_estimators}

Nonasymptotic confidence intervals and the thresholds $\overline{\alpha}_n$ and $\overline{n}_\alpha$ based on Theorem~\ref{thm:concentration_inequality_BC_general_case} rely on Assumptions~\ref{hyp:basic_dgp_XY} and~\ref{hyp:bt_dgp_XY}.
In practice, building such CIs or computing those thresholds require the knowledge of the constants $\lowerYn$, $\upperXn$ and $\upperYn$ that determine the class of distributions we consider.\footnote{Actually, the computation of $\overline{\alpha}_n$ and $\overline{n}_\alpha$ only require the knowledge of $\lowerYn$ and $\upperYn$.}
Therefore, we need to state some values for those constants.
Note that constructing nontrivial and nonasymptotic CIs that overcome the limitations of having to choose some a priori class of distributions is not possible.
Indeed, we would get back to \cite{bahadur1956} and \cite{dufour1997} type impossibility results.

\medskip

How to choose $\lowerYn$, $\upperXn$ and $\upperYn$ depends on the specific application.
Sometimes, stating values can be sensible if researchers do have control or expert knowledge of the variables.
Resuming an example started in the introduction, if the variable in the denominator is an indicator of being treated
in the setting of a Randomized Controlled Trial, researchers can have intuitions about reasonable values for the lower and upper bounds of the probability of being treated.

\medskip

The unknown constants are upper and lower bounds on moments that characterize the class~$\Pfam$.
As such, they can never be recovered from the data since observations are by construction drawn from a single distribution~${P \in \Pfam}$.
Under i.i.d. sampling, sample means converge to their corresponding theoretical moments, provided the latter are finite.
Hence, without prior information, a plug-in strategy has to be used which consists in: (i) using the moments of a single distribution instead of the bounds on the class, (ii) estimating those moments with their empirical counterparts.
As a consequence, this approach is valid pointwise only and not uniformly over~$\Pfam$ anymore.
Furthermore, it is only asymptotically justified.
On the other hand, for any sample provided~${\meanY \neq 0}$, this plug-in strategy enables us to construct our CIs and the quantity $\overline{n}_\alpha$ (or $\overline{\alpha}_n$), which can be a useful rule of thumb as explained below.
We stick to that principle in our simulations and application.

\medskip

For a given level~${1-\alpha}$ and a class of distributions satisfying Assumptions~\ref{hyp:basic_dgp_XY} and~\ref{hyp:bt_dgp_XY},
$\overline{n}_\alpha$ is the minimal sample size required to construct our nonasymptotic CIs.
In other words, for a sample size~${n < \overline{n}_\alpha}$, the data is not rich enough to construct the nonasymptotic CIs of Theorem~\ref{thm:bt_finite_ci} at this level.
Heuristically, the comparison of $\overline{n}_\alpha$ and $n$ can be used as a rule of thumb to assess whether the coverage of the CIs based on the delta method matches their nominal level.\footnote{Equivalently, we could compare $\overline{\alpha}_n$ and $\alpha$.
As a rule of thumb, $\overline{\alpha}_n$ can be seen as the lowest $\alpha$ (hence the highest nominal level $1-\alpha$) for which the asymptotic CIs based on the delta method are reliable given the sample size~$n$.}
Several simulations tend to confirm the practical interest of that rule of thumb as $\overline{n}_\alpha$ turns out to be very close to the sample size above which the gap between the coverage of the asymptotic CIs based on the delta method and their nominal level becomes negligible.
(see Section~\ref{sub:simulations} and Appendix~\ref{appendix:sec:additional_simulations}).

\section{Numerical applications}
\label{sec:numerical_app}

\subsection{Simulations}
\label{sub:simulations}

This section presents simulations that support the use of $\overline{n}_\alpha$, or equivalently $\overline{\alpha}_n$, as a rule of thumb to inspect the reliability of the asymptotic confidence intervals from the delta method.

\medskip

\begin{figure}[htbp]
    \centering
    \includegraphics[width=0.75\textwidth]{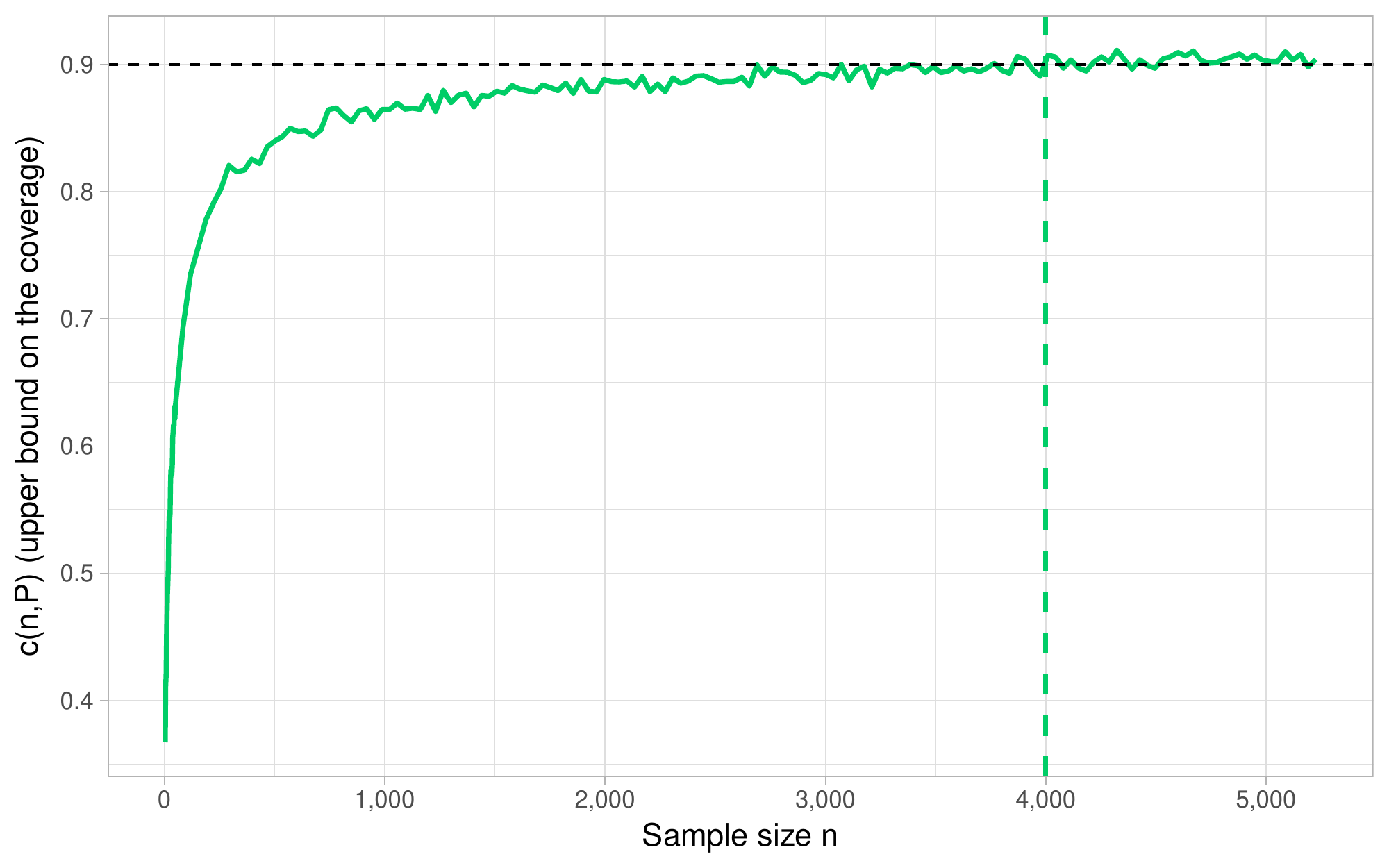}
    \caption{\small{$c(n,P)$ of the asymptotic CIs based on the delta method as a function of the sample size~$n$ and~$\overline{n}_\alpha$.\newline
    Specification: $\forall n \in \Nstar, {\distributionXYn = \mathcal{N}_2}$ (bivariate Gaussian) with $\expecX = 0.5$, $\expecY = 0.1$, $\VarX = 1$, $\VarY = 2$, $\Corr(X,Y) = 0.5$.
    The nominal pointwise asymptotic level is set to $0.90$.
    For a sample size $n$, the coverage is obtained as the mean over 5,000 repetitions.
    The dashed vertical line shows $\overline{n}_\alpha :=2\left(\upperYn-{\lowerYn}^2\right)/\big(\alpha\lowerYn^2\big)$, setting here $\alpha = 0.1$, $\lowerYn = \expecY$, $\upperYn = {\expecY}^2 + \VarY$.}}
    \label{fig:gr_newa_gaussians_n_alphapc_10nbrep_5000EYpc_10EXpc_50VXpc_100VYpc_200CorrXYpc_50dilanpc_130}
\end{figure}

\begin{figure}[htbp]
    \centering
    \includegraphics[width=0.75\textwidth]{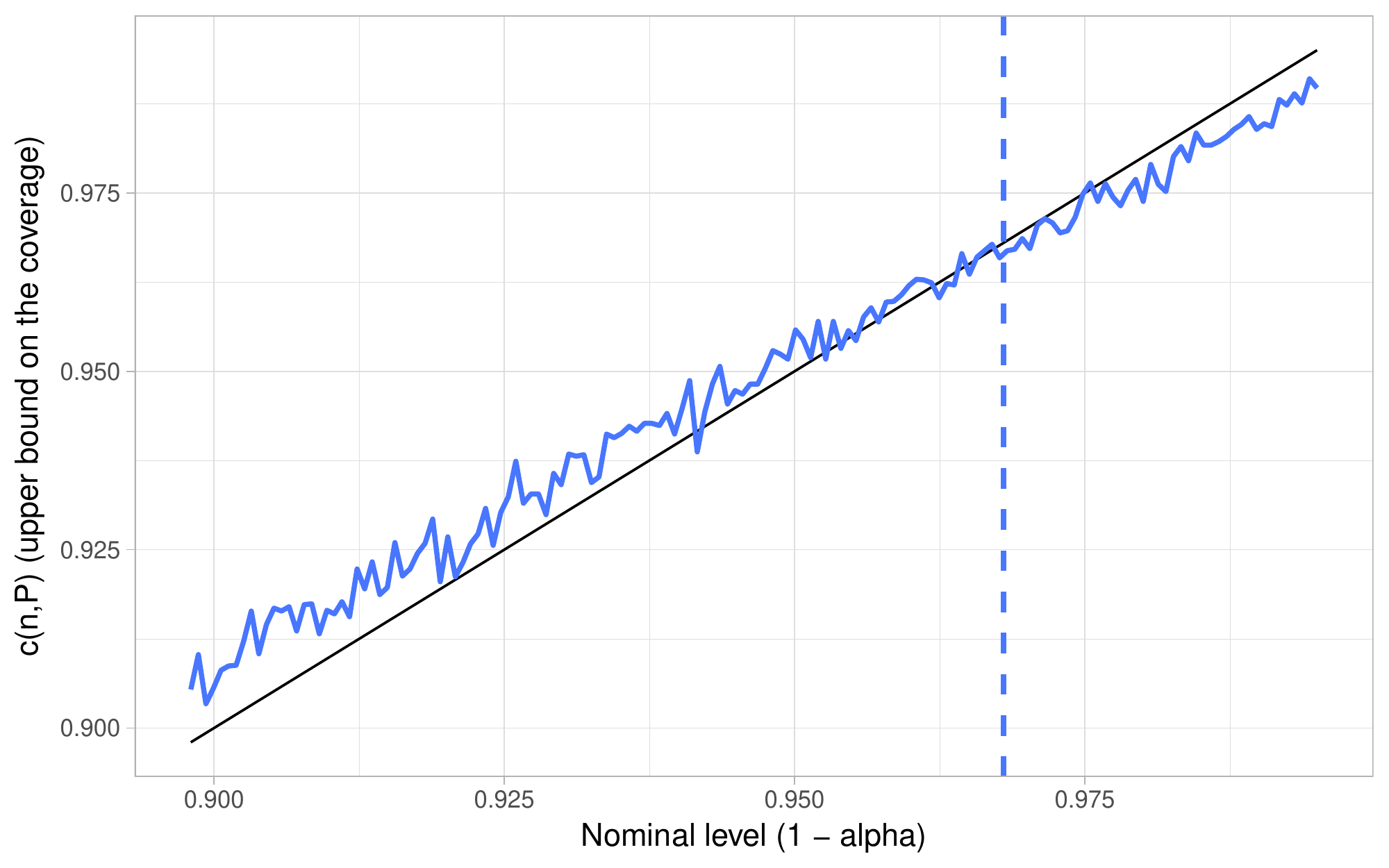}
    \caption{\small{$c(n,P)$ of the asymptotic CIs based on the delta method as a function of the sample size~$n$ and~$\overline{\alpha}_n$.\newline
    Specification: ${\forall n \in \Nstar, \distributionXYn = \mathcal{N}_2}$ (bivariate Gaussian) with $\expecX = 0.5$, $\expecY = 0.25$, $\VarX = 2$, $\VarY = 1$, $\Corr(X,Y) = 0.5$.
    The sample size is~$n =$ 1,000.
    For each nominal level $1-\alpha$ in the x-axis, we draw $10,000$ samples, compute the asymptotic CIs and see whether it covers or not the ratio of interest; we report the mean over the $10,000$ repetitions in the y-axis.
    The solid line is the first bisector $y=x$.
    The dashed vertical line shows ${\overline{\alpha}_n :=2\left(\upperYn-{\lowerYn}^2\right)/\big(n{\lowerYn}^2\big)}$, setting here ${\lowerYn = \expecY}$, ${\upperYn = {\expecY}^2 + \VarY}$.}}
    \label{fig:gr_gr_newa_gaussian_alpha_n_nbrep10000adila_3n_1000EY_025EX_05VX_2VY_1CorrXY_05}
\end{figure}

In Figure~\ref{fig:gr_newa_gaussians_n_alphapc_10nbrep_5000EYpc_10EXpc_50VXpc_100VYpc_200CorrXYpc_50dilanpc_130}, a nominal level ${1-\alpha}$ is fixed and we show the $c(n,P)$ of the CIs based on the delta method as a function of the sample size~$n$, as well as $\overline{n}_\alpha$ derived in Theorem~\ref{thm:concentration_inequality_BC_general_case}.
It happens that the coverage converges toward its nominal level for sample sizes around $\overline{n}_\alpha$, which supports $\overline{n}_\alpha$ as a rule of thumb of interest in practice.\footnote{This fact holds across various specifications (see additional simulations in Appendix~\ref{appendix:sec:additional_simulations}).}
%
%
In Figure~\ref{fig:gr_gr_newa_gaussian_alpha_n_nbrep10000adila_3n_1000EY_025EX_05VX_2VY_1CorrXY_05}, a sample size is fixed and we show the coverage
for different nominal levels, as well as the quantity $\overline{\alpha}_n$.
It is the converse of Figure~\ref{fig:gr_newa_gaussians_n_alphapc_10nbrep_5000EYpc_10EXpc_50VXpc_100VYpc_200CorrXYpc_50dilanpc_130} in that sense.
In this simulation,
$\overline{\alpha}_n$ turns out to fall close to the lowest $\alpha$ (hence highest $1-\alpha$) for which the coverage of the CIs based on the delta method attains their nominal level.

\medskip

All in all, Figures~\ref{fig:gr_newa_gaussians_n_alphapc_10nbrep_5000EYpc_10EXpc_50VXpc_100VYpc_200CorrXYpc_50dilanpc_130} and~\ref{fig:gr_gr_newa_gaussian_alpha_n_nbrep10000adila_3n_1000EY_025EX_05VX_2VY_1CorrXY_05} and additional simulations advocate the use of $\overline{n}_\alpha$ derived in Theorem~\ref{thm:concentration_inequality_BC_general_case} (or conversely $\overline{\alpha}_n$) as a rule of thumb to appraise the dependability of the CIs obtained with the delta method for ratios of expectations.

\subsection{Application to real data}
\label{sub:application_real_data}

We illustrate our methods with an application related to gender wage disparities.
The application resumes our canonical example of conditional expectations since we estimate the proportion of women within wage brackets that are defined as having a wage higher than a given threshold.
We use $n =$ 204,246 observations from the French Labor Survey data between 2010 and 2017.\footnote{Enquête Emploi en continu (version FPR) -- 2010-2017, INSEE [producteur], ADISP [diffuseur].}

\medskip

Let~$W$ be a real random variable that indicates the wage of an employee (expressed in euros per month) and $F$~an indicator variable equal to~$1$ if the employee is a woman and~$0$ otherwise.
For a given threshold wage~$w_0$, the parameter of interest is ${\Expec[F \mid W \geq w_0]}$.
It can be written as a ratio of expectations with $X = F\,\Indicator\{W \geq w_0\} = \Indicator\{F = 1, W \geq w_0\}$ in the numerator and $Y = \Indicator\{W \geq w_0\}$ in the denominator.
As we consider higher thresholds~$w_0$, the expectation in the denominator gets closer to~$0$.
As an illustration, out of $n=$ 204,246 observations, $355$ individuals have monthly wages higher than 10,000 euros (which corresponds to a mean in the denominator equal to $0.0017$); $44$ individuals above 20,000 ($\meanY = 2.2 \times 10^{-4}$); and only $17$ above 30,000 ($\meanY = 8.3 \times 10^{-5}$).\footnote{To give a sense of the wage distribution, note that the empirical quantiles of~$W$ at orders 90\%; 95\%; 99\%;
and 99.99\% are respectively: 2,989; 3,728; 6,000;
and 26,024.}

\begin{figure}[htb]
    \centering
    \includegraphics[width=0.9\textwidth]{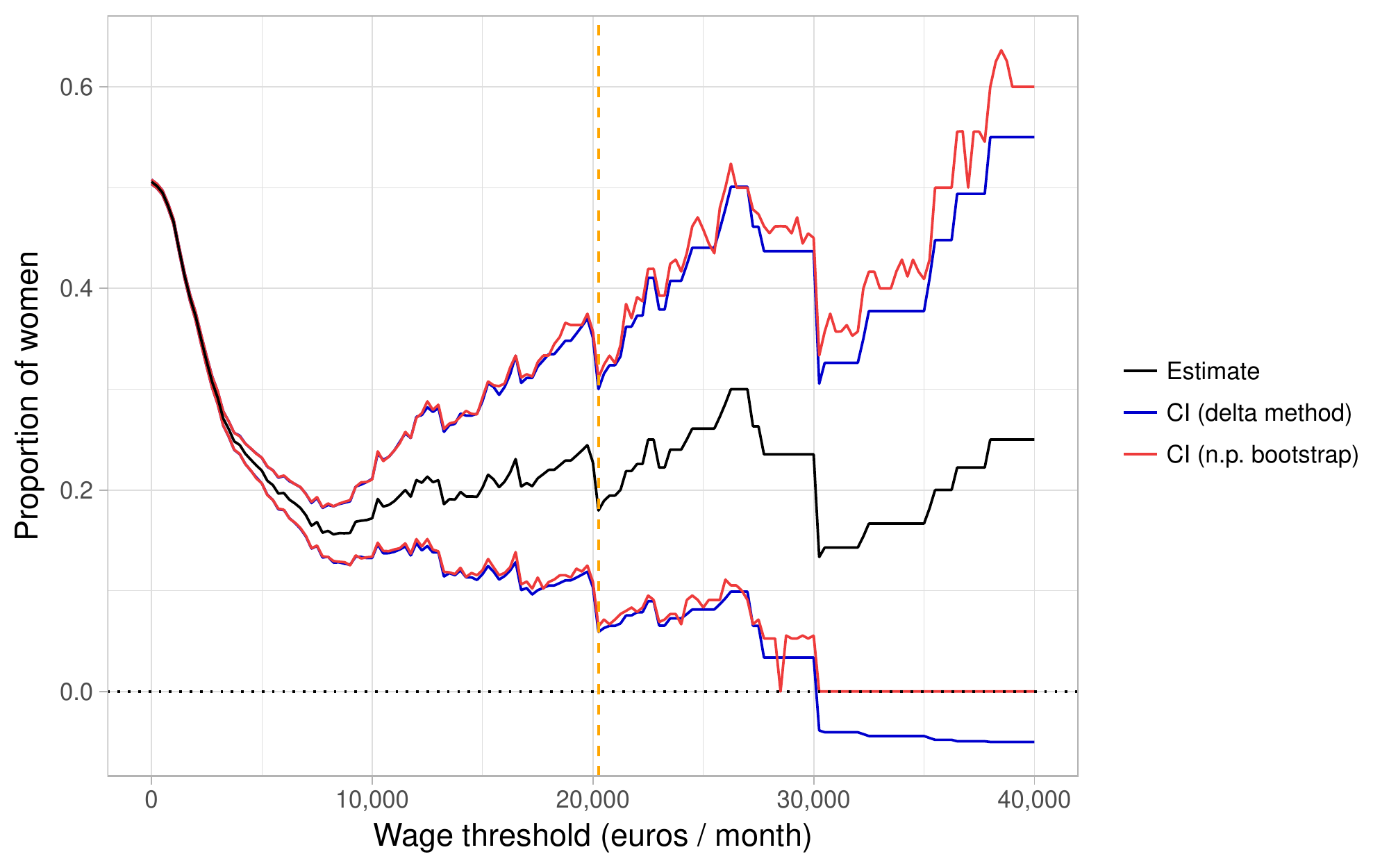}
    \caption{\small{Point estimate and confidence intervals for the parameter~${\Expec[F \mid W \geq w_0]}$ as a function of the wage threshold~$w_0$.
    The parameter is the proportion of women within the wage bracket~${[w_0, +\infty)}$.
    The nominal level of the CIs is set to~95\%.
    Efron's percentile bootstrap CIs are obtained using 2,000 bootstrap replications.
    The dashed vertical line represents the lowest wage threshold such that the plug-in counterpart of ${\overline{n}_\alpha}$ exceeds~$n$.}}
    \label{fig:application_prop_women_threshold_wage}
\end{figure}


For various thresholds~$w_0$, Figure~\ref{fig:application_prop_women_threshold_wage} presents the estimate~$\thetahatn$ and two 95\%-nominal-level confidence intervals for the parameter~${\Expec[F \mid W \geq w_0]}$: the one based on the delta method (see Section~\ref{sub:problem_slow_convergence}) and the one using Efron's percentile bootstrap (see Section~\ref{sub:bootstrap_CI}).
With higher thresholds, the expectation in the denominator is closer to~$0$ which results in wider confidence intervals.
For very high thresholds, the CIs become hardly informative.
In particular, the lower end of the interval based on the delta method is negative whereas the parameter of interest belongs to~$[0,1]$ by construction.

\medskip

The dashed vertical line relates to our rule of thumb introduced in Section~\ref{subsec:plug_in_estimators}.
More precisely, given the level~${1-\alpha=0.95}$, for each threshold~$w_0$, we compute the plug-in counterpart of~$\overline{n}_\alpha$ defined in Theorem~\ref{thm:concentration_inequality_BC_general_case}:
$2\left(n^{-1}\sum_{i=1}^n Y_i^2 - \meanY^{2}\right) /
\big(\alpha \meanY^{2} \big)$.
Given that~$Y$ is a binary variable, the latter quantity
is increasing with~$w_0$ and exceeds $n$ at some threshold represented by the dashed vertical line (here a little above~20,000).
Consequently, for higher thresholds, our rule of thumb suggests that the confidence intervals obtained with the delta method might undercover as the expectation in the denominator is ``too close to~$0$'' relative to the number of observations.
Actually, in the application, it is around this vertical line that the two CIs start to differ.
In particular, the upper end of Efron's percentile confidence interval becomes larger than the upper end of the interval based on the delta method.

\section{Conclusion}
\label{sec:conclusion}

This paper studies the construction of confidence intervals for ratios of expectations, which are frequent parameters of interest in applied econometrics.

\medskip

The most common method to do so is asymptotic and yields CIs based on the asymptotic normality of the empirical means that estimate the numerator and the denominator combined with the delta method.
We document on simulations that the coverage of the confidence intervals based on the delta method may fall short of their nominal level when the expectation in the denominator is close to~$0$, even with fairly large sample size.

\medskip

To further study the reliability of those CIs, we use a sequence-of-model framework, analogous to what a strand of the weak IV literature does.
Indeed, it enables to consider limiting cases, namely here 
denominators tending to~$0$.
In the weak IV case, the equivalent is to move closer to a null covariance between the endogenous regressor and the instrument.
At the limit, the coefficient of interest is not identified.
Our problem differs since the parameter is not even defined in the problematic case of a null denominator.
This issue underlies the impossibility type results presented in the paper.

\medskip

First, in an asymptotic perspective, the possibility of a denominator arbitrarily close to~$0$ explains why we need a sufficiently slow rate of convergence of the expectation in the denominator to~$0$ to conduct meaningful inference.
More precisely, our main asymptotic results basically show that the CIs based on the delta method are valid, as well as those obtained by Efron's percentile bootstrap, when this speed is lower than~$1/\sqrt{n}$ (the standard speed of the CLT).
Furthermore, on simulations, Efron's percentile bootstrap CIs reach their nominal level sooner (namely for smaller sample sizes) than the CIs based on the delta method.
It suggests that beyond the sequence-of-model rationalization, when confronted in practice to a mean in the denominator close to~$0$ relative to the size of the sample at hand, Efron's percentile bootstrap CIs may be more trustworthy than the delta method's ones.

\medskip

Obviously, those cases where the coverage of the CIs based on the delta method can be well below their nominal level do not self-signal to practitioners.
This is why the second part of the paper proposes a rule of thumb to detect those cases and thus assess the dependability of the asymptotic CIs based on the delta method on finite samples.
This index is based on the construction of nonasymptotic confidence intervals and on impossibility results that stem from the problematic null denominator case.

\medskip

In substance, even if we bound away from~$0$ the expectation in the denominator, there remains a partial impossibility result.
Indeed, we show that there exists a critical nominal level above which the coverage of any nonasymptotic confidence interval that is undefined when $\meanY = 0$ cannot uniformly attain its target level.
More precisely, we derive explicit upper and lower bounds on this critical level as a function of the characteristics of the considered class of distributions.
Then, the heuristic of our rule of thumb consists in estimating by plug-in a lower bound on this critical level (or equivalently, for a given level, an upper bound on the minimal required sample size).
The resulting index can thus be computed immediately on any sample.
In addition to its theoretical foundations, various simulations and an application to real data attest the practical usefulness of this rule of thumb.

\medskip

This paper can be seen as a first step towards nonasymptotic inference in econometric models where the issue of close-to-zero denominators arises.
Notable examples may include weak IV, Wald ratios, and difference-in-difference estimands.

\bigskip
\bibliographystyle{abbrv}
\bibliography{references}

\newpage
\appendix

\section{General definitions about confidence intervals}
\label{appendix:sec:definitions}

A standard situation in statistics or econometrics can be modelled as the observation of a sample of $n \in \Nstar$ i.i.d. observations valued in some measurable space $\left(\Zc,\Bz\right)$.
The statistical model is therefore $\left(\Zc,\Bz,\Pfam\right)^{\otimes n}$ with $\Pfam$ some specified set of distributions on $\left(\Zc,\Bz\right)$.
For every distribution $P \in \Pfam$, let $\theta(P)$ be a parameter of interest and the map $\theta:P\mapsto\theta(P)$ be valued in a metric space $\left(\Theta,d\right)$.

\medskip

We denote by $C_n$ a confidence set for $\theta(P)$. 
Formally, a confidence set $C_n$ can be defined as a measurable map from $\left(\Zc,\Bz\right)^{\otimes n}$ to the measurable space $\left(\Csetstheta
\sqcup \{ \text{undefined} \}, \BCsetstheta \sqcup \{ \text{undefined} \} \right)$, where $\Csetstheta$ is the family of all closed subsets of $\Theta$ and $\BCsetstheta$ is the sigma-algebra generated by $\left\{F\in\Csetstheta : F\cap K\neq\emptyset\right\}$ for $K$ running through the family of compact subsets of~$\Theta$.

\medskip

As the vocabulary may somewhat fluctuate between authors, 
we define below classical objects to fix the notations and terminology used in this paper. The goal is to build confidence sets for a targeted \textit{confidence level} $1-\alpha$ (also termed \textit{nominal level} of the confidence set).
For $n \in \Nstar$, for $\alpha \in (0,1)$, we say that a confidence set $C_n$ or a sequence of sets $(C_{n})_{n \in \Nstar}$ has:
\begin{enumerate}[i.]
    \item \emph{coverage} $1-\alpha$ over $\Pfam$ if:
    \vspace{-0.4cm}
    $$\inf_{P\in\Pfam}\Prob_{P^{\otimes n}}\left(C_n\ni\theta\left(P\right)\right)\geq 1-\alpha$$
    \vspace{-1.1cm}
    
    \item \emph{size} $1-\alpha$ over $\Pfam$ if the inequality is an equality:
    \vspace{-0.4cm}
    $$\inf_{P\in\Pfam}\Prob_{P^{\otimes n}}\left(C_n\ni\theta\left(P\right)\right) = 1-\alpha.$$
    \vspace{-1.1cm}
    
    \item \emph{asymptotic coverage $1-\alpha$ pointwise} over $\Pfam$ if:\footnote{Respectively \emph{pointwise asymptotic size} when the inequality is replaced by an equality.}
    \vspace{-0.4cm}
    $$\forall P\in\Pfam, \, \liminf_{n\to+\infty}\Prob_{P^{\otimes n}}\left(C_n\ni\theta\left(P\right)\right)\geq 1-\alpha.$$
    \vspace{-1.1cm}
    
    \item \emph{asymptotic coverage $1-\alpha$ uniformly} over $\Pfam$ if:\footnote{Respectively \emph{uniform asymptotic size} when the inequality is replaced by an equality.}
    \vspace{-0.4cm}
    $$\liminf_{n\to+\infty}\inf_{P\in\Pfam}\Prob_{P^{\otimes n}}\left(C_n\ni\theta\left(P\right)\right)\geq 1-\alpha.$$ 
    \vspace{-1.1cm}
\end{enumerate}

A confidence set with coverage $1-\alpha$ but size different from $1-\alpha$ over $\Pfam$ is said to be \emph{conservative} over $\Pfam$\footnote{Similarly, a confidence set is said to be \emph{asymptotically conservative pointwise} over $\Pfam$ (respectively \emph{uniformly} over $\Pfam$) if property iii. (resp. property iv.) holds with a strict inequality.}.
We further define a \textit{nontrivial confidence set} as a confidence set that is almost surely strictly included in $\Theta$ (whenever it is defined) under every distribution in $\Pfam$. For instance, if $\theta(P)$ is the expectation under $P$, $\Theta=\mathbb{R}$ and $\Pfam$ is the set of all distributions that admit a finite expectation, a nontrivial CI is any CI that is almost surely bounded under every distribution in $\Pfam$.
For ratios of expectations, ${\Theta = \Rb}$ too and we will use the term \emph{almost surely bounded} as a synonym of nontrivial, without stating ``under every distribution in~$\Pfam$'' when there is no ambiguity as regards the class $\Pfam$ considered.

\medskip

A family of confidence intervals $(C_{n, \alpha})_{n \in \Nstar, \, \alpha \in (0,1)}$ is said to be \emph{pointwise \emph{(resp.} uniformly\emph{)} consistent} if for every $\alpha \in (0,1)$, the sequence $(C_{n, \alpha})_{n \in \Nstar}$ has pointwise (resp. uniformly) asymptotic coverage at level $1 - \alpha$.

\section{Proofs of the results in Sections~\ref{sec:about_delta_method}, \ref{sec:construct_cis} and~\ref{sec:about_alpha}}
\label{appendix:sec:proofs_BC_case}

\subsection{Proof of Theorem~\ref{thm:delta_method_sequence_models}}
\label{proof:thm:delta_method_sequence_models}

Let $\theta_{X,n} := \expecXn$, $\theta_{Y,n} := \expecYn$.
Let $h_{X,n} := \sqrt{n} \gamma_{X,n} (\meanX - \expecXn)$ and $h_{Y,n} := \sqrt{n} \gamma_{Y,n} (\meanY - \expecYn)$ be the centered and normalized versions of $\meanX$ and $\meanY$.
We first rewrite Theorem~\ref{thm:delta_method_sequence_models} using this notation.
\begin{thm}
    Let Assumption~\ref{hyp:basic_dgp_XY} hold.
    Assume that
    $\Var[(\gamma_{X,n} X_{1,n} \, , \, \gamma_{Y,n} Y_{1,n})] \to V$ for some positive sequences $\gamma_{X,n}$ and $\gamma_{Y,n}$ where $V$ is a definite positive $2 \times 2$ matrix, that $\Prob(\meanY = 0) \to 0$, as $n \to \infty$ and that

    Then the sequence of random variables $A_n := \meanX/\meanY - \theta_{X,n}/\theta_{Y,n}$ satisfies as ${n \to \infty}$:
    \begin{enumerate}
        \item If $n^{-1/2} = o(\gamma_{Y,n} \theta_{Y,n})$, then $A_n$ is equivalent to
        $$n^{-1/2} \left( \frac{h_{X,n}}{\theta_{Y,n} \gamma_{X,n}} - \frac{h_{Y,n} \theta_{X,n}}{\gamma_{Y,n} \theta_{Y,n}^2} \right).$$

        \item If there exists a finite constant $C \neq 0$ such that $\sqrt{n} \gamma_{Y,n} \theta_{Y,n} \to C$ as $n \to \infty$, then $A_n$ is equivalent to
        $$\sqrt{n} \gamma_{Y,n} \theta_{X,n} \left( \frac{1}{C + h_{Y,n}} - \frac{1}{C} \right) + \frac{h_{X,n} \gamma_{Y,n}}{(C + h_{Y,n}) \gamma_{X,n}}.$$

        \item If $\gamma_{Y,n} \theta_{Y,n} = o(n^{-1/2})$, then $A_n$ is equivalent to
        $$ \frac{h_{X,n} \gamma_{Y,n}}{h_{Y,n} \gamma_{X,n}} - \frac{\theta_{X,n}}{\theta_{Y,n}}.$$
    \end{enumerate}
    \label{thm:delta_method_sequence_models_initial}
\end{thm}

Let us define $W_n :=
\Indicator\{ \theta_{Y,n} + h_{Y,n}/ (\sqrt{n} \gamma_{Y,n}) = 0 \}$ and remark that $W_n = 1$ whenever $\meanY = 0$. By assumption $\Prob(\meanY = 0) \to 0$, therefore $W_n \convD \delta_0$.
Moreover, by Lyapunov's central limit theorem applied to
\begin{equation*}
    (h_{X,n} , h_{Y,n})
    = \sqrt{n} \Big( \frac{1}{n}
    \sum_{i=1}^n (X_{i,n} \gamma_{X,n}, Y_{i,n} \gamma_{Y,n})
    - (\expecX \gamma_{X,n}, \expecY \gamma_{Y,n}) \Big),
\end{equation*}
using $V \neq \mathbf{0}$ and the boundedness of
$\Expec \big[|X_{1,n}|^3 \big] \gamma_{X,n}^3$ and
$\Expec \big[|Y_{1,n}|^3 \big] \gamma_{Y,n}^3$,
we obtain $(h_{X,n}, h_{Y,n}) \convD \Nc(0, V)$.
We also obtain $(h_{X,n}, h_{Y,n}, W_n) \convD \Nc(0, V) \otimes \delta_0$ by Slutsky's Lemma.
We can therefore apply Skorokhods's almost sure representation theorem, see \cite[Theorem 2.19]{van2000asymptotic}.
It means that there exists
a probability space $(\tilde \Omega, \tilde \Uc, \tilde\Prob)$,
a sequence of random vectors $(\tilde h_{X,n}, \tilde h_{Y,n}, \tilde W_n)$ such that for every $n \geq 1$, $(\tilde h_{X,n}, \tilde h_{Y,n}, \tilde W_n) \equalD (h_{X,n}, h_{Y,n}, W_n)$,
and a random vector $(\tilde h_{X,\infty}, \tilde h_{Y,\infty}, \tilde W_\infty)$ following the distribution $\Nc(0, V) \otimes \delta_0$ such that
$(\tilde h_{X,n}, \tilde h_{Y,n}, \tilde W_n)$ $\convAS (\tilde h_{X,\infty}, \tilde h_{Y,\infty}, \tilde W_\infty)$, where the convergence is to be seen as of a sequence of random vectors defined on $(\tilde \Omega, \tilde \Uc, \tilde\Prob)$.
Let us define
\begin{align*}
    \tilde A_n :=
    \frac{\theta_{X,n} + \tilde h_{X,n}
    / (\sqrt{n} \gamma_{X,n})}
    {\theta_{Y,n} + \tilde h_{Y,n}
    / (\sqrt{n} \gamma_{Y,n})}
    - \frac{\theta_{X,n}}{\theta_{Y,n}}
    & \equalD
    \frac{\theta_{X,n} + h_{X,n}/(\sqrt{n} \gamma_{X,n})}
    {\theta_{Y,n} + h_{Y,n}/(\sqrt{n} \gamma_{Y,n})}
    - \frac{\theta_{X,n}}{\theta_{Y,n}} \\
    &= \frac{\meanX}{\meanY} - \frac{\theta_{X,n}}{\theta_{Y,n}} = A_n.
\end{align*}
Moreover, we have $\tilde W_n = \Indicator\{ \theta_{Y,n} + \tilde h_{Y,n}/(\sqrt{n} \gamma_{Y,n}) = 0 \}$ and $\tilde W_\infty = 0$ almost surely. We can define
$$\tilde \Omega^*
= \{ \tilde \omega \in \tilde \Omega :
\tilde W_n(\tilde \omega) \to 0 \text{ and }
\exists N>0, \forall n \geq N, \tilde h_{Y,n}(\tilde \omega) \neq 0\}.$$
By the almost sure convergence of $(\tilde h_{Y,n}, \tilde W_n)$, we get $\tilde \Prob(\tilde \Omega^*) = 1$, and for every $\tilde \omega \in \tilde \Omega^*$, $\tilde W_n(\tilde \omega) = 0$ and $\tilde h_{Y,n}(\tilde \omega) \neq 0$ for every $n$ large enough.
This means that for every given $\tilde \omega \in \tilde \Omega^*$, and for every $n$ large enough, $\tilde A_n$ is well-defined.
In the rest of the proof, we will fix such a $\tilde \omega \in \tilde \Omega^*$, so that all random variables may be considered as deterministic. By the almost sure representation theorem, this means that the equivalents and limits that will be obtained will still be valid in law in the original spaces $\Omega_n$.

{\bf First case: }
We have
\begin{align*}
    \tilde A_n &= \frac{\meanX}{\meanY} - \frac{\theta_{X,n}}{\theta_{Y,n}}
    = \frac{\theta_{X,n} + \tilde h_{X,n}
    /(\sqrt{n} \gamma_{X,n})}
    {\theta_{Y,n} + \tilde h_{Y,n}
    /(\sqrt{n} \gamma_{Y,n})}
    - \frac{\theta_{X,n}}{\theta_{Y,n}} \\
    &= \frac{\theta_{X,n} + \tilde h_{X,n}/ (\sqrt{n} \gamma_{X,n})}{\theta_{Y,n}}
    \bigg(1 - \frac{\tilde h_{Y,n}}{\sqrt{n} \gamma_{Y,n} \theta_{Y,n}}
    + O\big( (\sqrt{n} \gamma_{Y,n} \theta_{Y,n})^{-2} \big) \bigg)
    - \frac{\theta_{X,n}}{\theta_{Y,n}} \\
    &\sim \frac{- \theta_{X,n} \tilde h_{Y,n}}{\sqrt{n} \gamma_{Y,n} \theta_{Y,n}^2}
    + \frac{\tilde h_{X,n}}{\sqrt{n} \gamma_{X,n} \theta_{Y,n}},
\end{align*}
as claimed.

\medskip

{\bf Second case: } We have
\begin{align*}
    \tilde A_n &\sim \frac{\theta_{X,n} + \tilde h_{X,n}/(\sqrt{n} \gamma_{X,n})}
    {C/(\sqrt{n} \gamma_{Y,n}) + \tilde h_{Y,n}/(\sqrt{n} \gamma_{Y,n})} - \frac{\theta_{X,n}}{C/(\sqrt{n} \gamma_{Y,n})} \\
    &= \frac{\sqrt{n} \gamma_{Y,n} \theta_{X,n} + \tilde h_{X,n} \gamma_{Y,n} / \gamma_{X,n}}
    {C + \tilde h_{Y,n}} - \frac{\sqrt{n} \gamma_{Y,n}\theta_{X,n}}{C}.
\end{align*}
We factorize by $\theta_{X,n}$ in the latter expression, which completes the proof.

\medskip

{\bf Third case: } We have
\begin{align*}
    \tilde A_n &= \frac{\theta_{X,n} + \tilde h_{X,n}/(\sqrt{n} \gamma_{X,n})}
    {\theta_{Y,n} + \tilde h_{Y,n}/(\sqrt{n} \gamma_{Y,n})} - \frac{\theta_{X,n}}{\theta_{Y,n}}
    = \frac{\theta_{X,n} + \tilde h_{X,n}/(\sqrt{n} \gamma_{X,n})}
    {\big( \tilde h_{Y,n} + o(1) \big) / (\sqrt{n} \gamma_{Y,n})}
    - \frac{\theta_{X,n}}{\theta_{Y,n}} \\
    &\sim \frac{\sqrt{n}\theta_{X,n} \gamma_{Y,n}} {\tilde h_{Y,n}}
    + \frac{\tilde h_{X,n} \gamma_{Y,n}} {\tilde h_{Y,n} \gamma_{X,n}}
    - \frac{\theta_{X,n}}{\theta_{Y,n}} \\
    &\sim \theta_{X,n} \bigg( \frac{\sqrt{n} \gamma_{X,n}} {\tilde h_{Y,n}}
    - \frac{1}{\theta_{Y,n}} \bigg)
    + \frac{\tilde h_{X,n} \gamma_{Y,n}} {\tilde h_{Y,n} \gamma_{X,n}},
\end{align*}
and the result follows from the fact that $\sqrt{n} \gamma_{X,n} / \tilde h_{Y,n}$ is negligible compared to $1/\theta_{Y,n}$.
\begin{flushright}
    $\Box$
\end{flushright}

\subsection{Proof of Theorem~\ref{thm:validity_bootstrap_sequenceModels}}
\label{proof:thm:validity_bootstrap_sequenceModels}

For $b=1,2$, let $h_{X,n}
:= \sqrt{n} \gamma_{X,n} (\meanX - \theta_{X,n})$ (resp. $h^{Y,n}$),
$S_n:=(h_{X,n}, h_{Y,n})'$ and
$S_n^{(b)} := (h_{X,n}^{(b)}, h_{Y,n}^{(b)})'$, where
$h_{X,n}^{(b)} := \sqrt{n} \gamma_{X,n} (\overline{X}_n^{(b)} - \meanX)$
is the $b$-th bootstrap replication of $h_{X,n}$ (resp. $h_{Y,n}^{(b)}$).

\begin{lemma}
    We have $d_{BL}\left(P_{S_n^{(1)} \, | \, \left(X_{i,n},Y_{i,n}\right)_{i=1}^n}, \mathcal{N}(0,V) \right) \convAS 0.$
    \label{lemma:BL_bootstrap}
\end{lemma}

By the Central Limit Theorem, we have $S_n\convD S$ with $S\sim\mathcal{N}(0,V)$ and by Lemma~\ref{lemma:BL_bootstrap} (proved in Section~\ref{proof:lemma:BL_bootstrap}) and the triangle inequality, we get $d_{BL}\left(P_{S_n^{(1)}\mid \left(X_l,Y_l\right)_{l=1}^n},P_{S_n}\right)\convP 0$.
Combining both results, Lemma 2.2 in \cite{buch_koja_2019} gives us $$d_{BL}\left(P_{(S_n, \, S_n^{(1)}, \, S_n^{(2)})}, \, P_{S}^{\otimes 3}\right)\to 0.$$

\medskip

Let us define
$W_n :=
\Indicator\{ \theta_{Y,n} + h_{Y,n}/ (\sqrt{n} \gamma_{Y,n}) = 0 \}$
and remark that $W_n = 1$ whenever $\meanY = 0$. By assumption $\Prob(\meanY = 0) \to 0$, therefore we have $W_n \convD \delta_0$.
We define also $W_n^{(b)} := \Indicator\{ \meanY^{(b)} \}
= \Indicator\{ \meanY + h_{Y,n}^{(b)}/ (\sqrt{n}\gamma_{Y,n}) = 0\},$
so that $W_n^{(1)} = 1$ whenever $\meanY^{(b)} = 0$.
In the same way as previously, $W_n^{(b)} \convD \delta_0$ holds by assumption.
Let $Z_n = (S_n, W_n, S_n^{(1)}, W_n^{(1)}, S_n^{(1)}, W_n^{(2)})$ be a random vector of size $9$, and let $Z$ be a random vector of size $9$ following~$(P_{S} \otimes \delta_0)^{\otimes 3}$.

\medskip

By Slutsky's lemma, we have $d_{BL}\left(P_{Z_n} , P_Z \right)\to 0$ with our new notation.
Using Skorokhods's almost sure representation theorem \cite[Theorem 2.19]{van2000asymptotic}, there exists a probability space $\Omega^+$, a sequence of random vectors $Z_n^+ \in \Rb^9$ and a vector $Z^+$ defined on $\Omega^+$ such that $Z_n^+ \convAS Z^+$, $Z_n \equalD Z_n^+$ and $Z \equalD Z^+$.
Let us use the notation
\begin{align*}
    Z_n^+ &= \big( S_n^+ , W_n^+ , S_n^{(1)+} , W_n^{(1)+},
    S_n^{(1)+}, W_n^{(2)+} \big) \\
    &= \Big( h_{X,n}^+ , h_{Y,n}^+ , W_n^+ ,
    h_{X,n}^{(1)+} , h_{Y,n}^{(1)+} , W_n^{(1)+},
    h_{X,n}^{(2)+} , h_{Y,n}^{(2)+} , W_n^{(2)+} \Big) \\
    \text{ and } Z^+ &= \big( Z_1^+ , \, Z_2^+ , \, Z_3^+ \big),
\end{align*}
where $S_n^+ , S_n^{(1)+} , S_n^{(2)+}$ are random vectors of dimension $2$ and $Z_1^+ , Z_2^+ , Z_3^+$ are random vectors of dimension $3$.
We define
\begin{align*}
    A_n &:= \frac{\meanX}{\meanY} - \frac{\theta_{X,n}}{\theta_{Y,n}}
    = \frac{\theta_{X,n}+h_{X,n}/ (\sqrt{n}\gamma_{X,n}) }
    {\theta_{Y,n} + h_{Y,n}/ (\sqrt{n}\gamma_{Y,n}) }
    - \frac{\theta_{X,n}}{\theta_{Y,n}} \\
    A_n^{(b)} &:= \frac{\meanX + h_{X,n}^{(b)}/ (\sqrt{n}\gamma_{X,n}) }
    {\meanY + h_{Y,n}^{(b)}/ (\sqrt{n}\gamma_{Y,n}) }
    - \frac{\meanX}{\meanY} \\
    &= \frac{\theta_{X,n}+h_{X,n}/ (\sqrt{n}\gamma_{X,n}) +h_{X,n}^{(b)}/ (\sqrt{n}\gamma_{X,n}) }
    {\theta_{Y,n }+h_{Y,n}/ (\sqrt{n}\gamma_{Y,n}) +h_{Y,n}^{(b)}/ (\sqrt{n}\gamma_{Y,n}) }
    - \frac{\theta_{X,n}+h_{X,n}/ (\sqrt{n}\gamma_{X,n}) }{\theta_{Y,n}+h_{Y,n}/ (\sqrt{n}\gamma_{Y,n}) },
\end{align*}
and respectively their counterparts $A_n^+$ and $A_n^{(b)+}$ defined on $\Omega^+$. The following lemma, proved in Section~\ref{proof:lemma:def:tilde_omega_prob1}, ensures the existence of an event of probability $1$ on which every quantity is well-defined.
\begin{lemma}
    There exists an event $\tilde \Omega \subset \Omega^+$ such that $\Prob(\tilde\Omega)=1$ and such that for every $\omega \in \tilde \Omega$, and for all $n$ large enough, $h_{Y,n}^+(\omega) \neq 0$, $h_{Y,n}^{(1)+}(\omega) \neq 0$, $h_{Y,n}^{(2)+}(\omega) \neq 0$ and $A_n^+(\omega)$, $A_n^{(1)+}(\omega)$ and $A_n^{(2)+}(\omega)$ are well-defined.
    \label{lemma:def:tilde_omega_prob1}
\end{lemma}
In the next step, we fix $\omega\in\widetilde{\Omega}$ and let $C:=\lim_{n\to+\infty}\theta_{X,n}\gamma_{X,n}/\theta_{Y,n}\gamma_{Y,n}$ and
\begin{align*}
    &\sigma_n := \sqrt{n}\theta_{Y,n}\left( \gamma_{X,n}\mathds{1}_{\{C\in\Rb\}} + \gamma_{Y,n}\theta_{Y,n}/\theta_{X,n}\mathds{1}_{\{|C|=+\infty\}} \right).
\end{align*}
We restrict ourselves to the case $n^{1/2} \gamma_{Y,n} \theta_{Y,n} \to + \infty$. Theorem~\ref{thm:delta_method_sequence_models} therefore yields
\begin{align} \label{eq:proof_boot_1}
    &\sigma_n A_n^+(\omega) = \left\{
    \begin{array}{ll}
        -Ch_{Y,n}^+(\omega) + h_{X,n}^+(\omega) + o(1)
        \text{ if } C\in\Rb \\
        -h_{Y,n}^+(\omega) + o(1) \text{ else.}
    \end{array}
\right.
\end{align}

Furthermore, the same tools as those used in the proof of Theorem~\ref{thm:delta_method_sequence_models} plus the fact that $\theta_{Y,n}+h_{Y,n}^+(\omega)/ (\sqrt{n}\gamma_{Y,n}) \sim \theta_{Y,n}$ imply
\begin{align*}
    &\sigma_n A_n^{(b)+}(\omega) \\
    &\sim \sigma_n \left(
    \frac{-\left(\theta_{X,n}+h_{X,n}^+(\omega)/ (\sqrt{n}\gamma_{X,n}) \right)}
    {\sqrt{n}\gamma_{Y,n} \left(\theta_{Y,n}+h_{Y,n}^+(\omega)/ (\sqrt{n}\gamma_{Y,n}) \right)^2}
    h_{Y,n}^{(b)+}(\omega) \right. \\
    & \left. \qquad\qquad\qquad + \frac{1}{\sqrt{n}\gamma_{Y,n}
    \left(\theta_{Y,n}+h_{Y,n}^+(\omega)/ (\sqrt{n}\gamma_{X,n}) \right)}
    h_{X,n}^{(b)+}(\omega)\right) \\
    &\sim \sigma_n
    \left( \frac{-\left(\theta_{X,n} + h_{X,n}^+(\omega)/ (\sqrt{n}\gamma_{Y,n}) \right)}
    {\sqrt{n}\gamma_{Y,n}\theta_{Y,n}^2}
    h_{Y,n}^{(b)+}(\omega) + \frac{1}{\sqrt{n}\gamma_{X,n} \theta_{Y,n}}
    h_{X,n}^{(b)+}(\omega) \right).
\end{align*}

We can also remark that when $\theta_{X,n}+h_{X,n}^+(\omega)/ (\sqrt{n}\gamma_{X,n}) \sim \theta_{X,n}$
\begin{align} \label{eq:proof_boot_2}
    &\sigma_nA_n^{(b)+}(\omega) = \left\{
    \begin{array}{ll}
        -Ch_{Y,n}^{(b)+}(\omega) + h_{X,n}^{(b)+}(\omega) + o(1)
        \text{ if } C\in\Rb \\
        -h_{Y,n}^{(b)+}(\omega) + o(1) \text{ else.}
    \end{array}
\right.
\end{align}

When $\theta_{X,n}+h_{X,n}^+(\omega)/ (\sqrt{n}\gamma_{X,n})
= O \big( h_{X,n}^+(\omega)/ (\sqrt{n}\gamma_{X,n}) \big),$
we have $C=0$ and we find again that
\begin{align} \label{eq:proof_boot_3}
    &\sigma_n A_n^{(b)+}(\omega) = h_{X,n}^{(b)+}(\omega) + o(1).
\end{align}

\medskip

Let $D_n^+:= \left(-Ch_{Y,n}^+ + h_{X,n}^+\right) \mathds{1}_{\{|C| < + \infty\}} - h_{Y,n}^+ \mathds{1}_{\{|C|=+\infty\}}$ (resp. $D_n$, $D_n^{(b)}$ and $D_n^{(b)+}$), which corresponds to the dominant terms in Equations~\eqref{eq:proof_boot_1}, \eqref{eq:proof_boot_2} and~\eqref{eq:proof_boot_3} above.
By construction of $Z_n^+$ and $Z_n$, we have $Z_n^+ \convAS Z^+$, so that the continuous mapping theorem ensures that
$\left(D_n^+,D_n^{(1)+},D_n^{(2)+}\right) \convAS (U_1,U_2,U_3)$,
where for every $i\in\{1,2,3\}$, we define
$U_i^+ := \left(-C Z_{i,2}^+ + Z_{i,1}^+\right) \mathds{1}_{\{C\in\Rb\}}
- Z_{i,2}^+ \mathds{1}_{\{|C|=+\infty\}}$
where $Z_{i,1}^+$ (resp. $Z_{i,2}^+$) is the first (resp. second) component of the vector~$Z_i^+$.
Combining the triangle inequality, Equations~\eqref{eq:proof_boot_1}, \eqref{eq:proof_boot_2} and \eqref{eq:proof_boot_3}, we get
$$\left(\sigma_n A_n^+, \sigma_n A_n^{(1)+}, \sigma_n A_n^{(2)+}\right)
\convAS \big( U_1^+,U_2^+,U_3^+ \big).$$
Using the fact that for all $n \in \INT$
$\left(A_n, A_n^{(1)}, A_n^{(2)}\right)
\equalD \left( A_n^+, A_n^{(1)+}, A_n^{(2)+} \right),$
we obtain
$$\left(\sigma_n A_n, \sigma_n A_n^{(1)},\sigma_n A_n^{(2)}\right)
\convD (U_1^+,U_2^+,U_3^+).$$
Therefore, $d_{BL} \Big(
P_{\big(\sigma_n A_n, \, \sigma_n A_n^{(1)}, \, \sigma_n A_n^{(2)} \big)}
\, , \, P_{U_1^+}^{\otimes 3} \Big) \to 0$ as $n \to +\infty$ and $\sigma_n A_n \convD U_1^+$.
Applying Lemma 2.2 of~\cite{buch_koja_2019}, we can conclude that
$$d_{BL} \Big(
P_{ \sigma_n A_n^{(1)} \, | \, \left(X_{i,n},Y_{i,n}\right)_{i=1}^n}
\, , \, P_{U_1^+} \Big) \convP 0.$$
The conclusion follows from Lemma 23.3 in~\cite{van2000asymptotic}.
\begin{flushright}
    $\Box$
\end{flushright}

\subsubsection{Proof of Lemma~\ref{lemma:BL_bootstrap}}
\label{proof:lemma:BL_bootstrap}

Let $t=(t_X, t_Y )'\in \Rb^2$, and denote
$T_{i,n} = t_X \gamma_{X,n} X_{i,n} + t_Y \gamma_{Y,n} Y_{i,n}$ for $i=1, \dots,n$ and $T_{i,n}^{(1)}$ its bootstrap counterpart. Let also $V_{T_{1,n}}:=t'\Var[(\gamma_{X,n} X_{1,n} \, , \, \gamma_{Y,n} Y_{1,n})]t$ and $V_T:=t'Vt$.
We start by showing that for every $t\in\Rb^2$, $P_{ \sqrt{n}\left(\overline{T}_n^{(1)}-\overline{T}_n\right) \, | \, \left(X_{i,n},Y_{i,n}\right)_{i=1}^n}$ converges weakly to $P_T=\mathcal{N}\left(0,V_T\right)$ almost surely conditionally on $(X_{i,n},Y_{i,n})_{i=1}^n$ in the sense of the Lévy criterion for weak convergence, i.e.
\begin{align} \label{eq:levy_crit}
    &\left| \Expec\left[ e^{iu\sqrt{n}\left(\overline{T}_n^{(1)}-\overline{T}_n\right)} \mid \left(X_{i,n},Y_{i,n}\right)_{i=1}^n \right] - e^{u^2V_T/2} \right|
    \convAS 0 \quad \forall u\in\Rb.
\end{align}

\medskip

To do so, we have to check the steps of the proof of Theorem 23.4 in \cite{van2000asymptotic}. We have
\begin{align*}
    &\Expec\Bigg[ \overline{T}_n^{(1)} \mid
    \left( X_{i,n},Y_{i,n} \right)_{i=1}^n \Bigg] = \overline{T}_n \quad \text{ and } \quad \\
    &\Expec\Bigg[ \left(T_{i,n}^{(1)}-\overline{T}_n \right)^2 \mid \left( X_{i,n},Y_{i,n} \right)_{i=1}^n \Bigg] =  \frac{1}{n}\sum_{i=1}^nT_{i,n}^2-\overline{T}_n^2 .
\end{align*}

The first requirement is to ensure almost sure convergence to~$0$ of both quantities $\left| \overline{T}_n - \Expec\left[ T_{1,n} \right] \right|$ and $\left| \frac{1}{n}\sum_{i=1}^nT_{i,n}^2-\overline{T}_n^2 - V_T \right|$. Under the assumption that $\sup_{n\in\INTST}\Expec\left[ \left| T_{1,n} \right|^{4+\delta} \right] < +\infty$, observe that all the conditions of Theorem 2.2 in \cite{gut1992complete} are satisfied with $p=1$.
We can thus conclude that $\left| \overline{T}_n - \Expec\left[ T_{1,n} \right] \right| \convAS 0$ and $\left| \frac{1}{n}\sum_{i=1}^n T_{i,n}^2 - \Expec\left[ T_{1,n}^2 \right] \right| \convAS 0$. Now using the fact that
\begin{align*}
    \left| \frac{1}{n}\sum_{i=1}^nT_{i,n}^2-\overline{T}_n^2 - V_T \right| \leq
    &\left| \frac{1}{n}\sum_{i=1}^n T_{i,n}^2 - \Expec\left[ T_{1,n}^2 \right] \right| + \left| \overline{T}_n - \Expec\left[ T_{1,n} \right] \right|^2 \\
    &+ 2\left| \Expec\left[ T_{1,n} \right]\left(\overline{T}_n-\Expec\left[T_{1,n}\right]\right)\right| + \left| \Expec\left[ T_{1,n}^2 \right] - \Expec\left[ T_{1,n} \right]^2 - V_T \right|,
\end{align*}
as well as $\left| \Expec\left[ T_{1,n} \right] \right| = O(1)$ and $\left|\left| \Var[(\gamma_{X,n} X_{1,n} \, , \, \gamma_{Y,n} Y_{1,n})] - V \right|\right| = o(1)$, to conclude that $\left| \frac{1}{n}\sum_{i=1}^nT_{i,n}^2-\overline{T}_n^2 - V_T \right| \convAS 0$.

\medskip

The second requirement is to check the Lindeberg condition for the bootstrap which writes
\begin{align*}
    \Expec\left[ \left|T_{1,n}^{(1)}\right|^2 \Indicator\left\{ \left|T_{1,n}^{(1)}\right|^2 > \epsilon\sqrt{n}\right\} \mid \left(X_{i,n},Y_{i,n}\right)_{i=1}^n \right] &= \frac{1}{n}\sum_{i=1}^n \left|T_{i,n}\right|^2 \Indicator\left\{ \left|T_{i,n}\right|^2 > \epsilon\sqrt{n}\right\} \\
    &\convAS 0 \quad \forall\epsilon>0.
\end{align*}

Let $M:\epsilon \mapsto M(\epsilon)$ be some function of $\epsilon$ to be defined later that does not depend on $n$ and satisfies $0<M(\epsilon)<+\infty$ $\forall \epsilon>0$. For such a function, there exists for every $\epsilon>0$, a $n_{\epsilon}$ such that for every $n>n_{\epsilon}$,
\begin{align*}
    &\frac{1}{n}\sum_{i=1}^n \left|T_{i,n}\right|^2 \Indicator\left\{ \left|T_{i,n}\right|^2 > \epsilon\sqrt{n}\right\} \leq \frac{1}{n}\sum_{i=1}^n \left|T_{i,n}\right|^2 \Indicator\left\{ \left|T_{i,n}\right|^2 > M(\epsilon) \right\} \, a.s.
\end{align*}

By the triangle inequality,
\begin{align*}
    &\frac{1}{n}\sum_{i=1}^n \left|T_{i,n}\right|^2 \Indicator\left\{ \left|T_{i,n}\right|^2 > M(\epsilon) \right\} \leq \\
    &\qquad\qquad \left| \frac{1}{n}\sum_{i=1}^n \left|T_{i,n}\right|^2 \Indicator\left\{ \left|T_{i,n}\right|^2 > M(\epsilon) \right\} - \Expec\left[ \left|T_{1,n}\right|^2 \Indicator\left\{ \left|T_{1,n}\right|^2 > M(\epsilon) \right\} \right] \right| \\
    &\qquad\qquad\qquad\qquad + \Expec\left[ \left|T_{1,n}\right|^2 \Indicator\left\{ \left|T_{1,n}\right|^2 > M(\epsilon) \right\} \right].
\end{align*}
The first term in the upper bound converges to 0 almost surely for every $\epsilon>0$ under the assumption $\sup_{n\in\INTST}\Expec\left[ \left| T_{1,n} \right|^{4+\delta} \right] < +\infty$ thanks to Theorem 2.2 in \cite{gut1992complete}. The second term in the upper bound can be bounded with the Cauchy-Schwarz and Markov inequalities
\begin{align*}
    &\Expec\left[ \left|T_{1,n}\right|^2 \Indicator\left\{ \left|T_{1,n}\right|^2 > M(\epsilon) \right\} \right] \leq \frac{ \sup_{n\in\INTST}\sqrt{\Expec\left[ \left| T_{1,n} \right|^4 \right] \Expec\left[ \left| T_{1,n} \right| \right]. }}{\sqrt{M(\epsilon)}}
\end{align*}

Picking $M(\epsilon)=\epsilon^{-1}\sup_{n\in\INTST}\Expec\left[ \left| T_{1,n} \right|^4 \right] \Expec\left[ \left| T_{1,n} \right| \right]$, we get that for every $\epsilon>0$
\begin{align*}
    &\limsup_{n\to +\infty}\frac{1}{n}\sum_{i=1}^n \left|T_{i,n}\right|^2 \Indicator\left\{ \left|T_{i,n}\right|^2 > \epsilon\sqrt{n}\right\} \leq \epsilon \quad a.s.
\end{align*}

Letting $\epsilon$ go to 0, we see that the Lindeberg condition is satisfied. This entails that \eqref{eq:levy_crit} is satisfied.

Arguments underpinning the Cramer-Wold device are valid as well so that we can claim that for every $t\in\Rb^2$
\begin{align} \label{eq:as_cramer_wold}
    &\left|\Expec\left[ e^{it'S_n^{(1)}} \mid \left(X_{i,n},Y_{i,n}\right)_{i=1}^n \right]-e^{-t'Vt/2}\right|\convAS 0,
\end{align}
where $S_n^{(1)} := \sqrt{n}\Bigg( \gamma_{X,n}\left(\overline{X}_n^{(1)} - \meanX \right) , \gamma_{Y,n}\left(\overline{Y}_n^{(1)} - \meanY \right) \Bigg).$

Let $\Omega$ be the set of probability one on which \eqref{eq:as_cramer_wold} occurs. For every $\omega\in\Omega$, $$\Bigg(P_{S_n^{(1)} \mid \left(X_{i,n},Y_{i,n}\right)_{i=1}^n = \left(X_{i,n}(\omega),Y_{i,n}(\omega)\right)_{i=1}^n }\Bigg)_{n\in\INTST}$$
is a sequence of nonrandom probability measures for which all weak convergence criteria are equivalent. In particular, for every $\omega\in\Omega$, the validity of the Lévy criterion due to \eqref{eq:as_cramer_wold} ensures that
\begin{align*}
    &d_{BL}\left( P_{S_n^{(1)} \mid \left(X_{i,n},Y_{i,n}\right)_{i=1}^n = \left(X_{i,n}(\omega),Y_{i,n}(\omega)\right)_{i=1}^n } , \mathcal{N}(0,V) \right) = o(1).
\end{align*}

This is enough to conclude.
\begin{flushright}
    $\Box$
\end{flushright}

\subsubsection{Proof of Lemma~\ref{lemma:def:tilde_omega_prob1}}
\label{proof:lemma:def:tilde_omega_prob1}
The vector $\big(W_n^+, W_n^{(1)+},W_n^{(2)+} \big)$ converges almost surely to $(0,0,0)$.
As a consequence, there exists an event $\tilde \Omega^1$ of probability $1$ such that $\forall \omega \in \tilde \Omega^1, \, \big(W_n^+(\omega),$ $W_n^{(1)+}(\omega),W_n^{(2)+}(\omega) \big) = (0,0,0)$ for $n$ large enough.
As $\big( h_{Y,n}^+ , h_{Y,n}^{(1)+} , h_{Y,n}^{(2)+} \big)$ converges almost surely to a continuous vector, there exists an event $\tilde \Omega^2$ of probability $1$ such that $\forall \omega \in \tilde \Omega^2$, the components of $\big( h_{Y,n}^+(\omega) , h_{Y,n}^{(1)+}(\omega) , h_{Y,n}^{(2)+}(\omega) \big)$ are all non-zero for $n$ large enough.
We finally define $\tilde \Omega := \tilde \Omega^1 \cap \tilde \Omega^2$, which is of probability $1$ and satisfies the stated conditions. $\Box$


\subsection{Proof of Example~\ref{example:sequence_Bernoulli_bootstrap}}
\label{proof:example:sequence_Bernoulli_bootstrap}

We have
\begin{align*}
    \Prob(\meanY^{(1)} = 0)
    &= \Expec \Big[\Prob \big(\meanY^{(1)} = 0 \, \big| \, (Y_{i,n})_{i=1}^n \big) \Big]\\
    &= \Expec \Big[\Prob \big(Y_{1,n}^{(1)} = 0, \dots, Y_{n,n}^{(1)} = 0 \, \big| \, (Y_{i,n})_{i=1}^n \big) \Big] \\
    &= \Expec \Big[\Prob \big(Y_{1,n}^{(1)} = 0 \, \big| \, (Y_{i,n})_{i=1}^n \big)^n \Big]
    = \Expec \Big[ \big( S_n / n) \big)^n \Big],
\end{align*}
where $S_n := \sum_{i=1}^n (1-Y_{i,n}) \sim Bin(n, 1-p_n)$.
Therefore, for any $x > 0$,
\begin{align*}
    \Prob(\meanY^{(1)} = 0)
    &= \sum_{k=1}^n (k/n)^n \Prob \big[S_n = k \big] \\
    &\leq \sum_{k=1}^{\lfloor n(1-p_n) + x \rfloor} (k/n)^n
    \Prob \big[S_n = k \big]
    + \Prob \big[ S_n \geq n(1-p_n) + x \big] \\
    &\leq \left( \frac{n (1-p_n) + x}{n} \right)^n
    + \Prob \big[S_n \geq n(1-p_n) + x \big] \\
    &\leq \big( 1 - p_n + x/n \big)^n
    + \Prob \big[S_n - n(1-p_n) \geq x \big].
\end{align*}
Let $\tilde S_n := \big(S_n - n(1-p_n)\big) / \sqrt{n p_n (1-p_n)} = O_P(1)$ be the renormalized version of $S_n$ and choose $x = n^a \sqrt{n p_n (1-p_n)}$ for $a = (1-b)/3 > 0$. Then
\begingroup \allowdisplaybreaks
\begin{align*}
    \Prob(\meanY^{(1)} = 0)
    &\leq \big( 1 - p_n + n^a \sqrt{p_n (1-p_n) / n}\big)^n
    + \Prob \big[\tilde S_n \geq n^a \big] \\
    &\leq \exp \Big(n \ln \big( 1 - p_n + n^a \sqrt{p_n (1-p_n) / n} + o(p_n) \big) \Big) + o(1) \\
    &\leq \exp \Big(n \big( n^{a - b/2 - 1/2} - n^{-b} + o(n^{-b})
    \big) \Big) + o(1) \\
    &\leq \exp \Big(n^{1/3 - b/3 - b/2 + 1/2} - n^{-b+1} + o(n^{-b+1})
    \Big) + o(1) \\
    &\leq \exp \Big(n^{(1-b)5/6} - n^{1-b} + o(n^{-b+1}) \Big) + o(1) \\
    &\leq \exp \big(- n^{1-b} \big) + o(1) = o(1),
\end{align*}
\endgroup
which completes the proof. $\Box$

\medskip

\subsection{Proof of Theorem~\ref{thm:bt_finite_ci_easy}}
\label{proof:thm:bt_finite_ci_easy}

We fix arbitrary $n \in \Nstar$ and $\eps \in \Rstarplus$.
Combining the triangle inequality, the bound $|\meanX| \leq |\meanX - \expecXn| + |\expecXn|$ and Assumptions~\ref{hyp:basic_dgp_XY} to~\ref{hyp:bt_easy_dgp_well_separated}, we get
\begin{align*}
    \bigg| \frac{\meanX}{\meanY} - \frac{\expecXn}{\expecYn} \bigg|
    &\leq |\meanX| \cdot
    \bigg| \frac{1}{\meanY} - \frac{1}{\expecYn} \bigg|
    + \frac{1}{\expecYn} \bigg| \meanX - \expecXn \bigg| \displaybreak[0] \\
    &\leq \frac{\left(\big| \meanX - \expecXn \big|+\sqrt{\upperXn}\right) \big| \meanY - \expecYn\big|}{\aYn\lowerYn} + \frac{\big|\meanX-\expecXn\big|}{\lowerYn}.
\end{align*}

Consequently, the event considered in Theorem \ref{thm:bt_finite_ci_easy} is included in the event
\begin{align}
\label{proof:eq:bt_easy_case_event}
    \frac{\left(\big| \meanX - \expecXn \big|+\sqrt{\upperXn}\right) \big| \meanY - \expecYn\big|}{\aYn\lowerYn} &+ \frac{\big|\meanX-\expecXn\big|}{\lowerYn} \nonumber \\
    & \qquad >\frac{\big(\varepsilon+\sqrt{\upperXn}\big)\varepsilon}{\aYn\lowerYn}
    + \frac{\varepsilon}{\lowerYn}.
\end{align}
If both ${|\meanX-\expecXn|}$ and ${|\meanY-\expecYn|}$ are inferior or equal to $\eps$, event \eqref{proof:eq:bt_easy_case_event} cannot happen.
By contraposition, we obtain:
\begin{align*}
    &\Prob\Bigg(\frac{\left(\big| \meanX - \expecXn \big|+\sqrt{\upperXn}\right) \big| \meanY - \expecYn\big|}{\aYn\lowerYn} + \frac{\big|\meanX-\expecXn\big|}{\lowerYn} \\
    & \qquad\qquad\qquad\qquad\qquad\qquad\qquad\qquad\qquad >
    \frac{\big(\varepsilon+\sqrt{\upperXn}\big)\varepsilon}{\aYn\lowerYn}
    + \frac{\varepsilon}{\lowerYn}\Bigg)  \\
    &\quad \leq \Prob\left( \left\{\big|\meanX-\expecXn\big| > \eps \right\}  \cup \left\{\big|\meanY-\expecYn\big| > \eps \right\} \right) \\
    &\quad \leq \Prob\left( \big|\meanX-\expecXn\big|> \eps \right) +
    \Prob \left(\big|\meanY-\expecYn\big|> \eps\right),
\end{align*}
where we use the union bound for the last inequality.
The first conclusion follows from using twice Bienaymé-Chebyshev's inequality applied to the variables $\meanX$ and $\meanY$ and the fact that under Assumptions~\ref{hyp:basic_dgp_XY} and \ref{hyp:bt_dgp_XY} and Jensen's inequality, $\Var\left[X_{1,n}\right]\leq\upperXn$ and $\Var\left[Y_{1,n}\right]\leq\upperYn-\lowerYn^2$. The second conclusion follows from solving $(\upperXn+\upperYn-\lowerYn^2)/(n\varepsilon^2)=\alpha$.
\begin{flushright}
    $\Box$
\end{flushright}

\subsection{Proof of Theorem~\ref{thm:concentration_inequality_BC_general_case}}
\label{proof:thm:concentration_inequality_BC_general_case}

We start by introducing and proving an intermediate lemma that is also used to prove Theorem~\ref{thm:hoef_finite_ci}. For a random variable~$U$, ${\eps > 0}$, and ${\epstilde \in (0,1)}$ we define the following events:
\begin{align*}
    A_\eps^U
    := \Big\{ \big|\meanU - \Expec[U] \big|
    \leq \eps \Big\}, \text{ and }
    \tilde A_{\epstilde}^U
    := \Big\{ \big|\meanU - \Expec[U] \big|
    \leq \epstilde \big| \Expec[U] \big| \Big\}.
\end{align*}

\begin{lemma}
\label{lemma:bound_general_ratioMean}
    Assume that Assumption~\ref{hyp:basic_dgp_XY} holds.
    Then for every ${n\in\Nstar}$, ${\eps>0}$ and ${\epstilde\in (0,1)}$, we have
    \begin{align*}
        & \Prob \Bigg( \bigg| \frac{\meanX}{\meanY} - \frac{\expecXn}{\expecYn} \bigg|
        >
        \bigg( \frac{\left(|\expecXn|+\eps\right) \epstilde}{(1 - \epstilde)^2} + \eps \bigg) \frac{1}{|\expecYn|} \Bigg) \\
        & \qquad\qquad \leq 1 - \Prob \left(\eventAXeps\right)
        + 1 - \Prob \left(\eventAYepstilde\right).
    \end{align*}
\end{lemma}

We fix arbitrary ${n\in\Nstar}$, ${\eps>0}$ and ${\epstilde\in(0,1)}$.
By Lemma~\ref{lemma:bound_general_ratioMean}, we have
\begin{align*}
    &\Prob \Bigg( \bigg| \frac{\meanX}{\meanY} - \frac{\expecXn}{\expecYn} \bigg|
    > \bigg( \frac{\left(|\expecXn|+\varepsilon\right) \tilde \varepsilon}{(1 - \tilde \varepsilon)^2} + \varepsilon \bigg) \frac{1}{|\expecYn|} \Bigg) \\
    &\leq 1 - \Prob \Big( \big| \meanX - \expecXn \big|
    \leq \varepsilon \Big)
    + 1 - \Prob \Big( \big| \meanY - \expecYn \big|
    \leq \tilde{\varepsilon} \big| \expecYn \big| \Big).
\end{align*}
Using Jensen's inequality and Assumption~\ref{hyp:bt_dgp_XY}, we have ${|\expecXn| \leq (\upperXn)^{1/2}}$, and Assumption~\ref{hyp:basic_dgp_XY} entails ${1/|\expecYn| \leq 1/\lowerYn}$.
Consequently, we get
\begin{align*}
    &\Prob \Bigg( \bigg| \frac{\meanX}{\meanY} - \frac{\expecXn}{\expecYn} \bigg|
    > \bigg( \frac{\left(\sqrt{\upperXn}+\varepsilon\right) \tilde \varepsilon}{(1 - \tilde \varepsilon)^2} + \varepsilon \bigg) \frac{1}{\lowerYn} \Bigg) \\
    &\leq 1 - \Prob \Big( \big| \meanX - \expecXn \big|
    \leq \varepsilon \Big)
    + 1 - \Prob \Big( \big| \meanY - \expecYn \big|
    \leq \tilde{\varepsilon} \big| \expecYn \big| \Big).
\end{align*}
Using Bienaymé-Chebyshev's inequality twice gives the bounds
\begin{align*}
  1 - \Prob \Big( \big| \meanX - \expecXn \big| \leq \eps \Big) &\leq \frac{\Var\left[X_{1,n}\right]}{n \eps^2}\\
  1 - \Prob \Big( \big| \meanY - \expecYn \big| \leq \epstilde \big| \expecYn \big| \Big) &\leq \frac{\Var\left[Y_{1,n}\right]}{n {\epstilde}^2 \left(\expecYn\right)^2}.
\end{align*}
For the numerator,
$\Var\left[X_{1,n}\right] = \Expec\left[X_{1,n}^2\right] - \left(\expecXn\right)^2 \leq \Expec\left[X_{1,n}^2\right] \leq \upperXn$ using Assumption~\ref{hyp:bt_dgp_XY}.
For the denominator, Assumption~\ref{hyp:basic_dgp_XY} immediately entails that ${1/\lowerYn^2}$ is an upper bound on ${1/\left(\expecYn\right)^2}$ and ${\lowerYn^2}$ a lower bound on $\left(\expecYn\right)^2$.
Therefore
\begin{equation*}
    \frac{\Var\left[Y_{1,n}\right]}{n {\epstilde}^2 \left(\expecYn\right)^2}
    \leq \frac{\Expec\left[Y_{1,n}^2\right] - {\lowerYn}^2}{n {\epstilde}^2 {\lowerYn}^2}
    \leq \frac{\upperYn - {\lowerYn}^2}{n {\epstilde}^2 {\lowerYn}^2},
\end{equation*}
where the second inequality uses Assumption~\ref{hyp:bt_dgp_XY}.

Combining the two bounds yields the following upper bound on the probability considered in Theorem~\ref{thm:concentration_inequality_BC_general_case}
\begin{equation}
\label{eq:upper_bound_proba_concentration_inequality_bt_case}
    \frac{\upperXn}{n \eps^2} +
    \frac{\upperYn - {\lowerYn}^2}{n {\epstilde}^2 {\lowerYn}^2},
\end{equation}
as claimed.

\medskip

For the second part of Theorem~\ref{thm:concentration_inequality_BC_general_case}, for a fixed $\alpha$, we equalize each of the two terms in \eqref{eq:upper_bound_proba_concentration_inequality_bt_case} to $\alpha/2$ and solve for $\eps$ and $\epstilde$, which yields:
\begin{equation*}
    \eps^2 = \frac{2\upperXn}{n\alpha}\, \text{ and }
    {\epstilde}^2 = \frac{2\left(\upperYn - {\lowerYn}^2\right)}{n\alpha\lowerYn^2}.
\end{equation*}
The bound $\overline{\alpha}_n$ comes from the fact that $\epstilde$ needs to be smaller than $1$.
\begin{flushright}
    $\Box$
\end{flushright}

\subsubsection{Proof of Lemma~\ref{lemma:bound_general_ratioMean}}

We fix arbitrary ${\eps > 0}$ and ${\epstilde\in(0,1)}$.
Without loss of generality, we can assume that $\expecYn > 0$ and $\expecXn \geq 0$.

First, using the union bound, note that the event $\eventAXeps \cap\, \eventAYepstilde$ holds with a probability bigger than $\Prob\left(\eventAXeps\right) + \Prob\left(\eventAYepstilde\right) - 1$.
Hence, its complement is of probability lower than $1 - \Prob\left(\eventAXeps\right)
+ 1 - \Prob\left(\eventAYepstilde\right)$.

Second, we show that the event considered in Lemma~\ref{lemma:bound_general_ratioMean}  is included in the complement of $\eventAXeps \cap\, \eventAYepstilde$, which concludes the proof.
To do so, we reason by contraposition and do the following computations on the event $\eventAXeps \cap\, \eventAYepstilde$.

\noindent By the triangle inequality, we get
\begin{align*}
    \bigg| \frac{\meanX}{\meanY} - \frac{\expecXn}{\expecYn} \bigg|
    \leq |\meanX| \cdot \bigg| \frac{1}{\meanY} - \frac{1}{\expecYn} \bigg|
    + \frac{1}{\expecYn} \bigg| \meanX - \expecXn \bigg|.
\end{align*}
We now bound the first term using the mean value theorem applied to the function $f(x) := 1 / (x + \expecYn)$
\begin{align*}
    \bigg| \frac{1}{\meanY} - \frac{1}{\expecYn} \bigg|
    = \Big| f(\meanY - \expecYn) - f(0) \Big|
    &\leq \frac{|\meanY - \expecYn|}{(1-\tilde \varepsilon)^2 \expecYn^2} \\
    &\leq \frac{\tilde \varepsilon \expecYn}{(1-\tilde \varepsilon)^2 \expecYn^2},
\end{align*}
where the first inequality uses the following observation: on the event $\eventAYepstilde$, a lower bound on ${|x + \expecYn|}$ with $x$ varying between $0$ and $\meanY-\expecYn$ is $(1-\epstilde)|\expecYn|$.
Therefore, on $\eventAXeps \cap\, \eventAYepstilde$,
\begin{align*}
    \bigg| \frac{\meanX}{\meanY} - \frac{\expecXn}{\expecYn} \bigg|
    &\leq |\meanX| \cdot \frac{\epstilde \expecYn}{(1-\epstilde)^2 \expecYn^2}
    + \frac{\varepsilon}{\expecYn} \\
    &\leq \big(|\expecXn| + |\meanX - \expecXn|\big) \frac{\epstilde}{(1 - \epstilde)^2 \expecYn} +  \frac{\eps}{\expecYn} \\
    &\leq \frac{ \big(|\expecXn| + \eps\big)\epstilde}{(1 - \epstilde)^2 \expecYn}  +  \frac{\eps}{\expecYn}\,,
\end{align*}
where we use the triangle inequality to get the second line.
It is indeed the complement of the event considered in the statement of Lemma~\ref{lemma:bound_general_ratioMean}.
\begin{flushright}
    $\Box$
\end{flushright}

\subsection{Proof of Theorem \ref{thm:bt_necessary_cond_alpha}}
\label{proof:thm:bt_necessary_cond_alpha}

This theorem relies crucially on the following lemma.

\begin{lemma}
    \label{lemma:lowerBound_probInfty}
    For each $\xi$ in the interval $\Big(0,1\wedge\big(\upperYn/\lowerYn^2-1\big)\Big)$, there exists a distribution $P_{n, \, \xi} \in \Pfam$ such that
    $\Prob\left( \meanY = 0\right)
    \geq \tilde \alpha_n(\xi),$
    where $\tilde \alpha_n(\xi)
    := \big(1 - (1+\xi) \lowerYn^2 / \upperYn \big)^n$.
\end{lemma}
\noindent
Note that the interval $\Big(0,1\wedge\big(\upperYn/\lowerYn^2-1\big)\Big)$ is not empty since we have assumed $\upperYn/\lowerYn^2>1$.

By Lemma~\ref{lemma:lowerBound_probInfty}, for every
$\xi < 1\wedge\big(\upperYn/\lowerYn^2-1\big)$, there exists a distribution $P_{n, \, \xi}$ such that
$\Prob\left( \meanY = 0 \right)
\geq \tilde \alpha_n(\xi).$
Taking the supremum over $\xi$, we deduce that
\begin{align*}
    \sup_{P_n \in \Pfam}
    \Prob \left( \meanY = 0 \right)
    \geq \sup_\xi \tilde \alpha_n(\xi) = \underline{\alpha}_n.
\end{align*}
Using the assumption that $I_n$ is undefined whenever $\meanY=0$, we deduce that
$\Prob \big( I_n \, undefined \big)
\geq \underline{\alpha}_n.$
\begin{flushright}
    $\Box$
\end{flushright}


\subsubsection{Proof of Lemma~\ref{lemma:lowerBound_probInfty}}
We consider the following distribution on $\Rb$
\begin{align*}
    &P_{n,\lowerYn,\upperYn,c,\xi}
    :=\left(\frac{c}{n}\right)^{1/n} \delta_{\lsin 0\rsin}
    + \frac{1}{2}\left(1-\left(\frac{c}{n}\right)^{1/n}\right)
    \delta_{\{y_{c-}\}}
    + \frac{1}{2}\left(1-\left(\frac{c}{n}\right)^{1/n}\right)
    \delta_{\{y_{c+}\}},
\end{align*}
where $c\in(0,n)$ is some constant to be chosen later,
$y_{c-} := \lowerYn(1-\sqrt{\xi}) / (1-(c/n)^{1/n})$
and $y_{c+} := \lowerYn(1+\sqrt{\xi}) / (1-(c/n)^{1/n})$.
Let $Y_{1,n}\sim P_{n,\lowerYn,\upperYn,c,\xi_n}$.
Observe that $\expecYn=\lowerYn$ and
$\expecYnsq = \lowerYn^2 (1+\xi_n) / \big(1-(c/n)^{1/n} \big)$.
With the choice
\begin{align*}
    &c = c_n :=
    n \left(1-\frac{\lowerYn^2}{\upperYn}\left(1+\xi\right)\right)^n,
\end{align*}
we have $\expecYnsq=\upperYn$.
Note that $\Cnalpha$ is strictly positive,
because $1-\frac{\lowerYn^2}{\upperYn}\left(1+\xi_n\right)$ is positive.
This is equivalent to $\upperYn/\lowerYn^2>1+\xi_n$, which is true by assumption.

\medskip

Consider now the following product measure on $\Rb^2$
defined by $P_n:=\delta_{\lsin\sqrt{\upperXn}\rsin}
\otimes P_{n,\lowerYn,\upperYn,c_n,\xi}.$
Let $(X_{i,n},Y_{i,n})_{i=1}^n\simiid P_n$. These random vectors satisfy $\expecXnsq=\upperXn$, $\expecYn=\lowerYn$ and $\expecYnsq=\upperYn$.
The next step is to build a lower bound on the event $\{ \meanY = 0 \}$.


The assumption that $(X_{i,n},Y_{i,n})_{i=1}^n\simiid P_n$ and the construction of $P_{n,\lowerYn,\upperYn,c_n,\xi}$ imply that
\begin{align*}
    \Prob\left( \meanY=0 \right)
    = \frac{c_n}{n}
    = \left(1-\frac{\lowerYn^2}{\upperYn}\left(1+\xi\right)\right)^n
    = \tilde \alpha_n(\xi).
\end{align*}


\begin{flushright}
    $\Box$
\end{flushright}

\subsection{Proof of Theorem \ref{thm:bt_l_bound_ci_length}}
\label{proof:thm:bt_l_bound_ci_length}

To prove Theorem~\ref{thm:bt_l_bound_ci_length}, we need the following lemma.

\begin{lemma}\label{lem:tech_lemma_lower_bounds}
    For every integer $n\geq 7$ and every $x\in\left(0,1\right)$, $x\left(1-x/n\right)^{n-1}\geq x/3$.
\end{lemma}

We start using arguments developed in the proof of \cite{catoni2012challenging}[Proposition 6.2]. We detail those for the sake of clarity. For every $n\in\INTST$ and $\eta>\sqrt{\upperXn}/n$, let us define the following distribution on $\Rb$, which will be used for the
variable in the numerator\footnote{The notation $\delta$ denotes the Dirac distribution.}:
\begin{align*}
    &P_{n,\upperXn,\eta}:=\frac{\upperXn}{2n^2\eta^2}\delta_{\lsin -n\eta\rsin}+\left(1-\frac{\upperXn}{n^2\eta^2}\right)\delta_{\lsin 0\rsin}+\frac{\upperXn}{2n^2\eta^2}\delta_{\lsin n\eta\rsin}.
\end{align*}

This distribution is symmetric, centered and has variance $\upperXn$. As shown in \cite{catoni2012challenging}, every i.i.d. sample $(X_{i,n})_{i=1}^n$ drawn from $P_{n,\upperXn,\eta}$ satisfies
\begin{align*}
    &\Prob\left(\overline{X}_n\leq-\eta\right)
    = \Prob\left(\overline{X}_n\geq\eta\right)
    \geq \Prob\left(\overline{X}_n=\eta\right) \\
    &\geq\sum_{i=1}^n\Prob\left(X_{i,n}=n\eta,X_{j,n}=0, \; \forall j\neq i\right)
    =\frac{\upperXn}{2 n\eta^2}\left(1-\frac{\upperXn}{\eta^2n^2}\right)^{n-1}.
\end{align*}

Note further that for every integer $n\geq 2$,  $\Prob\left(\overline{X}_n\geq\eta\right)\geq\Prob\left(\overline{X}_n=\eta\right)$ becomes a strict inequality strict and for every $\xi\in(0,1)$ $\left\{\left|\meanX\right|\geq\eta\right\}\subseteq\left\{\left|\meanX\right|>\xi\eta\right\}$. As a result, if $(X_{i,n})_{i=1}^n \simiid P_{n,\upperXn,\eta}$, for every $\eta>0$, we have
\begin{align}
    \label{eq:catoni_lower_bound}
    &\Prob\left(\left|\overline{X}_n\right|>\xi\eta\right)
    > \frac{\upperXn}{n\eta^2}\left(1-\frac{\upperXn}{\eta^2n^2}\right)^{n-1}.
\end{align}

The following steps do not show up in \cite{catoni2012challenging} since they are specific to controlling ratios of expectations and sample averages. For every $n\in\INTST$, let us define the following distribution on $\Rb$, which will be used for the variable in the denominator
\begin{align*}
    &P_{n,\lowerYn,\upperYn}:=\frac{1}{2}\delta_{\lsin\lowerYn-\sqrt{\upperYn-\lowerYn^2}\rsin}+\frac{1}{2}\delta_{\lsin\lowerYn+\sqrt{\upperYn-\lowerYn^2}\rsin}.
\end{align*}

Let $(X_{i,n},Y_{i,n})_{i=1}^n\simiid P_n:=P_{n,\upperXn,\eta}\otimes P_{n,\lowerYn,\upperYn}$. Observe that $\expecYn=\lowerYn$ and $\expecYnsq=\upperYn$. Furthermore, $\left|\meanY\right|\leq \lowerYn+\sqrt{\upperYn-\lowerYn^2}$ almost surely. This implies that for every $\eta>0$ and $\xi\in(0,1)$, the following holds
\begin{align*}
    &\left\{\left|\meanX\right|> \left(\lowerYn+\sqrt{\upperYn-\lowerYn^2}\right)\xi\eta\right\}\subseteq\left\{\left|\frac{\meanX}{\meanY}\right|>\xi\eta\right\}.
\end{align*}



For fixed $n \geq 7$ and $\alpha \in
\left(0,1\wedge n/\left(\lowerYn+\sqrt{\upperYn-\lowerYn^2}\right)^2\right)$,
we choose $\eta = \eta(\alpha) = \sqrt{v_n / 3 n \alpha}$.
Combining the above inclusion with \eqref{eq:catoni_lower_bound}, and Lemma~\ref{lem:tech_lemma_lower_bounds} (with the choice $x = 3 \alpha$),
we conclude that there exists a distribution on $\Rb^2$, namely $P_n$,
that fulfills Assumptions~\ref{hyp:basic_dgp_XY} and~\ref{hyp:bt_dgp_XY}
such that
\begin{align*}
    &\Prob\left(\left|\frac{\overline{X}_n}{\overline{Y}_n}-\frac{\expecXn}{\expecYn}\right|>\xi\sqrt{\frac{v_n}{3n\alpha}}\right)>\alpha,
\end{align*}
which completes the proof.
\begin{flushright}
    $\Box$
\end{flushright}

\subsubsection{Proof of Lemma \ref{lem:tech_lemma_lower_bounds}}
Under our assumptions on $n$ and $x$, $\ln\left(1-x/n\right)$ is well-defined.
Using Taylor-Lagrange formula on the function ${[0,x] \ni t \mapsto\ln\left(1-t/n\right)}$ yields:
\begin{align*}
    &\left(1-\frac{x}{n}\right)^{n-1}= \exp\left((n-1)\ln\left(1-\frac{x}{n}\right)\right) \\
    &= \exp\left(-(n-1)\left(\frac{x}{n}+\frac{1}{2\left(1-\tau x/n\right)^2}\frac{x^2}{n^2}\right)\right)
\end{align*}
for some $\tau\in(0,1)$.
Using the fact that $\frac{n-1}{n}\leq 1$, $x\leq 1$ and
$\frac{1}{2\left(1-\tau x/n\right)^2}\leq\frac{1}{2\left(1-n^{-1}\right)^2}$,
we get that under our assumptions
$\left(1-\frac{x}{n}\right)^{n-1}
\geq \exp\left(-\left(1+\frac{1}{2n\left(1-n^{-1}\right)^2}\right)\right)$.
This bound is actually valid for every $x\in(0,1)$ and every $n\in\Nstar$.
The computation of
$\exp\left(-\left(1+\frac{1}{2n\left(1-n^{-1}\right)^2}\right)\right)$
shows that the latter is larger than $1/4$ whenever $n\geq 3$ and larger than $1/3$ whenever $n\geq 7$.
\begin{flushright}
    $\Box$
\end{flushright}

\section{Adapted results for ``Hoeffding'' framework}
\label{appendix:sec:results_hoef_framework}


\begin{hyp}
\label{hyp:hoef_dgp_XY}
    For every ${n \in \Nstar}$,
    there exist finite constants
    $\aXn$, $\bXn$, $\aYn$, $\bYn$ and $\lowerYn$ such that $X_{1,n}$ (respectively $Y_{1,n}$) lies $\distributionXYn$-almost surely in the interval $[\aXn , \bXn]$ (resp. $[\aYn , \bYn]$) and ${\left|\expecYn\right| \geq \lowerYn}$.
\end{hyp}

The support of $X_{1,n}$ and $Y_{1,n}$ is allowed to change with $n$, even though in many examples of interest, the former can be chosen independent from $n$. Assumptions~\ref{hyp:basic_dgp_XY} and \ref{hyp:hoef_dgp_XY} together correspond to the \textit{Hoeffding case} because under these two assumptions, we can use the Hoeffding inequality to build nonasymptotic CIs.

\subsection{Concentration inequality in an easy case: the support of the denominator is well-separated from \texorpdfstring{$0$}{0}}

\begin{hyp}\label{hyp:hoef_easy_dgp_well_separated}
    For every $n \in \Nstar$, the lower bound $\aYn$ is strictly positive.
\end{hyp}

\begin{thm}\label{thm:hoef_finite_ci_easy}
    Let $\upperXn:=\left(\bXn-\aXn\right)^2$ and $\upperYn:=\left(\bYn-\aYn\right)^2$. Under Assumptions~\ref{hyp:basic_dgp_XY}, \ref{hyp:hoef_dgp_XY} and \ref{hyp:hoef_easy_dgp_well_separated}, we have for every $n\in\INTST$ and $\varepsilon\in\Rstarplus$
    \begin{align*}
        \sup_{P \in \Pfam} \Prob_{P^{\otimes n}}  \Bigg( \bigg| \frac{\meanX}{\meanY} - \frac{\expecXn}{\expecYn} \bigg|
        &>\frac{\varepsilon}{\lowerYn}\left\{1+\frac{1}{\aYn}\left(\left|\aXn\right|\vee\left|\bXn\right|+\varepsilon\right)\right\}\Bigg) \\
        &\leq 4\exp\left(-\frac{2n\varepsilon^2}{\upperXn\vee\upperYn}\right).
    \end{align*}
    As a consequence, $\inf_{P \in \Pfam} \Prob_{P^{\otimes n}} \Big( \expecXn / \, \expecYn \in
    \left[\, \meanX /\,\meanY \pm t \right] \Big) \geq 1-\alpha$,
    with the following choice for $t$: $$\sqrt{\frac{\left(\upperXn\vee\upperYn\right)\ln\left(4/\alpha\right)}{2n {\lowerYn}^2}}\left(1+\frac{1}{\aYn}\left(\left|\aXn\right|\vee\left|\bXn\right|+\sqrt{\frac{\left(\upperXn\vee\upperYn\right)\ln\left(4/\alpha\right)}{2n}}\right)\right),$$ for every $\alpha \in (0,1)$.
\end{thm}


The theorem shows that it is possible to construct nonasymptotic CIs for ratios of expectations at every confidence level that are almost surely bounded.
However, it requires the additional Assumption \ref{hyp:hoef_easy_dgp_well_separated}, that in particular does not allow for binary $\{0,1\}$ random variables in the denominator which may limit its applicability for various applications.
In Section~\ref{sub:hoef_robust_ci_general_case}, we give an analogous result that only requires Assumptions~\ref{hyp:basic_dgp_XY} and \ref{hyp:hoef_dgp_XY} to hold, so that it encompasses the case of $\{0,1\}$-valued denominators.
However, the cost to pay will be an upper bound on the achievable coverage of the confidence intervals.

\subsection{Concentration inequality in the general case}
\label{sub:hoef_robust_ci_general_case}

We seek to build nontrivial nonasymptotic CIs under Assumptions~\ref{hyp:basic_dgp_XY} and \ref{hyp:hoef_dgp_XY} only. Under Assumption \ref{hyp:basic_dgp_XY}, $\expecYn\neq 0$, so that there is no issue in considering the fraction $\expecXn/\expecYn$. However, without Assumption~\ref{hyp:hoef_easy_dgp_well_separated}, $\left\{\meanY=0\right\}$ has positive probability in general so that $\meanX/\meanY$ is well-defined with probability less than one and undefined else. Note that when $P_{Y,n}$ is continuous wrt to Lebesgue's measure, there is no issue in defining $\meanX/\meanY$ anymore since the event $\left\{\meanY = 0\right\}$ has probability zero. This is not an easier case to establish concentration inequalities though, since without more restrictions, $\meanY$ can still be arbitrarily close to~$0$ with positive probability.

\begin{thm}
\label{thm:concentration_inequality_H_general_case}
    Assume that Assumptions~\ref{hyp:basic_dgp_XY} and~\ref{hyp:hoef_dgp_XY} hold.
    For every $n\in\INTST$, $\varepsilon > 0, \tilde \varepsilon \in (0,1)$,
    we have
    \begin{align*}
        \sup_{P \in \Pfam} \Prob_{P^{\otimes n}}  &\Bigg( \bigg| \frac{\meanX}{\meanY} - \frac{\expecXn}{\expecYn} \bigg|
        > \bigg( \frac{(|\aXn| \vee |\bXn|+\varepsilon) \tilde \varepsilon}{(1 - \tilde \varepsilon)^2} + \varepsilon \bigg) \frac{1}{\lowerYn}  \Bigg) \\
        &\hspace{4cm} \leq 2 \exp(-n \varepsilon^2 \gamma\left(X_{1,n}\right))
        + 2 \exp(-n \tilde \varepsilon^2 \gamma\left(Y_{1,n}\right)),
    \end{align*}
    where $\gamma\left(X_{1,n}\right) = 2 / (\bXn-\aXn)^2$ and
    $\gamma\left(Y_{1,n}\right) = 2 \lowerYn^2 / (\bYn-\aYn)^2$.

    As a consequence, $\inf_{P \in \Pfam} \Prob_{P^{\otimes n}} \Big( \expecXn / \, \expecYn \in
    \left[\, \meanX /\,\meanY \pm t \right] \Big) \geq 1-\alpha$,
    with the choice
    \begin{align*}
        t
        %
        &:= \sqrt{\frac{\ln(4/\alpha) }{n \gamma\left(X_{1,n}\right) \wedge \gamma\left(Y_{1,n}\right)}}
        \bigg( \frac{ |\aXn| \vee |\bXn| + \sqrt{ \ln(4/\alpha) / \left(n \gamma\left(X_{1,n}\right)\right) } }
        {\Big(1 - \sqrt{ \ln(4/\alpha) / \left(n \gamma\left(Y_{1,n}\right)\right) }\Big)^2}
        + 1 \bigg) \frac{1}{\lowerYn},
    \end{align*}
    for every $\alpha > \overline{\alpha}_{n,H} := 4 e^{-n \gamma\left(Y_{1,n}\right)}$.\footnote{Equivalently, it means that for a given level $\alpha$, the choice of $t$ is valid for every integer~$n > \overline{n}_{\alpha,H} := \ln(4/\alpha) / \gamma\left(Y_{1,n}\right)$.}
    \label{thm:hoef_finite_ci}
\end{thm}
This theorem is proved in Section~\ref{proof:thm:hoef_finite_ci}. It states that when $\lowerYn> 0$, it is possible to build valid nonasymptotic CIs with finite length up to the confidence level $1-\overline{\alpha}_{n,H}$. This is a more positive result than \cite{dufour1997} which claims that it is not possible to build nontrivial nonasymptotic CIs when $\lowerYn$ is taken equal to 0, no matter the confidence level. 
Note that Theorem~\ref{thm:hoef_finite_ci} is not an impossibility theorem since it only claims that considering confidence levels smaller than $1-\overline{\alpha}_{n,H}$ is \textit{sufficient} to build nontrivial CIs under Assumptions~\ref{hyp:basic_dgp_XY} and \ref{hyp:hoef_dgp_XY}. The remaining question is to find out whether it is \textit{necessary} to focus on confidence levels that do not exceed a certain threshold under Assumptions~\ref{hyp:basic_dgp_XY} and~\ref{hyp:hoef_dgp_XY}. We answer this in Section~\ref{subsec:hoef_upper_bound_level}.

Theorem~\ref{thm:hoef_finite_ci} has two other interesting consequences: for every confidence level up to $1-\overline{\alpha}_{n,H}$, a nonasymptotic CI of the form $\left[\meanX/\meanY\pm \tilde{t}\right]$ with $\tilde{t}>t$ has good coverage but is too conservative. What is more, if the DGP does not depend on $n$ (\textit{i.e} in the standard i.i.d. set-up), for every fixed $\alpha>\overline{\alpha}_{n,H}$, the length of the confidence interval shrinks at the optimal rate $1/\sqrt{n}$.


\subsection{An upper bound on testable confidence levels}\label{subsec:hoef_upper_bound_level}

\begin{thm}
    \label{thm:hoef_necessary_cond_alpha}
	For every $n\in\INTST$, and
	every $\alpha \in \left(0, \underline{\alpha}_{n,H} \right),$
	where $\underline{\alpha}_{n,H} := \big(1-l_{Y,n}/(\bYn-\aYn)\big)^n$,
	if~$(\bYn-\aYn)/\lowerYn> 1$, there is no finite $t>0$ such that $\left[\meanX/\meanY\pm t\right]$ has coverage $1-\alpha$ over $\Pfam_H$, where $\Pfam_H$ is the class of all distributions satisfying Assumptions~\ref{hyp:basic_dgp_XY} and~\ref{hyp:hoef_dgp_XY} for a fixed lower bound $\lowerYn$ and fixed lengths $\bXn-\aXn$ and $\bYn-\aYn$.
\end{thm}



This theorem asserts that confidence intervals of the form $\left[\meanX/\meanY\pm t\right]$ with coverage higher than $1-\underline{\alpha}_{n,H}$ under Assumptions~\ref{hyp:basic_dgp_XY} and \ref{hyp:hoef_dgp_XY} are not defined (or are of infinite length) with positive probability for at least one distribution in $\Pfam_H$. The additional restriction $(\bYn-\aYn)/\lowerYn> 1$ is rather mild in practice: it is equivalent to $\bYn-\aYn>\lowerYn$ and is satisfied as soon as $\aYn\leq 0$ and $\bYn>\lowerYn>0$. This encompasses all DGPs where the denominator is $\{0,1\}$-valued and the probability that the denominator equals 1 is bounded from below by $\lowerYn\in(0,1)$.

\medskip

Note that for Theorems~\ref{thm:hoef_finite_ci_easy} and~\ref{thm:concentration_inequality_H_general_case}, it is required to know not only the length ${\bXn - \aXn}$ but also the actual endpoints of the support, $\aXn$ and $\bXn$.
On the contrary, Theorem~\ref{thm:hoef_necessary_cond_alpha} does not require the latter.
In that respect, the class of Theorem~\ref{thm:hoef_necessary_cond_alpha} is larger than the one of the two preceding theorems.

\subsection{Proof of Theorems~\ref{thm:hoef_finite_ci_easy} and \ref{thm:hoef_finite_ci}}
\label{proof:thm:hoef_finite_ci}

The proofs are identical to those of Theorems~\ref{thm:bt_finite_ci_easy} and \ref{thm:bt_finite_ci}, except for the Bienaymé-Chebyshev inequality that has to be replaced with the Hoeffding inequality. The latter can be used under Assumption~\ref{hyp:hoef_dgp_XY}. Note also that $\expecXn$ is now bounded by $|\aXn| \vee |\bXn|$.
\begin{flushright}
    $\Box$
\end{flushright}

\subsection{Proof of Theorem~\ref{thm:hoef_necessary_cond_alpha}}
\label{proof:thm:hoef_necessary_cond_alpha}

We need the subsequent lemma.

\begin{lemma}
    \label{lemma:hoef:lowerBound_probInfty}
    For each $\xi$ in the interval
    $\Big(0 , 1 \wedge \big((\bYn-\aYn)/\lowerYn - 1\big) \Big)$, there exists a distribution $P_{n, \, \xi} \in \Pfam_H$ such that
    $\Prob\left( \meanY = 0\right) \geq \tilde \alpha_{n,H}(\xi),$
    where $\tilde \alpha_{n,H}(\xi)
    := \big(1 - (1+\xi) \lowerYn / (\bYn-\aYn) \big)^n$.
\end{lemma}
\noindent
Note that the interval
$\Big(0 , 1\wedge\big((\bYn-\aYn)/\lowerYn-1\big)\Big)$ is not empty since we have assumed $(\bYn-\aYn)/\lowerYn > 1$.

By Lemma~\ref{lemma:hoef:lowerBound_probInfty}, for every
$\xi < 1\wedge\big((\bYn-\aYn)/\lowerYn - 1\big)$,
there exists a distribution $P_{n, \, \xi} \in \Pfam_H$ satisfying Assumptions~\ref{hyp:basic_dgp_XY} and~\ref{hyp:hoef_dgp_XY} such that
$\Prob\left( \meanY = 0 \right)
\geq \tilde \alpha_{n,H}(\xi)$.
Denote its marginal distributions by $P_{X, n, \, \xi}$
and $P_{Y, n, \, \xi}$.
Therefore, $P_{n, \, \xi}$ satisfies Assumptions~\ref{hyp:basic_dgp_XY} and \ref{hyp:hoef_dgp_XY}, and $\meanX/\meanY$ is undefined with probability greater than $\tilde \alpha_{n,H}(\xi)$.
Taking the supremum over $\xi$, we deduce that
\begin{align*}
    \sup_{P_n \in \Pfam_H}
    \Prob \left( \meanY = 0 \right)
    \geq \sup_\xi \tilde \alpha_n(\xi) = \underline{\alpha}_{n,H}.
\end{align*}
This means that the random interval
$I_n^* := \left[\meanX/\meanY\pm t\right]$ cannot have coverage higher than $1-\underline{\alpha}_{n,H}$ since it may be undefined with a probability higher than $\underline{\alpha}_{n,H}$.
\begin{flushright}
    $\Box$
\end{flushright}

\subsubsection{Proof of Lemma~\ref{lemma:hoef:lowerBound_probInfty}}
We consider the following distribution on $\Rb$
\begin{align*}
    &P_{n,\lowerYn,c,\xi}
    :=\left(\frac{c}{n}\right)^{1/n} \delta_{\lsin 0\rsin}
    + \frac{1}{2}\left(1-\left(\frac{c}{n}\right)^{1/n}\right)
    \delta_{\{y_{c-}\}}
    + \frac{1}{2}\left(1-\left(\frac{c}{n}\right)^{1/n}\right)
    \delta_{\{y_{c+}\}},
\end{align*}
where $c\in(0,n)$ is some constant to be chosen later,
$y_{c-} := \lowerYn(1-\xi) / (1-(c/n)^{1/n})$
and $y_{c+} := \lowerYn(1+\xi) / (1-(c/n)^{1/n})$.
Let $Y_{1,n}\sim P_{n,\lowerYn,c,\xi_n}$.
Observe that $\expecYn=\lowerYn$. With the choice
\begin{align*}
    &c = c_n :=
    n \left(1-\frac{\lowerYn}{\bYn-\aYn}\left(1+\xi\right)\right)^n,
\end{align*}
we have $y_{c+} = \bYn-\aYn$.
Note that $\Cnalpha$ is strictly positive,
because $1-\frac{\lowerYn}{\bYn-\aYn}\left(1+\xi_n\right) > 0$.
This is equivalent to $\bYn-\aYn / \lowerYn > 1 + \xi_n$, which is true by assumption.

\medskip

Consider now the following product measure on $\Rb^2$
defined by $P_n:=
\big(0.5 \delta_{\lsin 0 \rsin}
+ 0.5 \delta_{\lsin \bXn - \aXn \rsin} \big)
\otimes P_{n,\lowerYn,c_n,\xi}.$
Let $(X_{i,n},Y_{i,n})_{i=1}^n\simiid P_n$.
These random vectors satisfy $\expecYn=\lowerYn$, $(\max - \min)[Y_{1,n}] = \bYn - \aYn$ and $(\max - \min)[X_{1,n}] = \bXn - \aXn$.
The next step is to build a lower bound on the event $\{ \meanY = 0 \}$.

The assumption that $(X_{i,n},Y_{i,n})_{i=1}^n\simiid P_n$ and the construction of $P_{n,\lowerYn,c_n,\xi}$ imply that
\begin{align*}
    \Prob\left( \meanY=0 \right)
    = \frac{c_n}{n}
    = \left(1-\frac{\lowerYn}{\bXn - \aXn}\left(1+\xi\right)\right)^n
    = \tilde \alpha_{n,H}(\xi).
\end{align*}
\begin{flushright}
    $\Box$
\end{flushright}

\section{Additional simulations}
\label{appendix:sec:additional_simulations}

This section complements the simulations presented in the main body of the article using different distributions for the variables in the numerator and in the denominator.

\medskip
In this setting of simulations, we use the best bounds by setting the constants $\lowerYn$ and $\upperYn$ that define our class of distributions equal to the actual corresponding moments (respectively the expectation for $\lowerYn$ and the second moment 
for $\upperYn$).
That is we use ${\overline{n}_\alpha = 2\VarY/(\alpha \expecY^2)}$ or ${\overline{\alpha}_n = 2\VarY/(n\expecY^2)}$.
In practice, our rule-of-thumb uses the plug-in version of those quantities replacing the theoretical unknown moments by their empirical counterparts as explained in Section~\ref{subsec:plug_in_estimators}.

\medskip

The following Figures are similar to Figures~\ref{fig:gr_newa_gaussians_n_alphapc_10nbrep_5000EYpc_10EXpc_50VXpc_100VYpc_200CorrXYpc_50dilanpc_130} and~\ref{fig:gr_gr_newa_gaussian_alpha_n_nbrep10000adila_3n_1000EY_025EX_05VX_2VY_1CorrXY_05}.
They show the $c(n,P)$ of the asymptotic CIs based on the delta method as a function of the sample size~$n$ and also reports $\overline{n}_\alpha :=2\left(\upperYn-{\lowerYn}^2\right)/\left(\alpha\lowerYn^2\right)$, with $\alpha$ chosen according to the desired nominal level (equal to~${1-\alpha}$) and  $\lowerYn = \expecY$, $\upperYn = {\expecY}^2 + \VarY$.
Consequently, the titles of the figures only indicate the specification used for~$\distributionXYn$, the nominal pointwise asymptotic level~${1-\alpha}$, and the number of repetitions used to approximate the probability~$c(n,P)$.

\medskip

With discrete distributions for the variable in the denominator, it may happen that ${\meanY = 0}$, all the more so as the expectation and the sample size are low typically.
As discussed at the end of Section~\ref{sec:framework}, confidence intervals are said to be undefined when ${\meanY = 0}$.
In such cases, for any value~${a \in \Rb}$, it is undefined whether $a$ belongs or not to the CIs.
Consequently, whenever the sample drawn is such that $\meanY = 0$ in the simulations, we count the draw as a no coverage occurrence in the Monte Carlo estimation of~$c(n,P)$.
In other words, this quantity is approximated as an average over~$M$ repetitions and the repetitions for which ${\meanY = 0}$ account for~$0$ in this average.\footnote{Note that in some specifications, a substantial part of the repetitions yield $\meanY = 0$. For instance, with Bernoulli distributions, for $n$ smaller than $10$ and the expectation at the denominator equal to $0.01$, around 10\% only of the repetitions display $\meanY \neq 0$.}

\subsection{Gaussian distributions}


\begin{figure}[ht]
    \centering
    \includegraphics[width=0.6\textwidth]{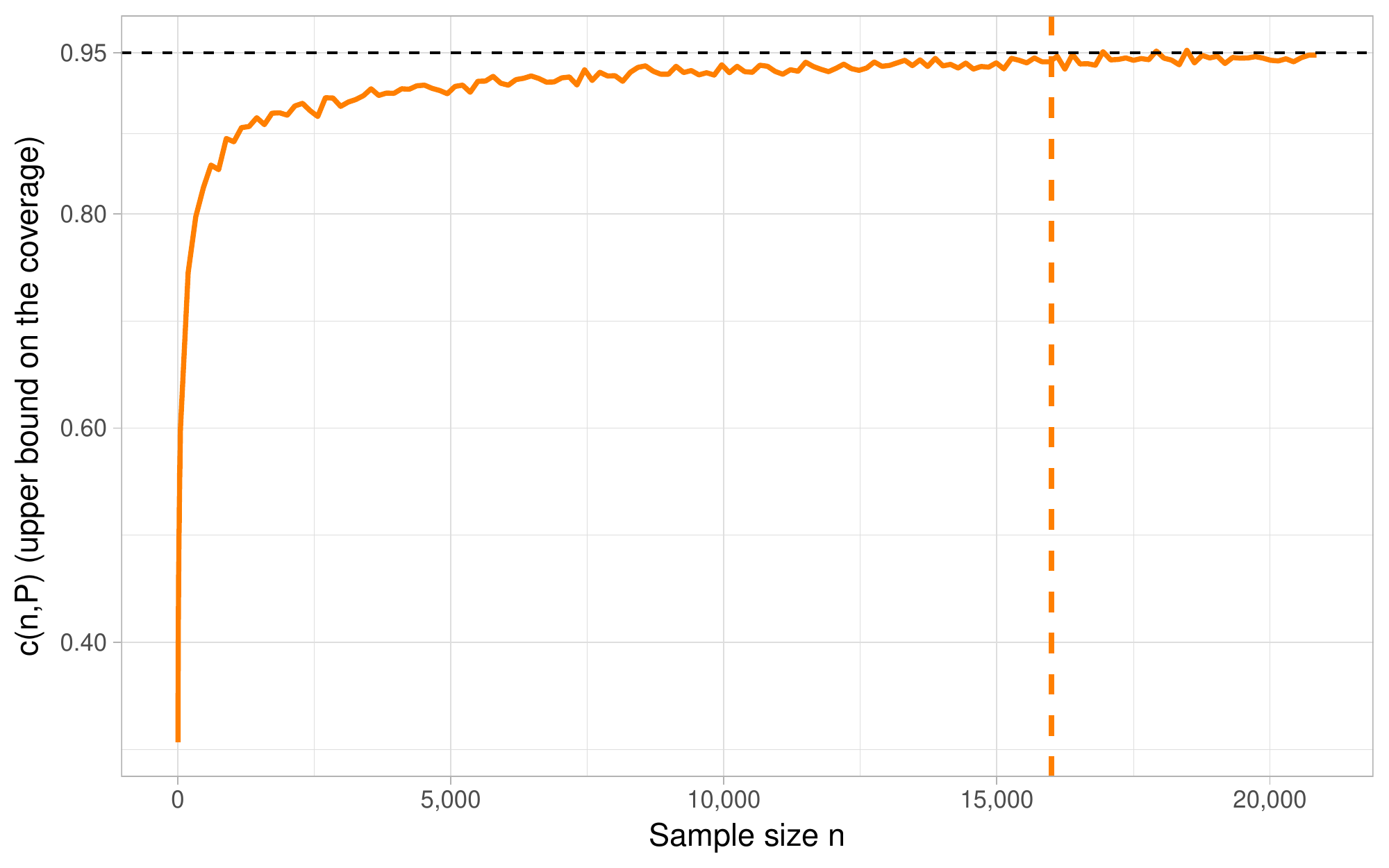}
    \caption{\small{Specification: $\forall n \in \Nstar, {\distributionXYn = \mathcal{N}(1,1) \otimes \mathcal{N}(0.05,1)}$;
    $1-\alpha=0.95$;
    5,000 repetitions used.}}
    \label{fig:appendix:gr_newa_gaussian_n_alpha_nbrep5000EY_005EX_1VX_1VY_1CorrXY_0}
\end{figure}


\begin{figure}[ht]
    \centering
    \includegraphics[width=0.6\textwidth]{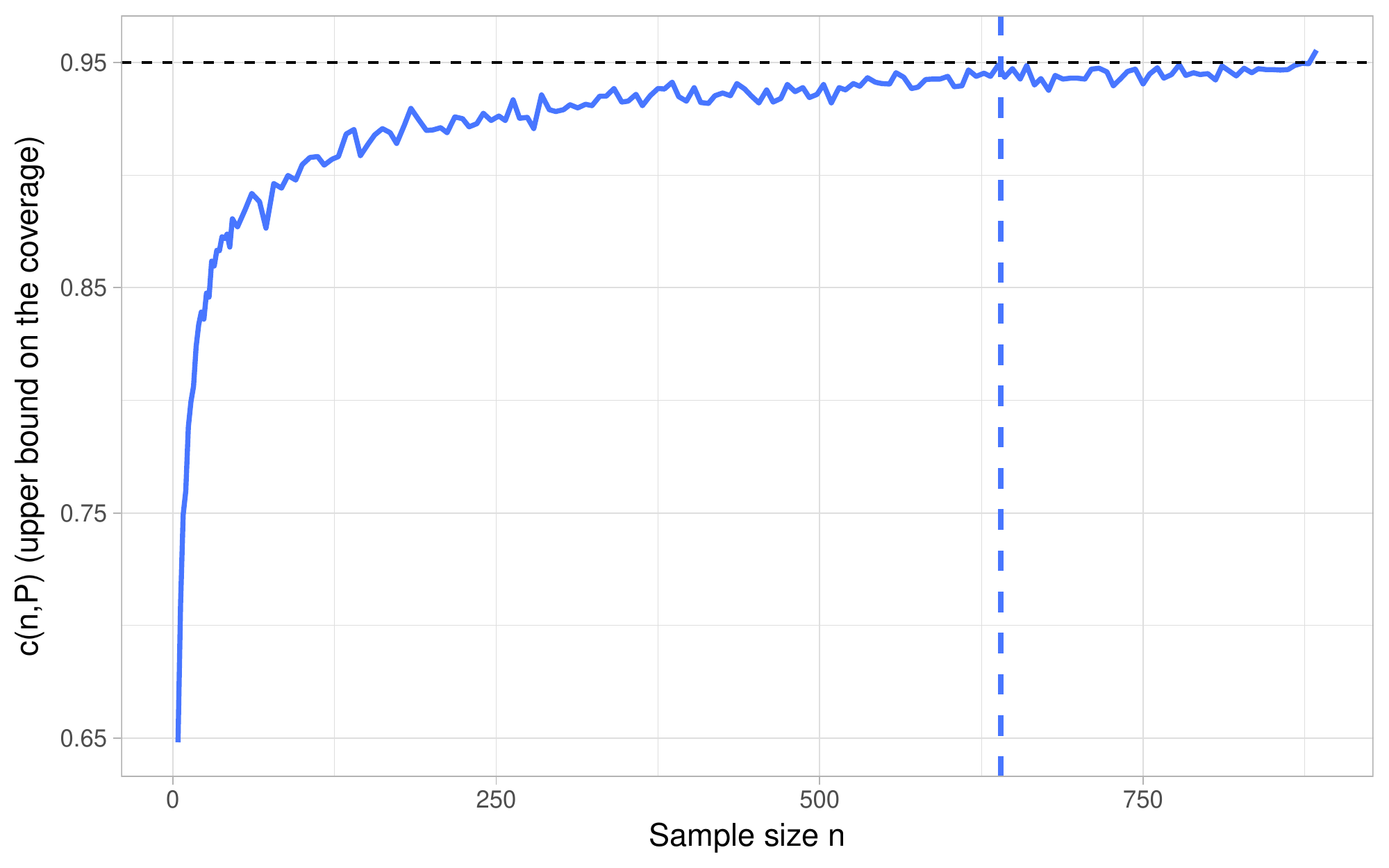}
    \caption{\small{Specification: $\forall n \in \Nstar, {\distributionXYn = \mathcal{N}(1,1) \otimes \mathcal{N}(0.25,1)}$;
    $1-\alpha=0.95$;
    5,000 repetitions used.}}
    \label{fig:appendix:gr_newa_gaussian_n_alpha_nbrep5000EY_025EX_1VX_1VY_1CorrXY_0}
\end{figure}


\begin{figure}[ht]
    \centering
    \includegraphics[width=0.6\textwidth]{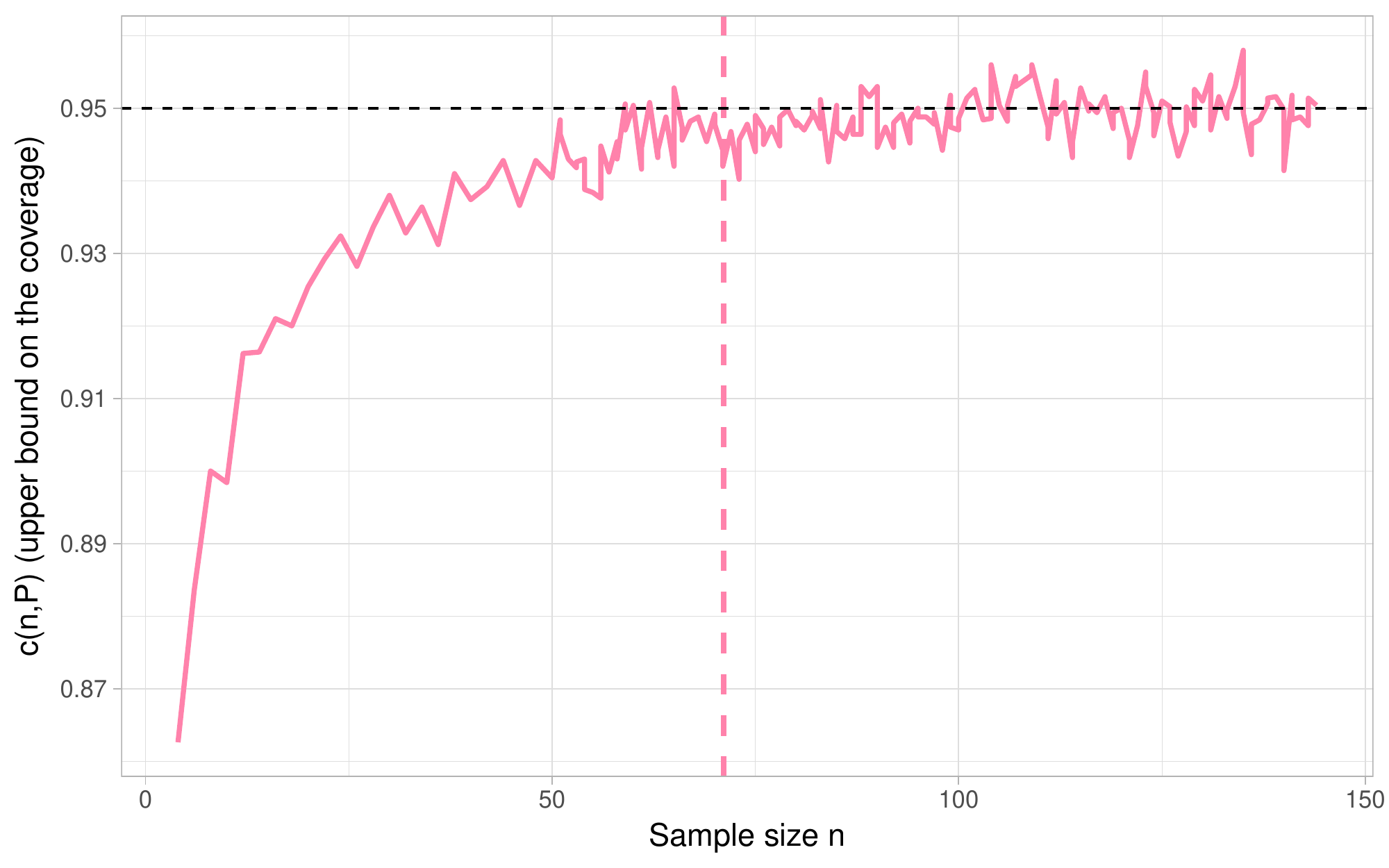}
    \caption{\small{Specification: $\forall n \in \Nstar, {\distributionXYn = \mathcal{N}(1,1) \otimes \mathcal{N}(0.75,1)}$;
    $1-\alpha=0.95$;
    5,000 repetitions used.}}
    \label{fig:appendix:gr_newa_gaussian_n_alpha_nbrep5000EY_075EX_1VX_1VY_1CorrXY_0}
\end{figure}


\begin{figure}[ht]
    \centering
    \includegraphics[width=0.6\textwidth]{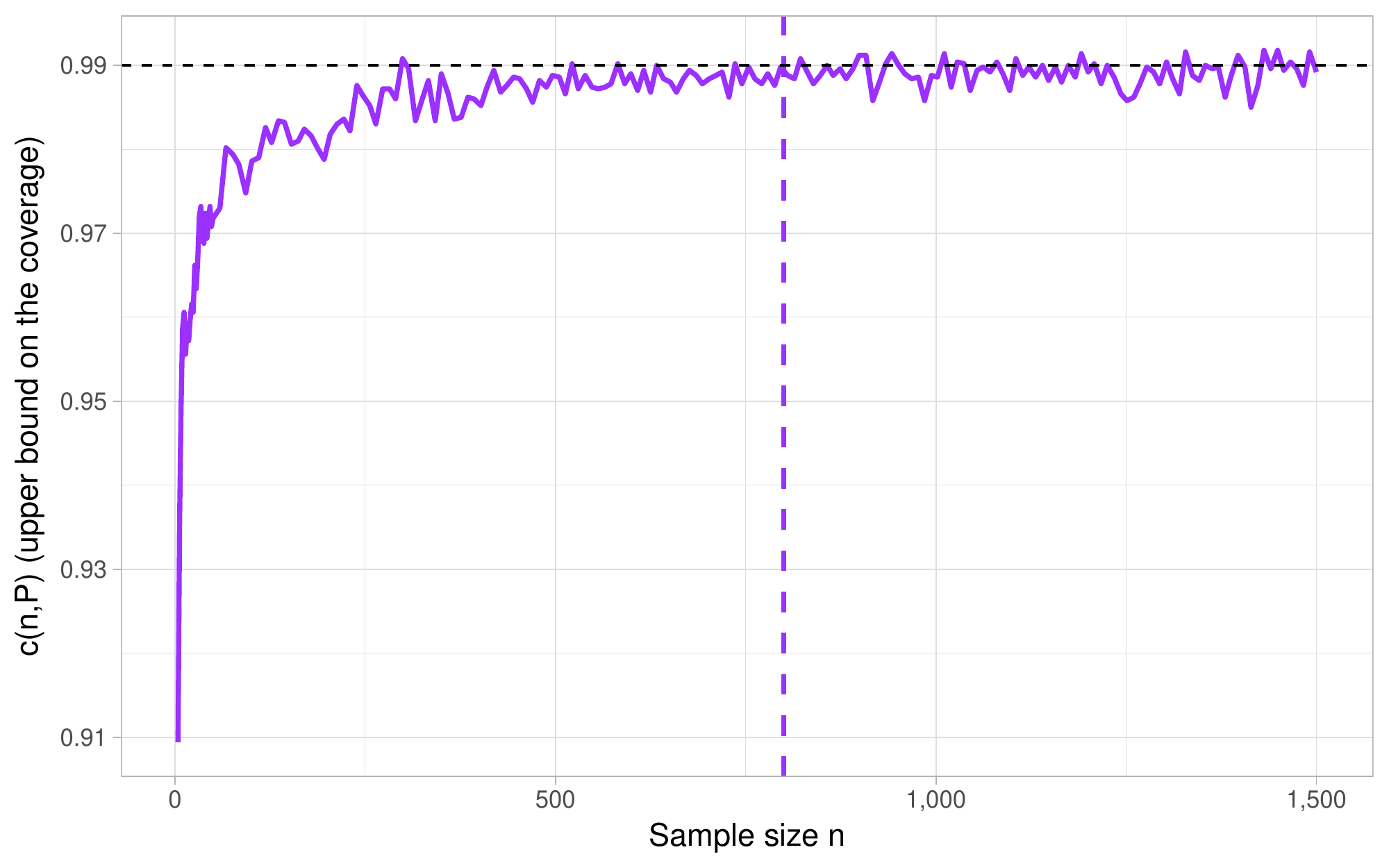}
    \caption{\small{Specification: $\forall n \in \Nstar, {\distributionXYn = \mathcal{N}_2}$ bivariate Gaussian with $\expecX = 0.5$, $\expecY = 0.5$, $\VarX = 2$, $\VarY = 1$, $\Corr(X,Y) = -0.3$;
    $1-\alpha=0.99$;
    5,000 repetitions used.}}    
    \label{fig:appendix:gr_newa_gaussians_n_alphapc_1nbrep_5000EYpc_50EXpc_50VXpc_200VYpc_100CorrXYpc_-30dilanpc_150}
\end{figure}




\begin{figure}[ht]
    \centering
    \includegraphics[width=0.6\textwidth]{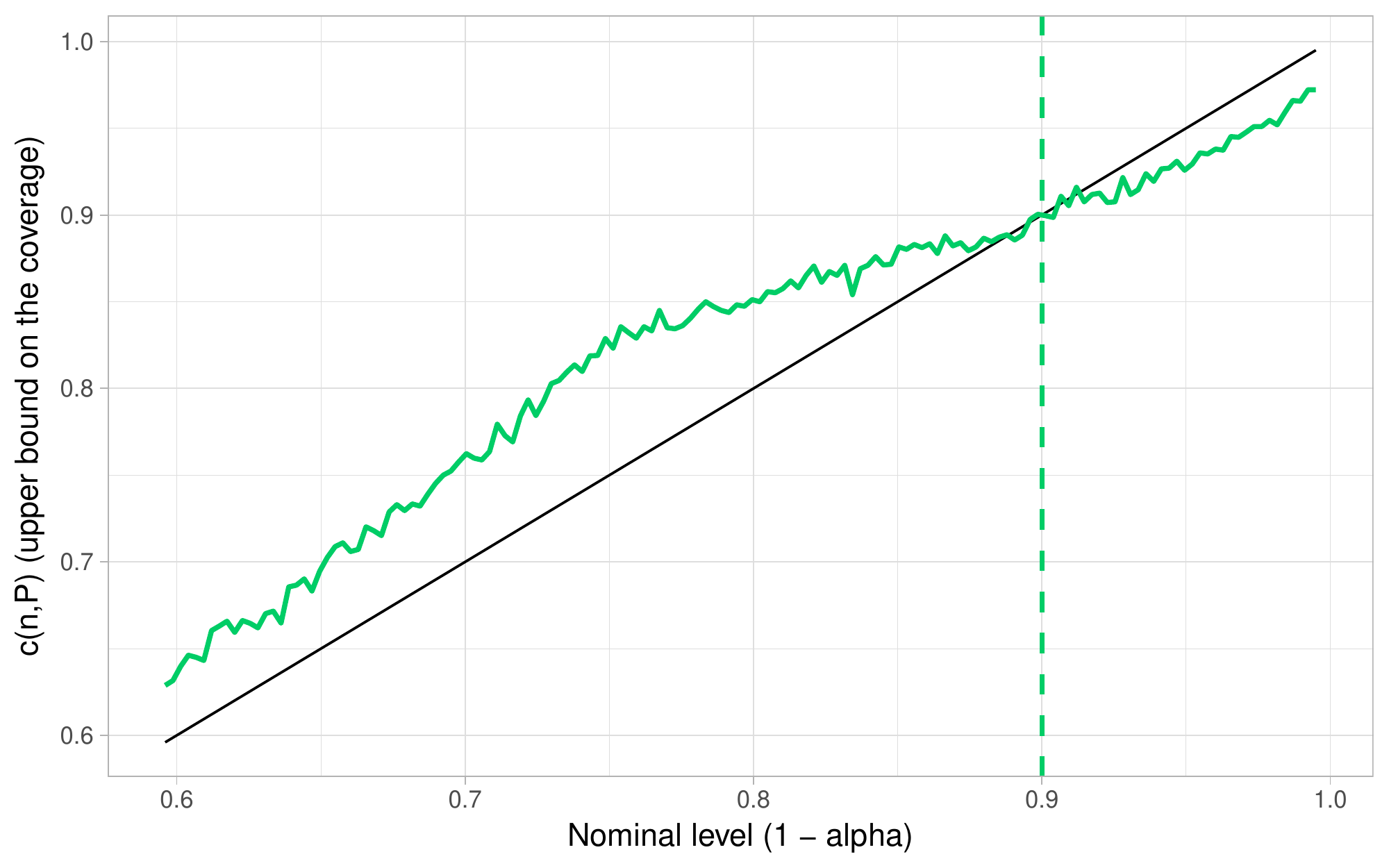}
    \caption{\small{Specification: $\forall n \in \Nstar, {\distributionXYn = \mathcal{N}(1,1) \otimes \mathcal{N}(0.1,1)}$;
    $n = 2,000$;
    5,000 repetitions used.}}  
    \label{fig:appendix:gr_gr_newa_gaussian_alpha_n_nbrep5000adila_4n_2000EY_01EX_1VX_1VY_1CorrXY_0}
\end{figure}

\subsection{Student distributions}

The specification here is two Student distributions, both in the numerator and in the denominator. 
Standard Student distributions are centered.
We use therefore translated versions by simply adding the expectations in order to avoid a null denominator for the ratio of expectations of interest.
Below, $\mathcal{T}(\mu, \nu)$ denotes the distribution of a translated standard Student variable: $\mu + T$ where $T$ is distributed according to a Student distribution with $\nu$ degrees of freedom.
To satisfy Assumption \ref{hyp:basic_dgp_XY}, we need finite variance: we use degrees of freedom strictly higher than $2$ for this purpose.

\begin{figure}[ht]
    \centering
    \includegraphics[width=0.6\textwidth]{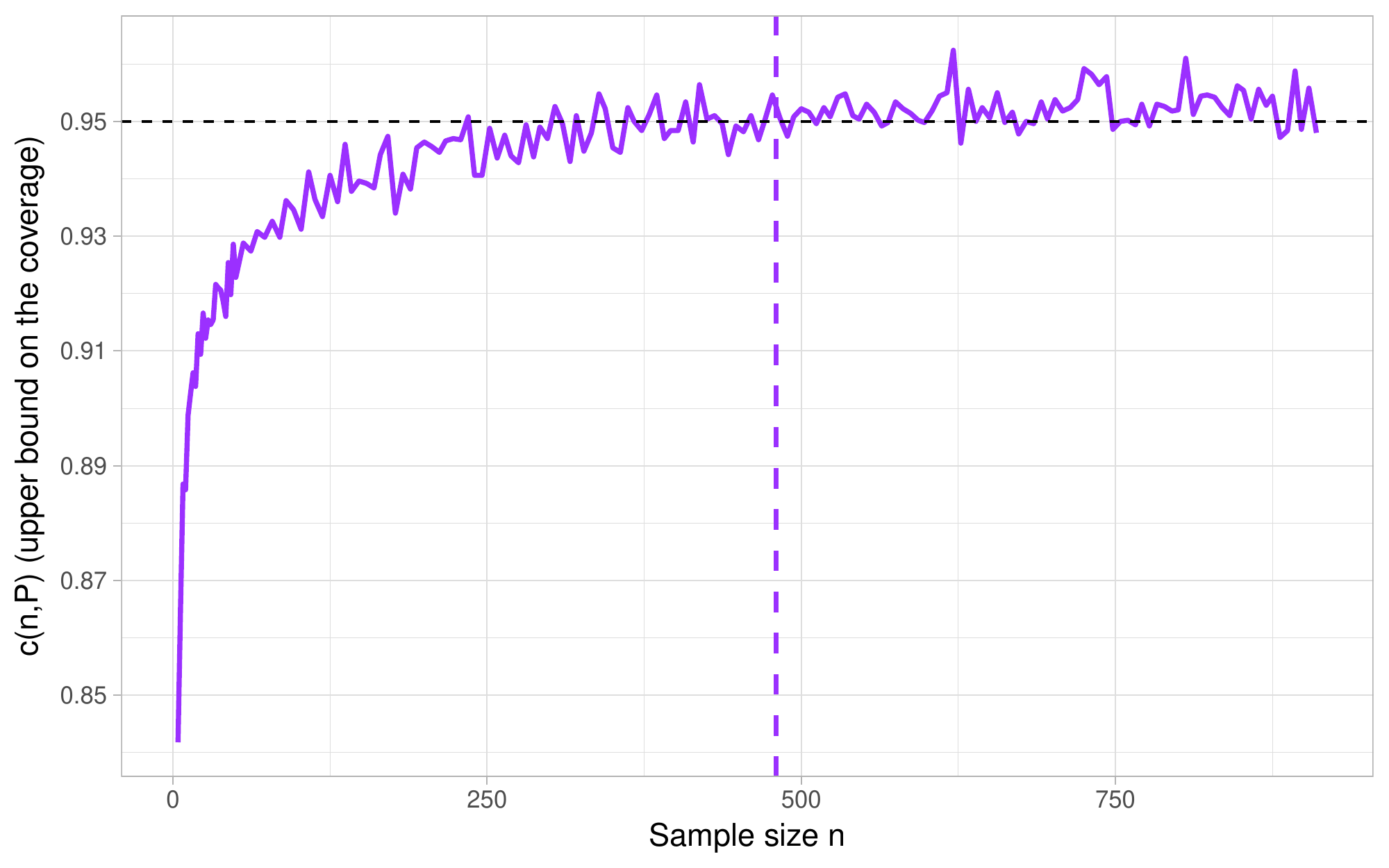}
    \caption{\small{Specification: $\forall n \in \Nstar, {\distributionXYn = \mathcal{T}(0.5, 3) \otimes \mathcal{T}(0.5, 3)}$;
    $1-\alpha=0.95$;
    5,000 repetitions used.}}    
    \label{fig:appendix:gr_newa_students_n_alphapc_5nbrep_5000EYpc_50dfY_3EXpc_50dfX_3CorrXYpc_0}
\end{figure}


\begin{figure}[ht]
    \centering
    \includegraphics[width=0.6\textwidth]{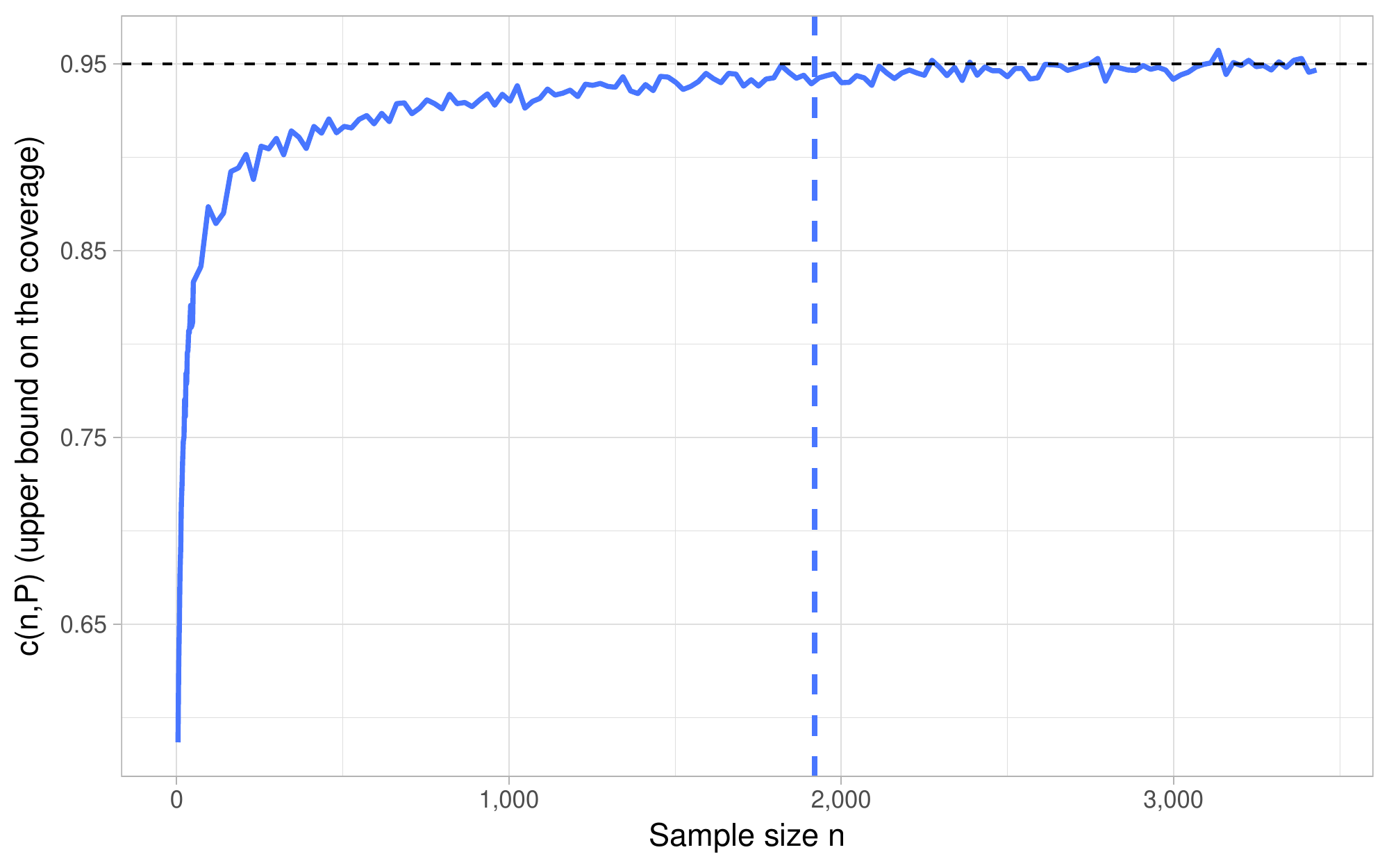}
    \caption{\small{Specification: $\forall n \in \Nstar$, the marginal distributions of $X$ and $Y$ are $\mathcal{T}(1, 3)$ and $\mathcal{T}(0.25, 3)$ respectively and are simulated using a Gaussian copula to have ${\Corr(X,Y) \approx 0.5}$.;
    $1-\alpha=0.95$;
    5,000 repetitions used.}} 
    \label{fig:appendix:gr_newa_students_n_alphapc_5nbrep_5000EYpc_25dfY_3EXpc_100dfX_3CorrXYpc_55}
\end{figure}

\subsection{Exponential distributions}
\label{appendix:sec:exponentials}

The specification here is two exponential distributions, both in the numerator and in the denominator.
The case of the exponential is specific as a unique parameter determines both the expectation and the variance of the distribution.

\begin{figure}[htb]
    \centering
    \includegraphics[width=0.6\textwidth]{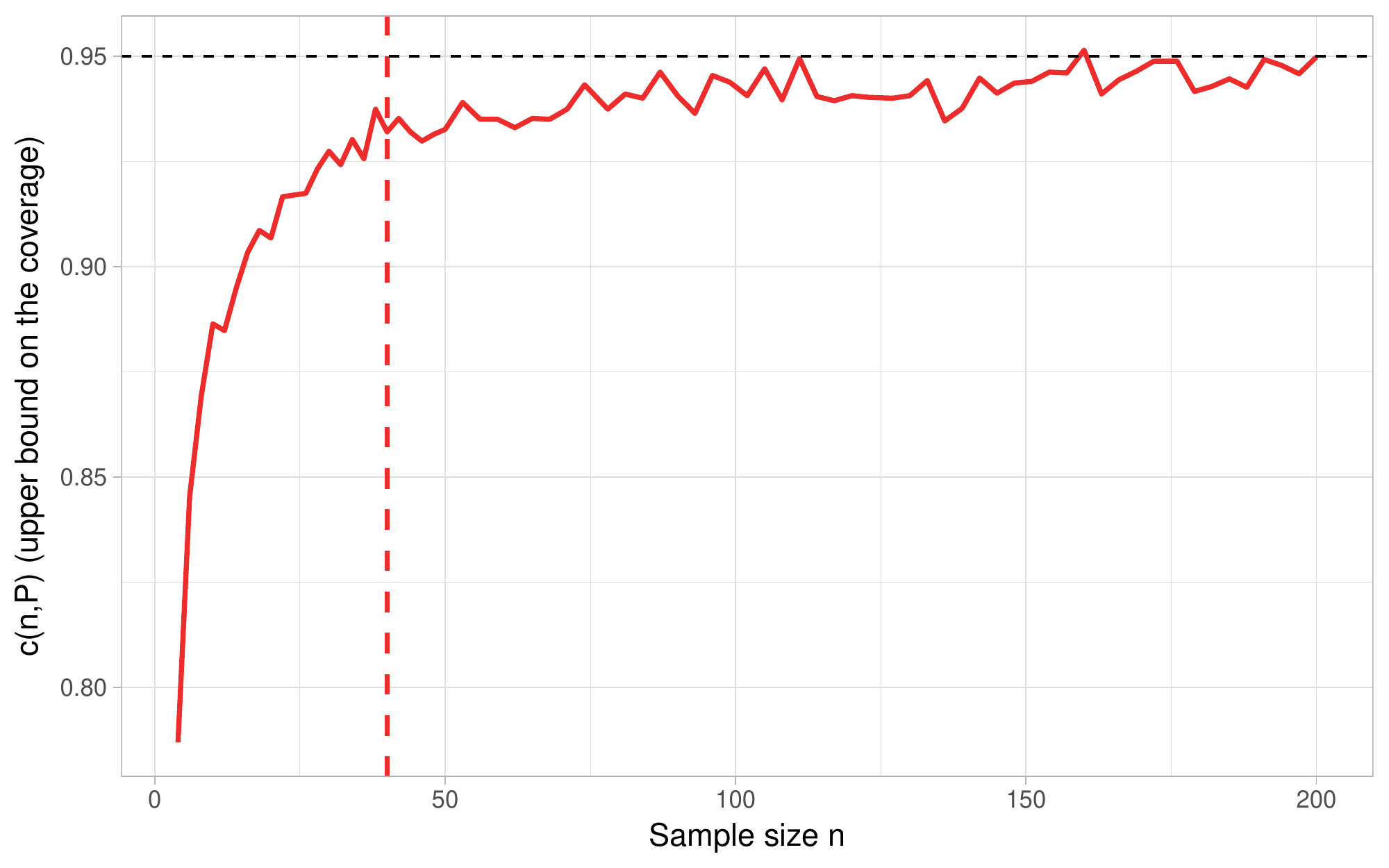}
    \caption{\small{Specification: $\forall n \in \Nstar, {\distributionXYn = \mathcal{E} \otimes \mathcal{E}}$ with ${\expecX=1}$ and ${\expecY=0.01}$;
    $1-\alpha=0.95$;
    5,000 repetitions used.}}      
    \label{fig:appendix:gr_newa_exponentials_n_alphapc_5nbrep_5000EYpc_1EXpc_100CorrXYpc_0}
\end{figure}

More precisely, the variance is equal to the square of the expectation.
Consequently, whatever the parameter of the exponential distribution in the denominator, we have ${\overline{n}_\alpha = 4/\alpha}$.
Previous simulations suggest that the closer the expectation in the denominator to $0$, the larger the sample size required for the asymptotic approximation to hold.
At first sight, we might thus be worried for the usefulness of our rule-of-thumb to obtain $\overline{n}_\alpha$ independent of $\expecY$.
Yet, with exponential distributions, the lower the expectation, the lower is the variance too. 
Intuitively, the lower variance will compensate having an expectation closer to~$0$.
The previous statement that links the closeness to $0$ of the expectation in the denominator and the sample size required to reach the asymptotic approximation presupposes keeping fixed the variance.
It cannot be anymore for exponential distributions. 

\begin{figure}[htb]
    \centering
    \includegraphics[width=0.6\textwidth]{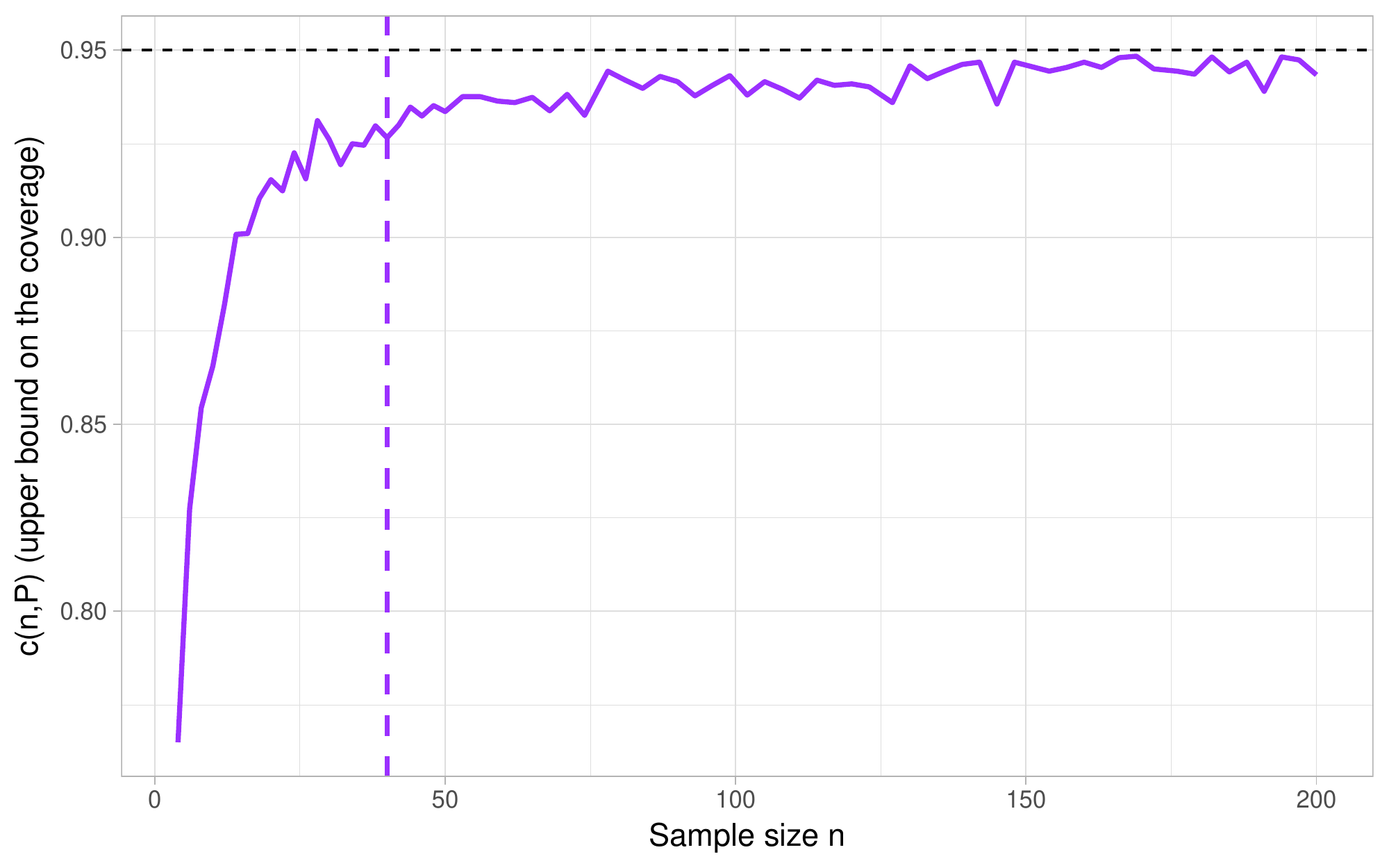}
    \caption{\small{Specification: $\forall n \in \Nstar$, the marginal distributions of $X$ and $Y$ are two exponentials with ${\expecX = 1}$ and ${\expecY = 0.5}$ and are simulated using a Gaussian copula to have ${\Corr(X,Y) \approx 0.75}$.;
    $1-\alpha=0.95$;
    5,000 repetitions used.}} 
    \label{fig:appendix:gr_newa_exponentials_n_alphapc_5nbrep_5000EYpc_50EXpc_100CorrXYpc_80}
\end{figure}

\medskip
The simulations reveal that the convergence of the coverage of the asymptotic confidence intervals toward their nominal level happens for $n$ around one hundred fifty and has the same pattern whatever the expectation of the exponential distribution in the denominator.
Our rule-of-thumb $\overline{n}_\alpha$ appears to be a bit small. Nonetheless, it is coherent that it is constant across the value of $\expecY$.

\subsection{Pareto distributions}

 The specification here is two Pareto distributions, both in the numerator and in the denominator.
Pareto distributions have support in $\Rstarplus$. 
They would fall in the easier case when the support of the denominator is well separated from $0$.
To assess the dependability of our rule-of-thumb in the general case, we use translated Pareto distributions.
In what follows, the notation $\mathcal{P}areto(\expecY, \tau, \gamma)$ denotes the distribution of a random variable that follows 
a Pareto distribution with shape parameter equal to~$\gamma$ translated such that its support is~${(\tau, +\infty)}$ and its expectation is~$\expecY$.
A variable that is distributed according to  $\mathcal{P}areto(\expecY, \tau, \gamma)$ is equal in distribution to $P + (\expecY - \gamma t_Y)/(\gamma-1)$ with $t_Y = (\expecY - \tau) \times (\gamma-1)$ and $P$ a usual Pareto distribution with support or scale parameter $t_Y$ and shape parameter $\gamma$, that is $P$ has the density $x \mapsto \Indicator\{x \geq t_Y\} \times \gamma t_Y^{\gamma} / x^{\gamma+1}$ with respect to Lebesgue measure.

\begin{figure}[ht]
    \centering
    \includegraphics[width=0.6\textwidth]{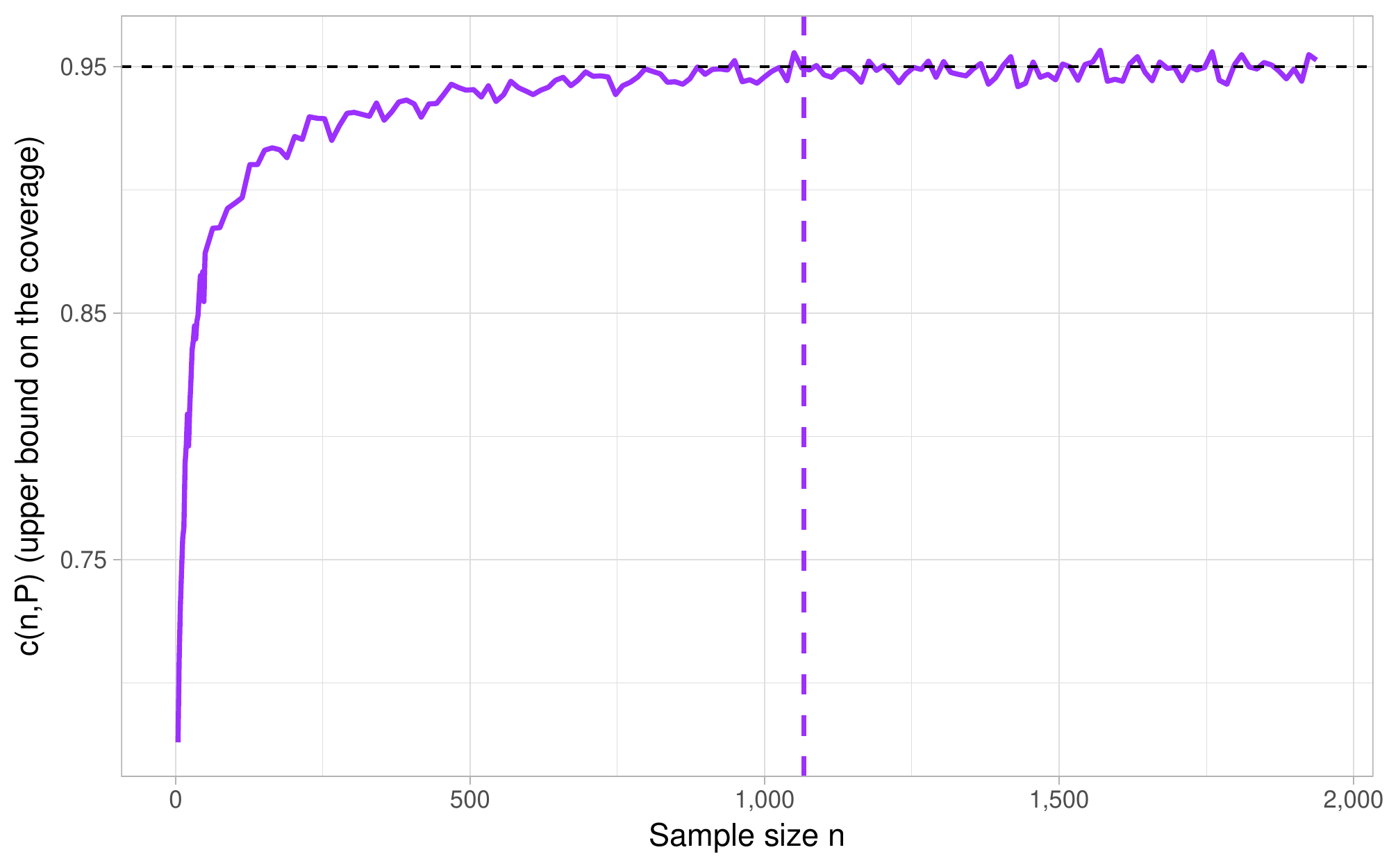}
    \caption{\small{Specification: $\forall n \in \Nstar, {\distributionXYn = \mathcal{P}areto(1, -1.5, 5) \otimes \mathcal{P}areto(\expecY, -1.5, 5)}$, with $\expecY = 0.5$;
    $1-\alpha=0.95$;
    5,000 repetitions used.}}       
    \label{fig:appendix:gr_newa_paretos_n_alphapc_5nbrep_5000EYpc_50thrsYpc_-150shpY_5EXpc_100thrsXpc_-150shpX_5CorrXYpc_0}
\end{figure}


\begin{figure}[ht]
    \centering
    \includegraphics[width=0.6\textwidth]{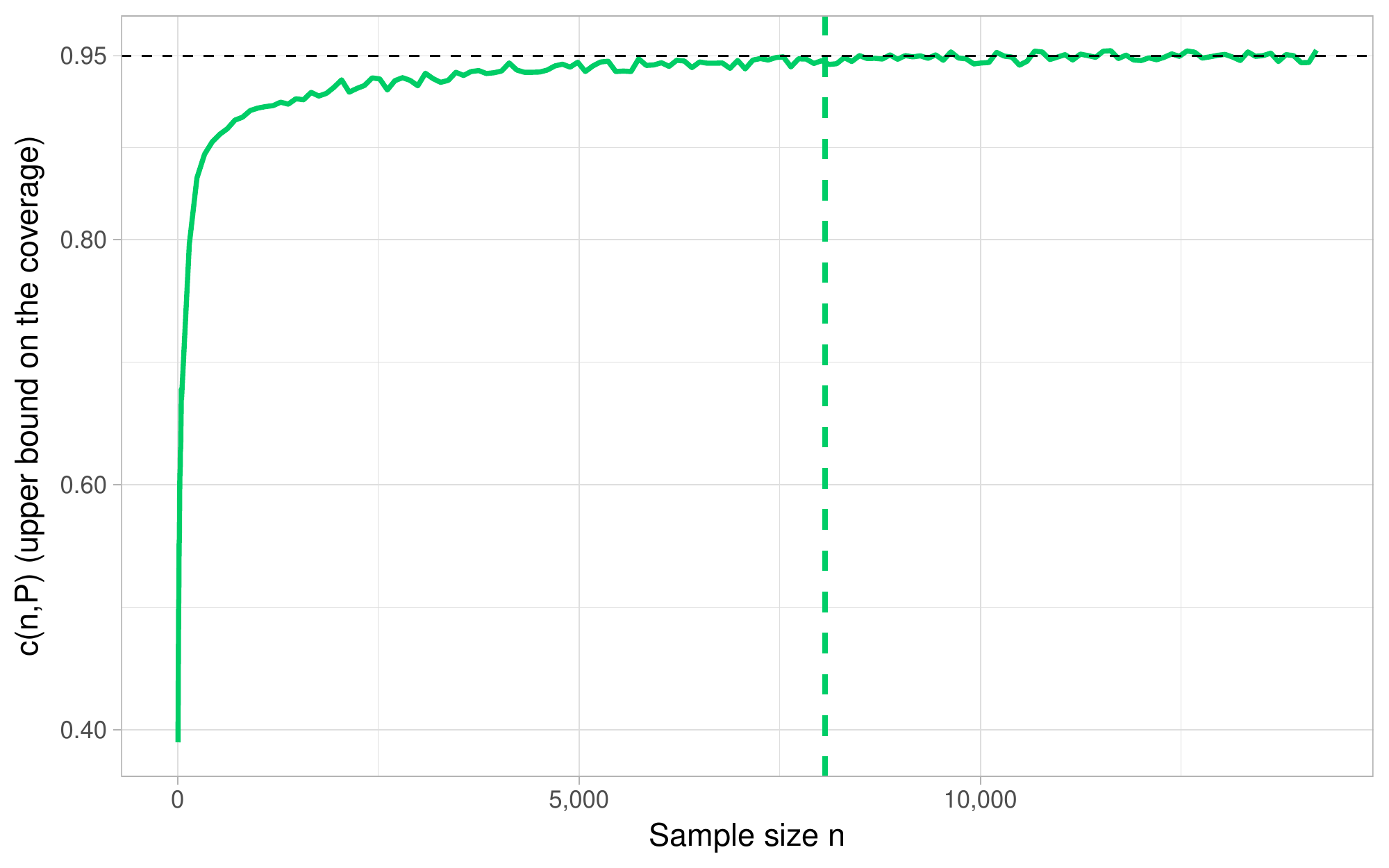}
    \caption{\small{Specification: $\forall n \in \Nstar, {\distributionXYn = \mathcal{P}areto(1, -1.5, 5) \otimes \mathcal{P}areto(\expecY, -1.5, 5)}$, with $\expecY = 0.1$;
    $1-\alpha=0.95$;
    5,000 repetitions used.}}   
    \label{fig:appendix:gr_newa_paretos_n_alphapc_5nbrep_5000EYpc_10thrsYpc_-100shpY_5EXpc_100thrsXpc_-150shpX_5CorrXYpc_0}
\end{figure}

\subsection{Bernoulli distributions}

\begin{figure}[ht]
\centering
\includegraphics[width=0.7\textwidth]{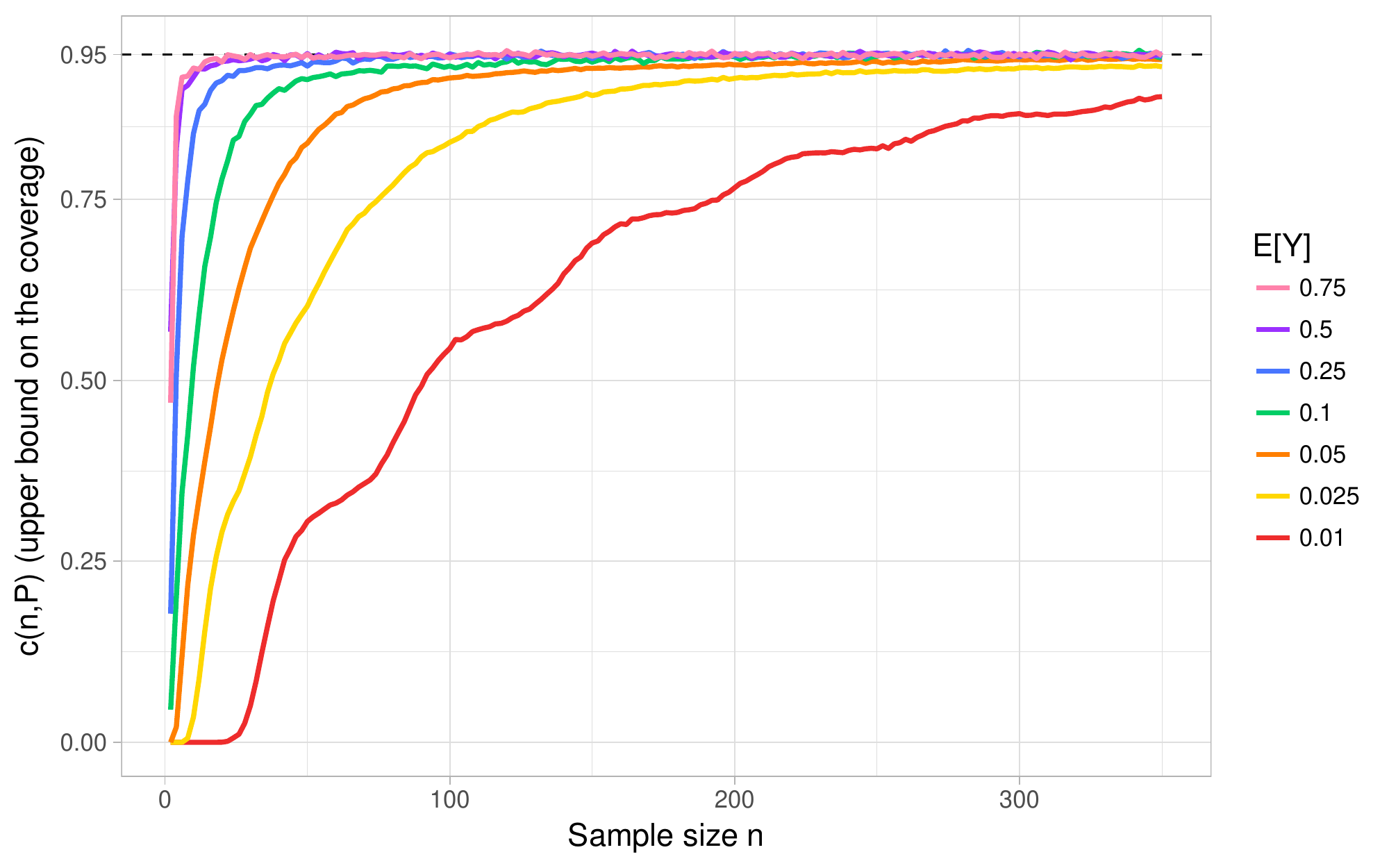}
\caption{\small{$c(n,P)$ of asymptotic CIs based on the delta method as a function of the sample size~$n$. \newline Specification: $\forall n \in \Nstar, {\distributionXYn = \mathcal{B}(0.5) \otimes \mathcal{B}(\expecY)}$.
The nominal pointwise asymptotic level is set to $0.95$.
For each pair $(\expecY, n)$, the coverage is obtained as the conditional mean over repetitions such that ${\meanY \neq 0}$ (starting from 10,000 repetitions).}}
\label{appendix:fig:slow_convergence_bernoulli_case_EX05_CorrXY0}
\end{figure}

Figure~\ref{appendix:fig:slow_convergence_bernoulli_case_EX05_CorrXY0} is the equivalent of Figure~\ref{fig:slow_convergence_normal_case_EX1_VX1_VY1_CorrXY0} with Bernoulli distributions.
The following graphs illustrate the use of $\overline{n}_\alpha$ to appraise the reliability of the asymptotic confidence based on the delta method.
In practice a plug-in strategy has to be used to compute $\overline{n}_\alpha$ and, in the setting of simulations, we simply use the known moments and bounds of the DGP used in the simulation.
With two Bernoulli variables in the numerator and the denominator, we are both in the BC and the ``Hoeffding'' cases.
Thus, we show both the one obtained in the BC case $\overline{n}_\alpha :=2\left(\upperYn-{\lowerYn}^2\right)/\left(\alpha\lowerYn^2\right)$ with a dashed vertical line (Theorem~\ref{thm:concentration_inequality_BC_general_case})  and the one obtained in the ``Hoeffding'' case $\overline{n}_{\alpha,H} := \ln(4/\alpha) / \gamma(Y_{1,n})$, setting here $\aYn = 0$, $\bYn = 1$ and $\lowerYn = \expecY$, with a dotted vertical line (Theorem~\ref{thm:concentration_inequality_H_general_case}).

\begin{figure}[ht]
    \centering
    \includegraphics[width=0.6\textwidth]{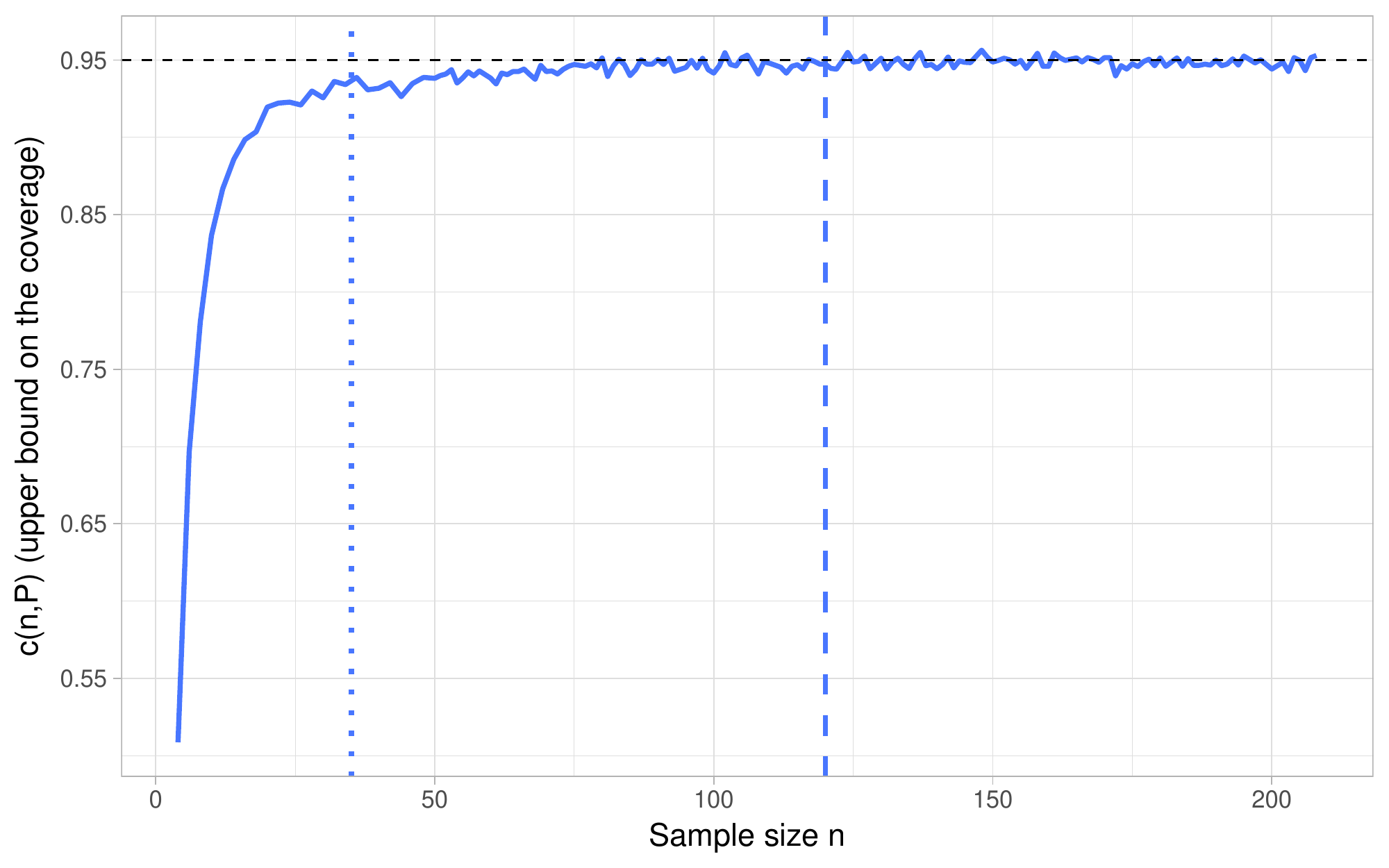}
    \caption{\small{Specification: $\forall n \in \Nstar, {\distributionXYn = \mathcal{B}(0.5) \otimes \mathcal{B}(0.25)}$;
    $1-\alpha=0.95$;
    5,000 repetitions used.}}   
    \label{fig:appendix:gr_newlab_uncondcov_bernoulli_n_alphapc_5nbrep_5000EYpc_25EXpc_50CorrXYpc_0}
\end{figure}

\begin{figure}[ht]
    \centering
    \includegraphics[width=0.6\textwidth]{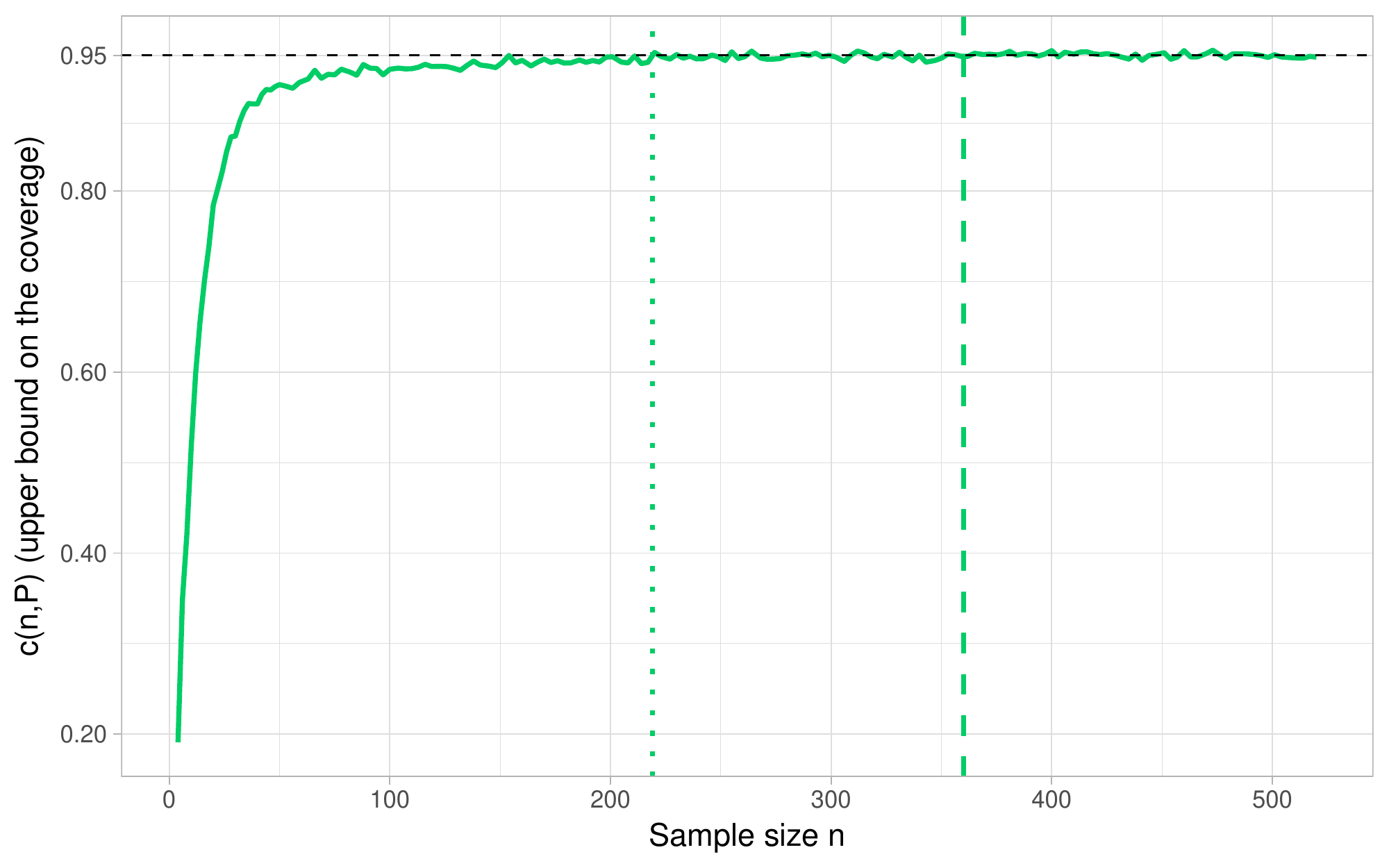}
    \caption{\small{Specification: $\forall n \in \Nstar, {\distributionXYn = \mathcal{B}(0.5) \otimes \mathcal{B}(0.1)}$;
    $1-\alpha=0.95$;
    5,000 repetitions used.}}   
    \label{fig:appendix:gr_newlab_uncondcov_bernoulli_n_alphapc_5nbrep_5000EYpc_10EXpc_50CorrXYpc_0}
\end{figure}

\subsection{Poisson distributions}

The specification here considers two variables distributed according to a Poisson distribution, both in the numerator and in the denominator.

\begin{figure}[ht]
    \centering
    \includegraphics[width=0.6\textwidth]{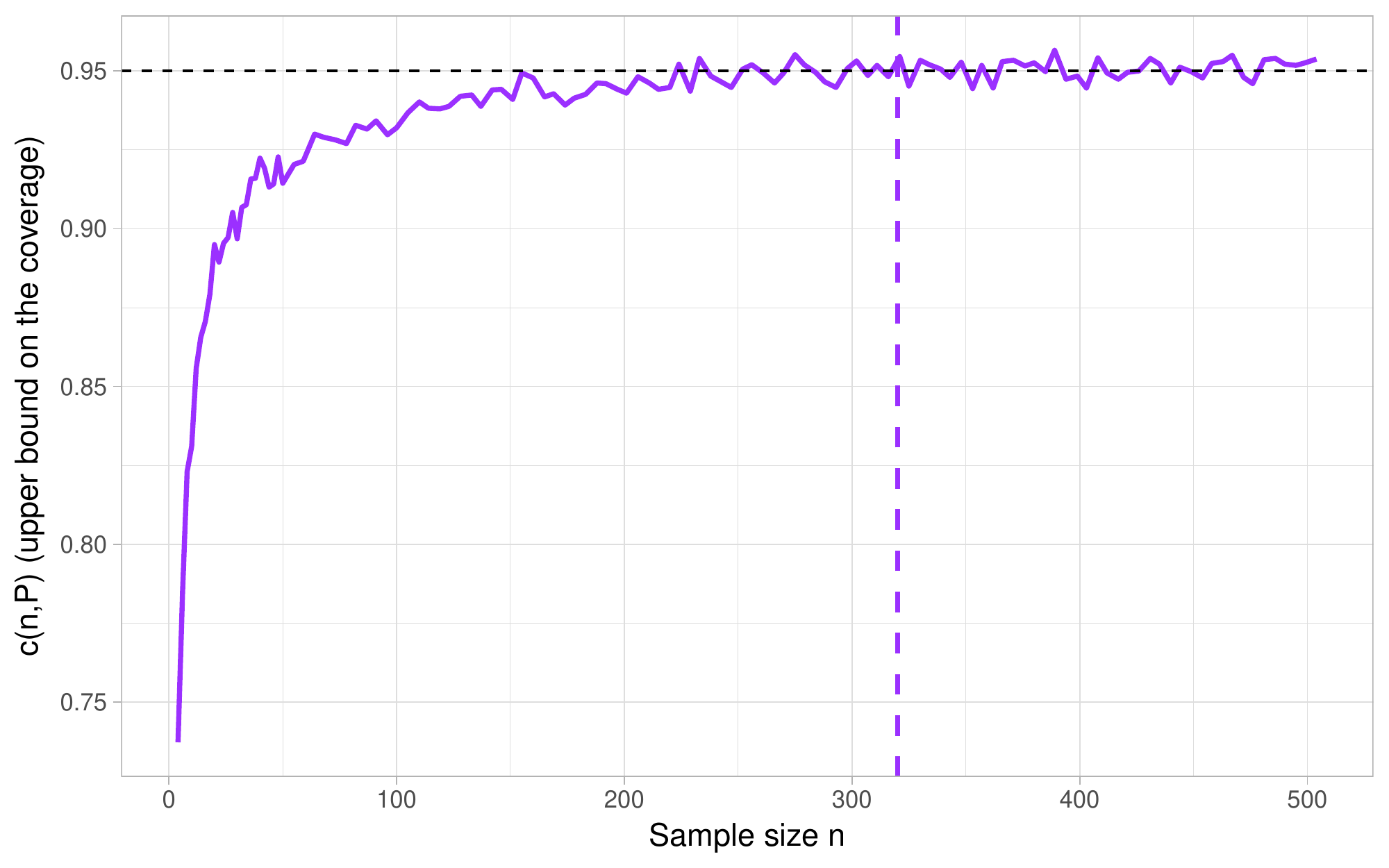}
    \caption{\small{Specification: $\forall n \in \Nstar, {\distributionXYn = \mathcal{P}oisson(0.5, 2) \otimes \mathcal{P}oisson(0.5, 2)}$;
    $1-\alpha=0.95$;
    5,000 repetitions used.}}   
    \label{fig:appendix:gr_newaylab_poissons_n_alphapc_5nbrep_5000EYpc_50VYpc_200EXpc_50VXpc_200CorrXYpc_0dilafornpc_140}
\end{figure}

A Poisson distribution is entirely defined by its positive real parameter, which is equal to both its expectation and its variance.
Consequently, to have denominator close to $0$, we would need small variance too, as in the exponential specification (see Section \ref{appendix:sec:exponentials}).
In order to disentangle expectation and variance, we use below translated Poisson variables.
More precisely, the notation ${\mathcal{P}oisson(\mu, \sigma^2)}$, $\mu \in \Rb$, $\sigma^2 \in \Rstarplus$, denotes a distribution alike to a Poisson, with parameter and variance equal to $\sigma^2$ but translated such that its expectation is $\mu$.
That is a variable distributed according to ${\mathcal{P}oisson(\mu, \sigma^2)}$ is equal in distribution to $P + (\mu - \sigma^2)$ with $P$ a standard Poisson distribution with parameter $\sigma^2$ - that is with density with respect to the counting measure equal to $(\sigma^2)^k \exp(-\sigma^2)/(k!)$ for every ${k \in \mathbb{N}}$.
Thus, a ${\mathcal{P}oisson(\mu, \sigma^2)}$ has expectation $\mu$ and variance $\sigma^2$.

\begin{figure}[ht]
    \centering
    \includegraphics[width=0.6\textwidth]{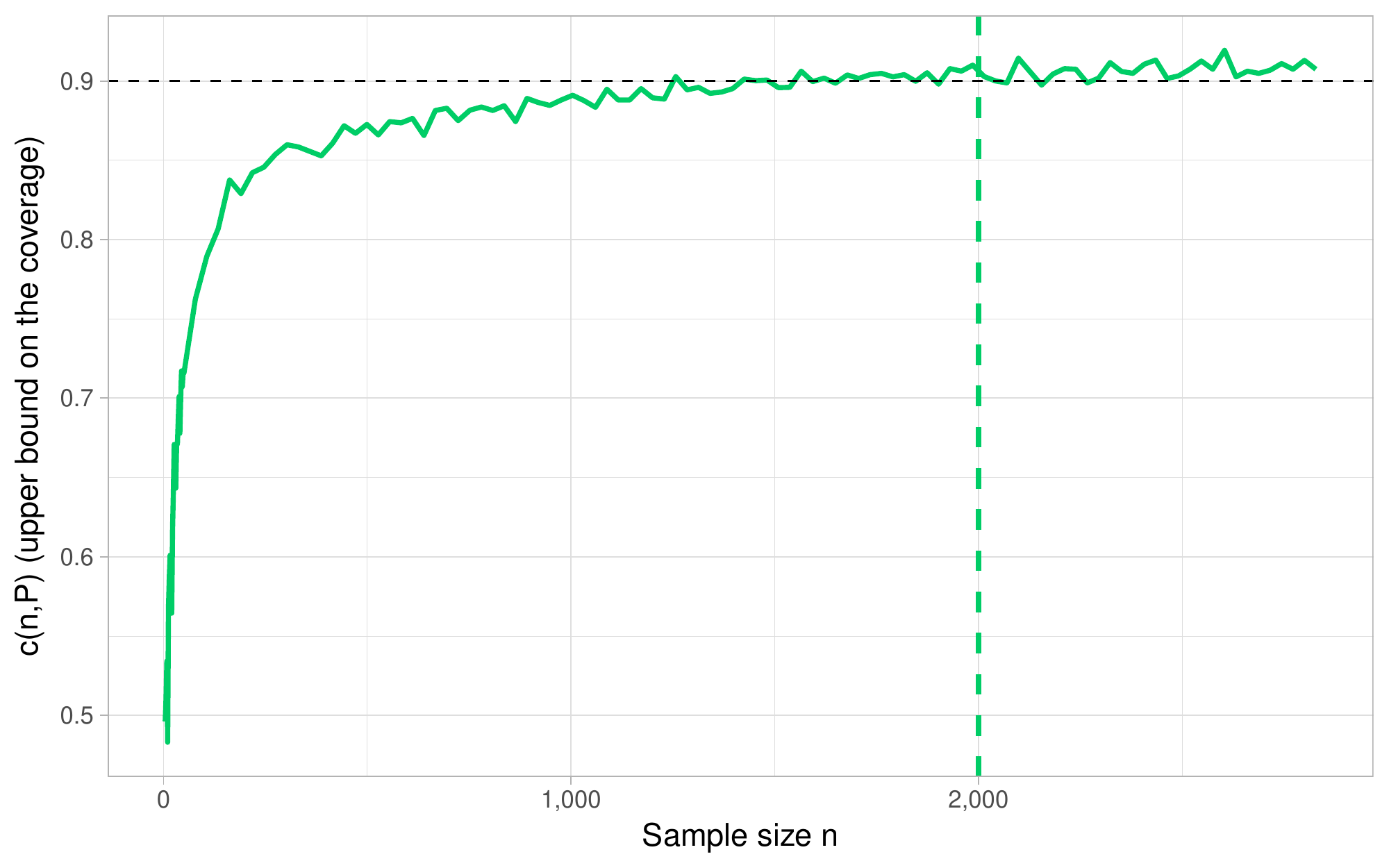}
    \caption{\small{Specification: $\forall n \in \Nstar$, the marginal distributions of $X$ and $Y$ are respectively ${\mathcal{P}oisson(0.5, 2)}$ and ${\mathcal{P}oisson(0.1, 1)}$ and are simulated using a Gaussian copula to have ${\Corr(X,Y) \approx 0.6}$;
    $1-\alpha=0.9$;
    5,000 repetitions used.}}   
    \label{fig:appendix:gr_newaylab_poissons_n_alphapc_10nbrep_5000EYpc_10VYpc_100EXpc_50VXpc_200CorrXYpc_70dilafornpc_140}
\end{figure}

\subsection{Delta method and nonparametric percentile bootstrap confidence intervals}
\label{appendix:sub:additional_simulations_DM_bootstrap}

\begin{figure}[ht]
    \caption{\small{delta method in blue; Efron's percentile bootstrap in red; $C=0.2$.}}
    \label{fig:bootstrap_DM_sequences_of_models_C_02} 
    \begin{minipage}[b]{0.5\linewidth}
    \includegraphics[width=1\linewidth]{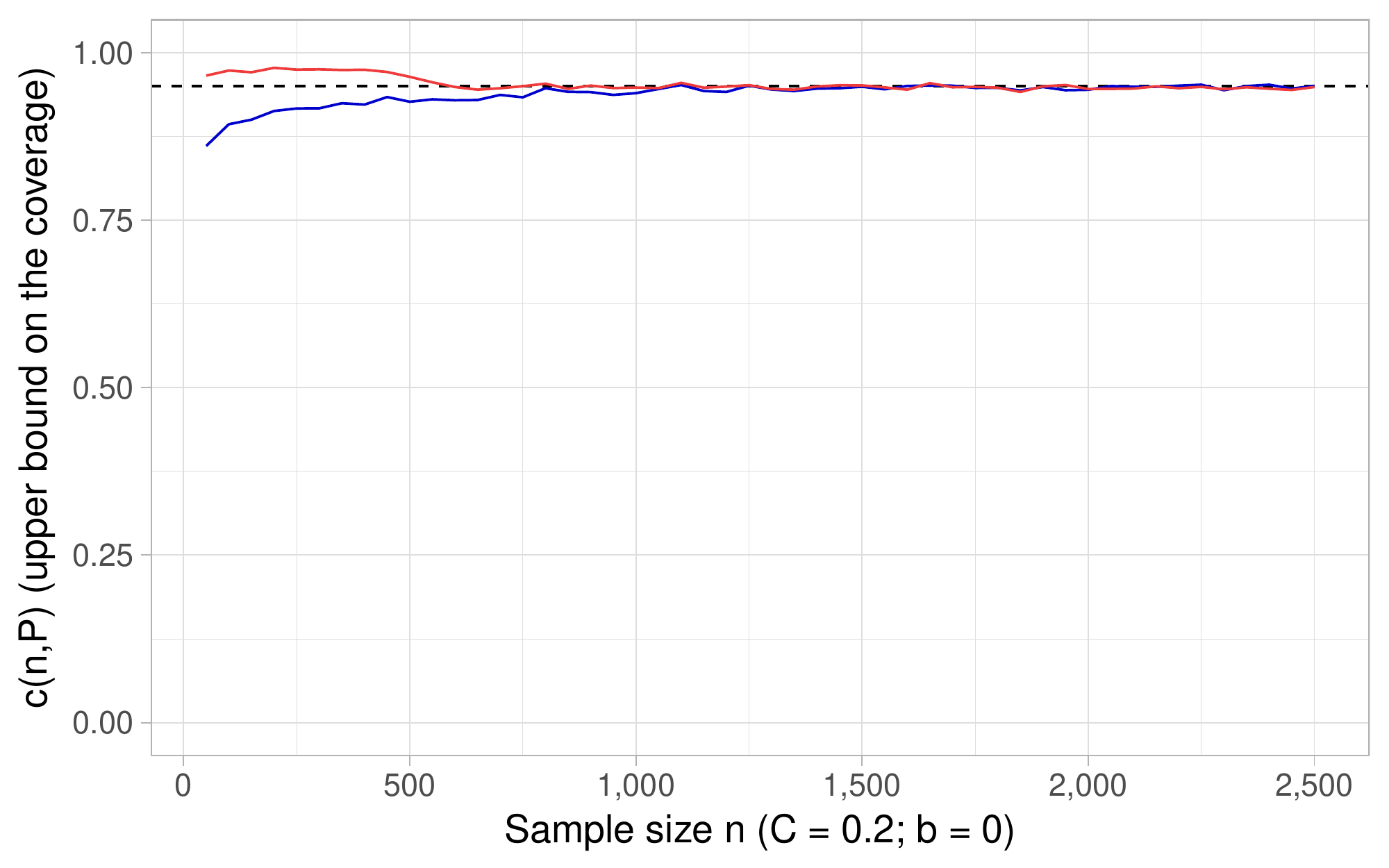}
    \end{minipage} 
    \begin{minipage}[b]{0.5\linewidth}
    \includegraphics[width=1\linewidth]{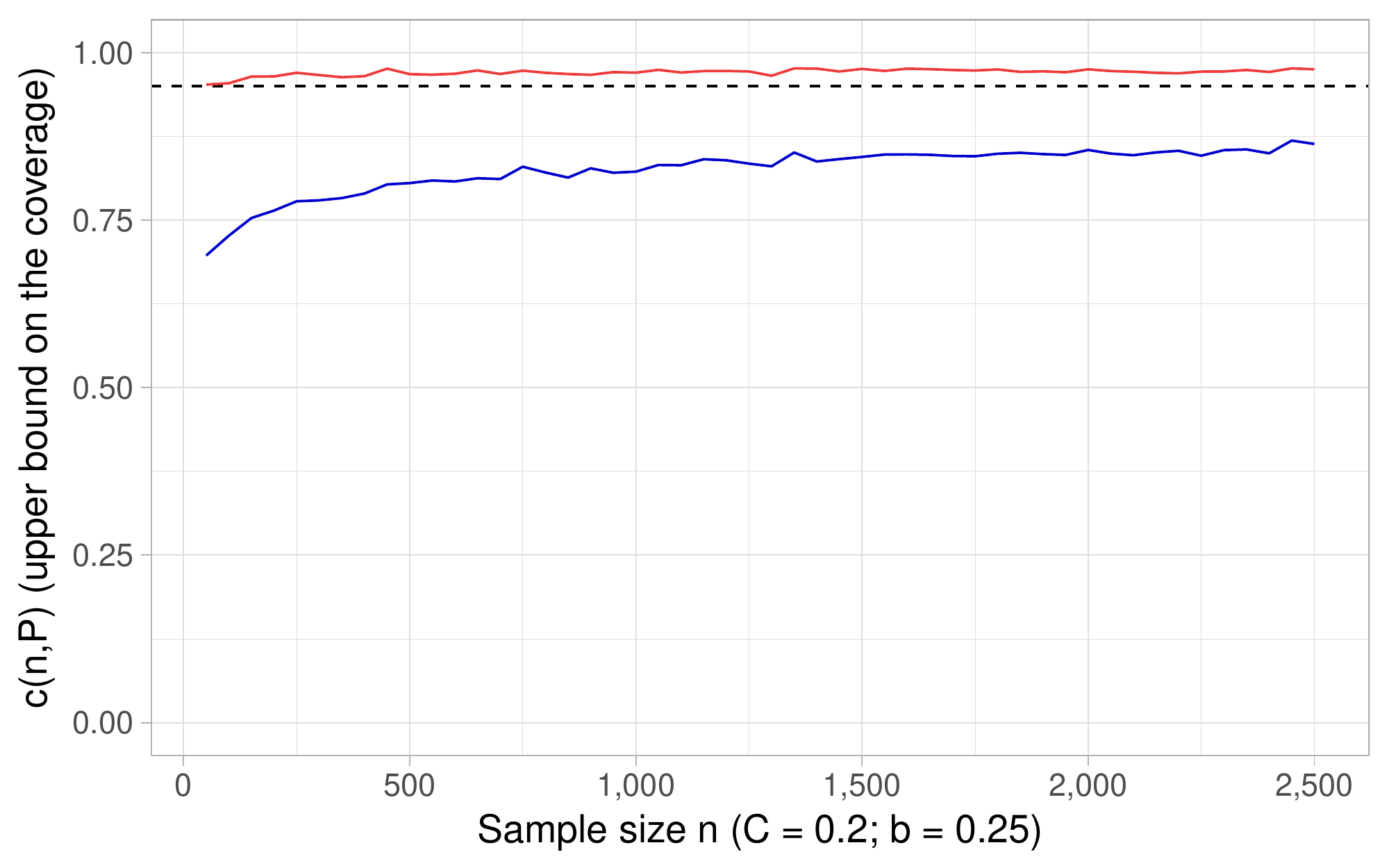}
    \end{minipage} 
    \begin{minipage}[b]{0.5\linewidth}
    \includegraphics[width=1\linewidth]{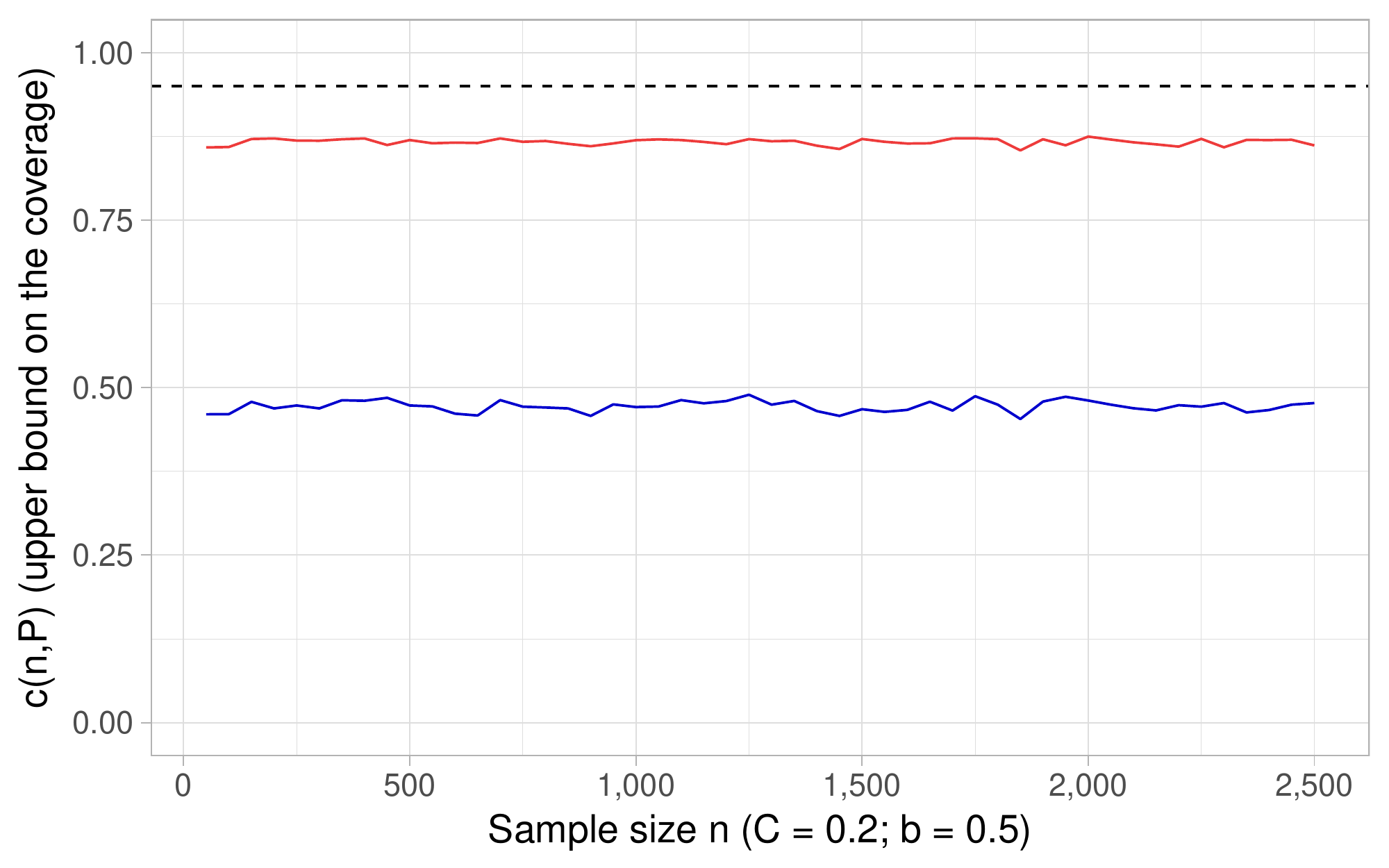}
    \end{minipage}
    \hfill
    \begin{minipage}[b]{0.5\linewidth}
    \includegraphics[width=1\linewidth]{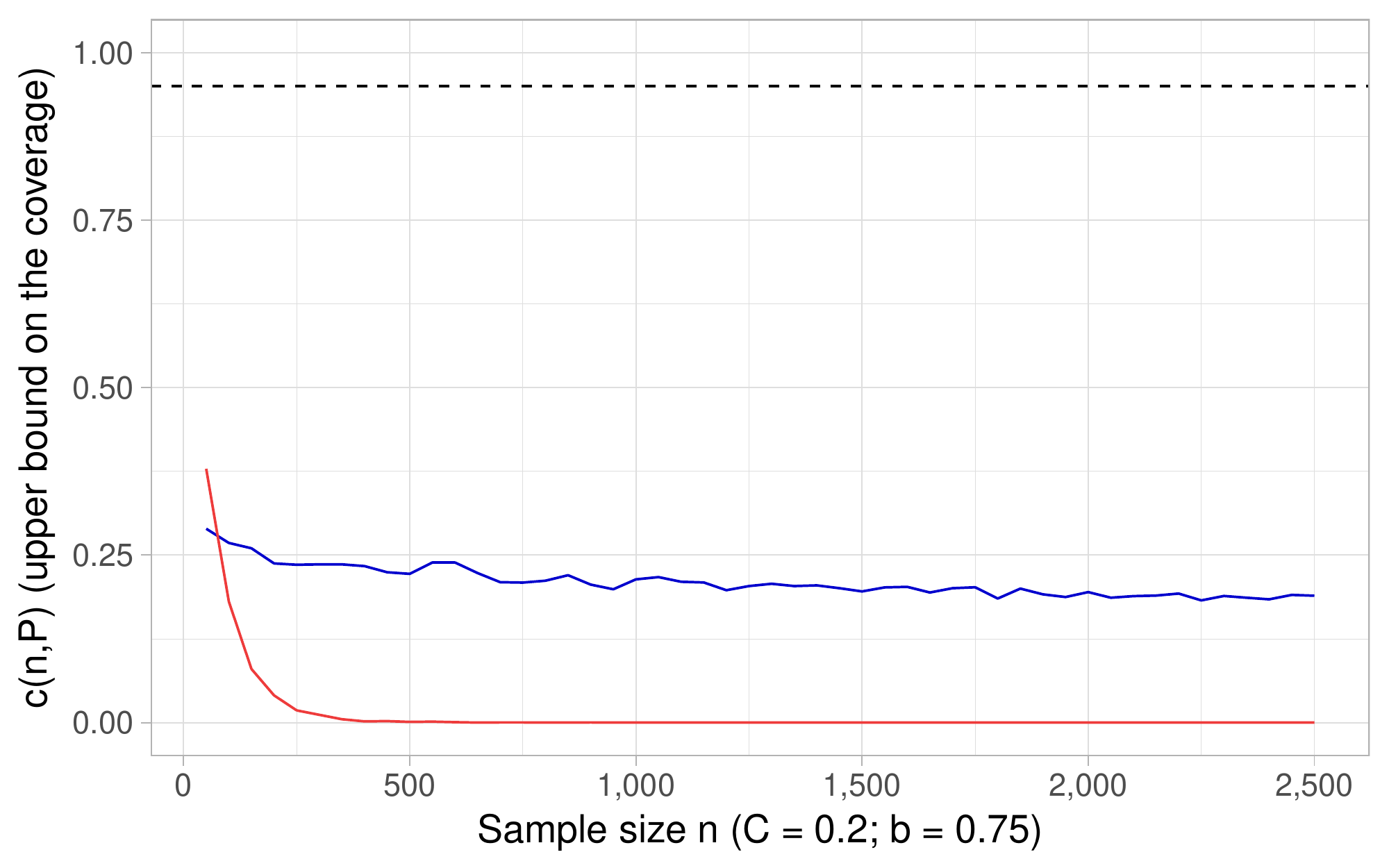}
    \end{minipage} 
\end{figure}

\begin{figure}[ht]
    \begin{minipage}[b]{0.5\linewidth}
    \includegraphics[width=1\linewidth]{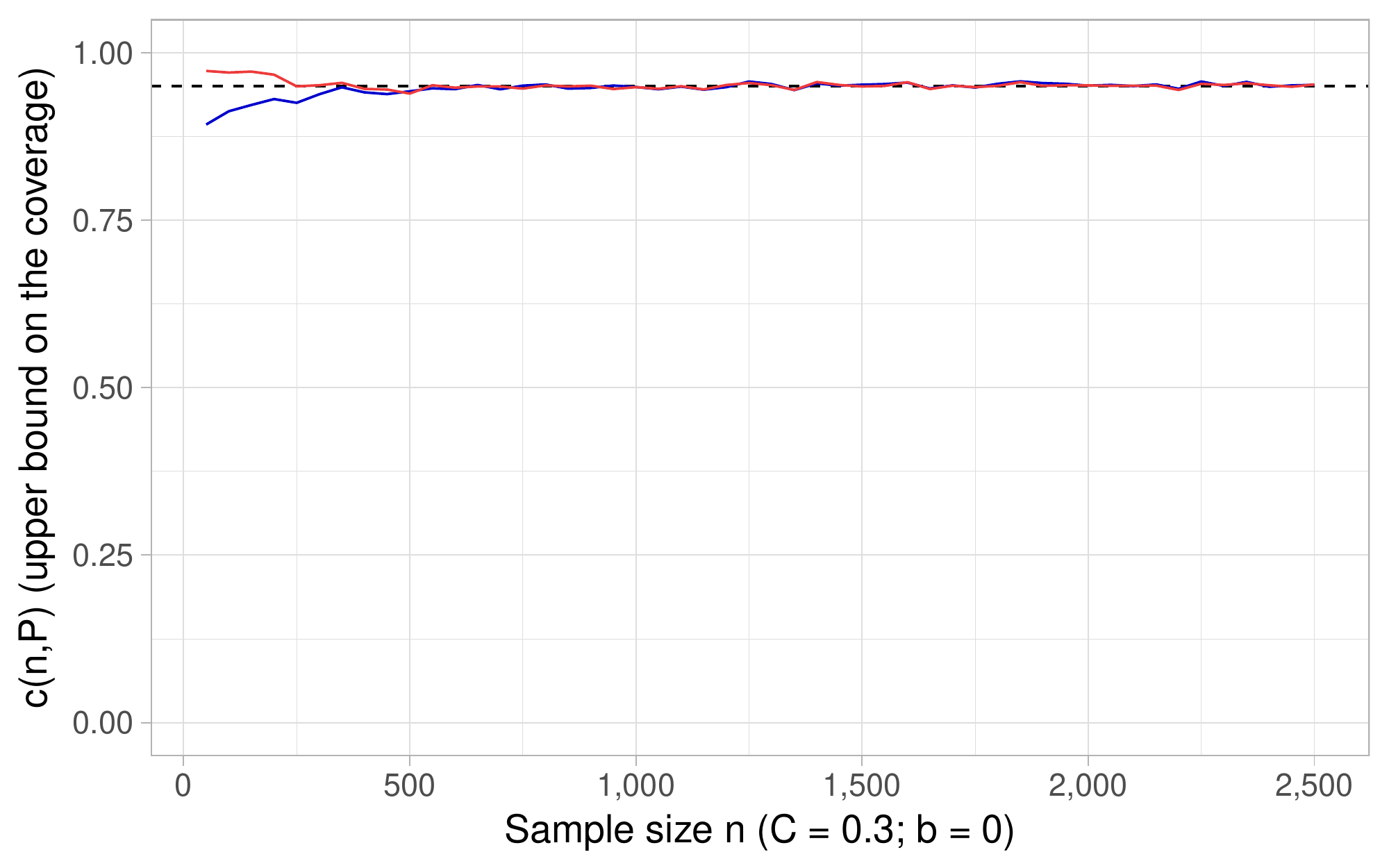}
    \end{minipage} 
    \begin{minipage}[b]{0.5\linewidth}
    \includegraphics[width=1\linewidth]{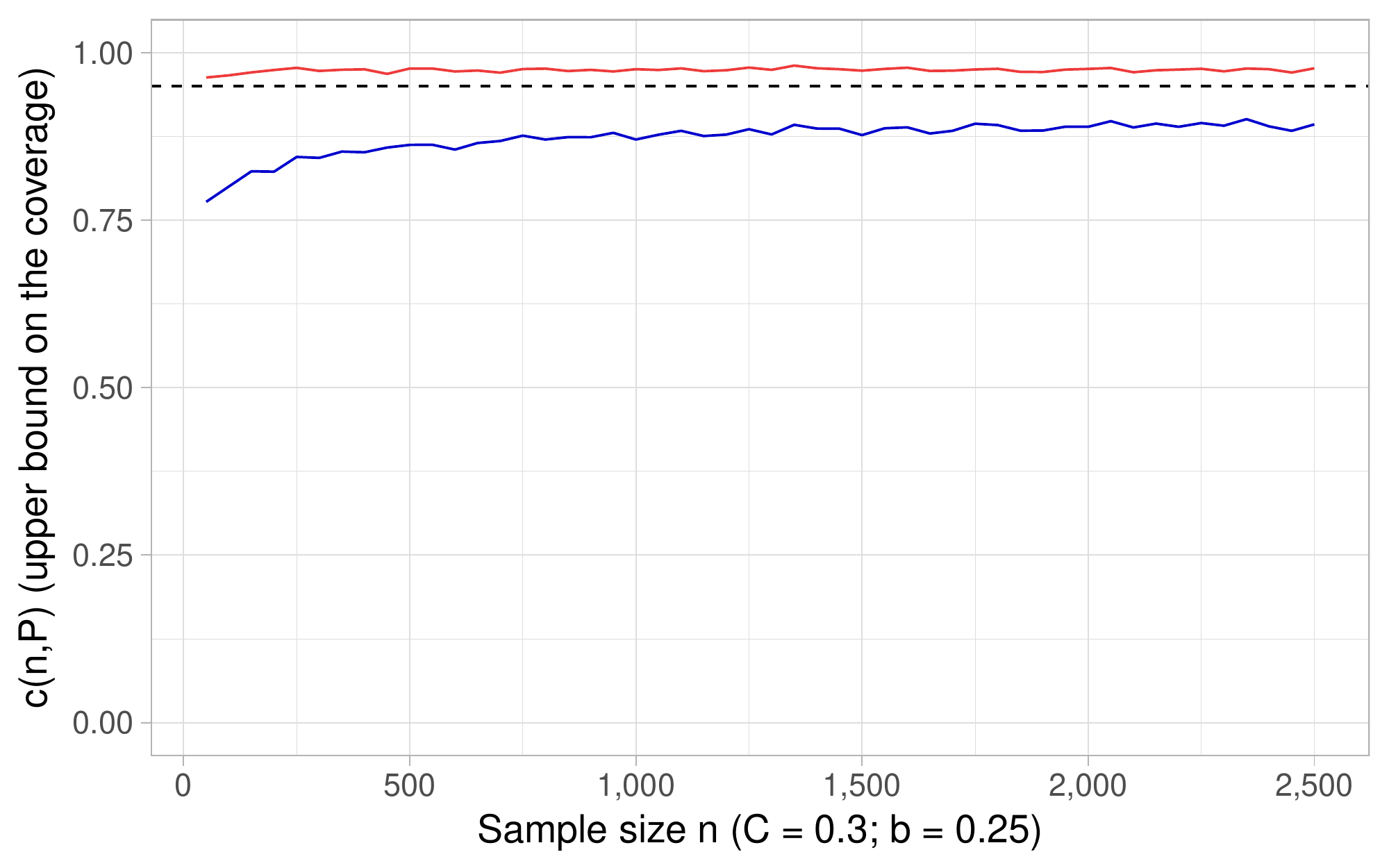}
    \end{minipage} 
    \begin{minipage}[b]{0.5\linewidth}
    \includegraphics[width=1\linewidth]{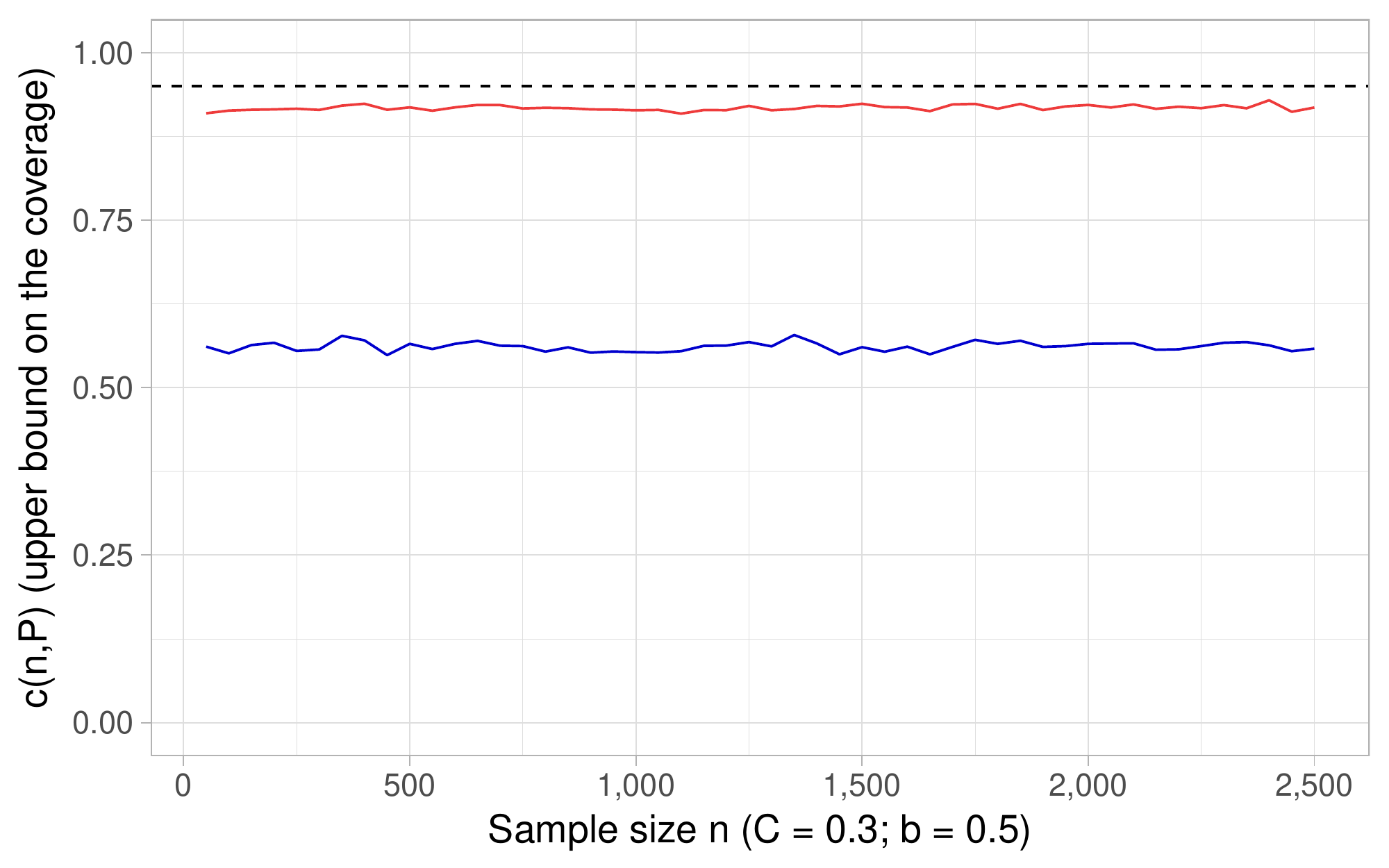}
    \end{minipage}
    \hfill
    \begin{minipage}[b]{0.5\linewidth}
    \includegraphics[width=1\linewidth]{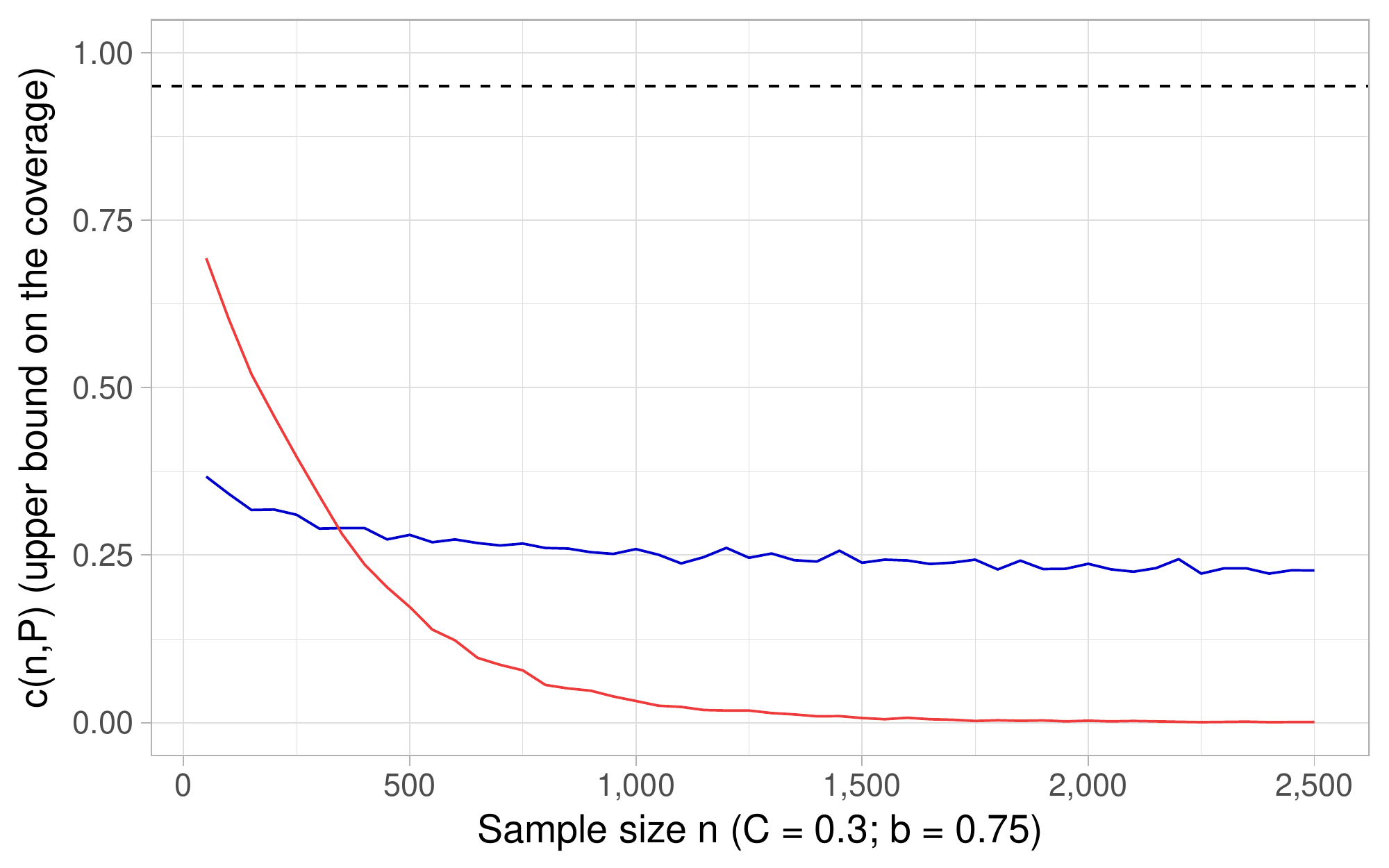}
    \end{minipage} 
    \caption{\small{delta method in blue; Efron's percentile bootstrap in red; $C=0.3$.}}
    \label{fig:bootstrap_DM_sequences_of_models_C03} 
\end{figure}

The two following figures are the equivalent to Figure~\ref{fig:bootstrap_DM_sequences_of_models_C01} with different values of~$C$.
They illustrate that the lower~$C$, the lower the signal-to-noise ratio in the denominator, hence the more difficult in some sense is the estimation of~$\theta_n$.
This is illustrated by the fact that, all other things equal, larger~$C$ basically translates $c(n,P)$ upward as revealed by the series of Figures~\ref{fig:bootstrap_DM_sequences_of_models_C01}, \ref{fig:bootstrap_DM_sequences_of_models_C_02}, and~\ref{fig:bootstrap_DM_sequences_of_models_C03}.

These three figures all report the $c(n,P)$ of the CIs based on the delta method (in blue) and of the CIs constructed with Efron's non parametric bootstrap using 2,000 bootstrap replications (in red) with the specification $\forall n \in \Nstar, \distributionXYn = \mathcal{N}(1,1) \otimes \mathcal{N}(Cn^{-b},1)$, with ${b \in \{0, 0.25, 0.5, 0.75\}}$.
For the three of them, the nominal pointwise asymptotic level is set to $0.95$ and for each pair~${(b, n)}$, the coverage is obtained as the mean over 5,000 repetitions.

\vspace{-1em}

\end{document}